\documentclass[a4paper,12pt]{amsart}

\usepackage[
  margin=30mm,
  marginparwidth=25mm,     
  marginparsep=2mm,       
  bottom=25mm,
  ]{geometry}

\usepackage[T1]{fontenc}
\usepackage[utf8]{inputenc}
\usepackage[english]{babel}
\usepackage{amsmath,amssymb,amsthm,mathtools}

\expandafter\let\expandafter\xbf\csname bfseries \endcsname
\expandafter\let\expandafter\xmd\csname mdseries \endcsname
\let\xbar\bar
\usepackage{allrunes}
\let\bar\xbar
\expandafter\let\csname bfseries \endcsname\xbf
\expandafter\let\csname mdseries \endcsname\xmd

\usepackage{latexsym}
\usepackage{delarray}
\usepackage{bbm}
\usepackage{scrextend}
\usepackage{hyperref}
\usepackage[pdftex,usenames,dvipsnames]{xcolor}
\usepackage{datetime}
\usepackage{mathrsfs}
\usepackage{enumitem}  
\usepackage{tikz}
\usetikzlibrary{fadings}

\newtheorem{Th}{Theorem}[section]
\newtheorem{Prop}[Th]{Proposition}
\newtheorem{Lem}[Th]{Lemma}

\newtheorem{Rem}[Th]{Remark}
\newtheorem{Def}[Th]{Definition}

\newcommand{\eps}{\varepsilon}

\newcommand{\xxi}{\langle \xi\rangle}
\newcommand{\Rl}{\mathbb{R}}

\newcommand{\Rn}{\mathbb{R}^{n}}
\newcommand{\Z}{\mathbb{Z}}

\renewcommand{\d}{\partial}

\def\bra#1{{\langle{#1}\rangle}}

\newcommand{\FF}{\mathcal{F}}

\newcommand{\jap}[1]{\left\langle #1\right\rangle}

\newcommand{\phase}{\varphi}

\newcommand{\dd}{\, \mathrm{d}}
\newcommand{\ddd}{\,\text{\rm{\mbox{\dj}}}}

\renewcommand{\SS}{{\mathscr{S}}}

\newcommand{\at}{\mathfrak{a}}

\newcommand{\brkt}[1]{\Big({#1}\Big)}
\newcommand{\set}[1]{\left\{{#1}\right\}}
\newcommand{\norm}[1]{\Big\|#1\Big\|}
\newcommand{\abs}[1]{\Big |#1\Big |}
\newcommand{\diffcases}[1]{\begin{cases}#1\end{cases}}

\newcommand{\eq}[1]{
    \begin{align*}
        #1
    \end{align*}
}

\newcommand{\nm}[2]{\begin{align}\label{#1}
#2
\end{align}}

\DeclareMathOperator{\supp}{supp}

\allowdisplaybreaks 

\title[Boundedness of FIOs]
{Boundedness of Fourier integral operators on classical function spaces}
\dedicatory{\emph{To the memory of our friend Gustav Hammarhjelm}}
\author[A. Israelsson]{Anders Israelsson}
\author[T. Mattsson]{Tobias Mattsson}
\author[W. Staubach]{Wolfgang Staubach}

\address{\newline
      Anders Israelsson, Tobias Mattsson,
       Wolfgang Staubach \newline
       Department of  Mathematics, Uppsala University, \newline
       S-751 06 Uppsala, Sweden}
       \email{anders.israelsson@math.uu.se, tobias.mattsson@math.uu.se, wulf@math.uu.se}

 \thanks{
 The second author is supported by the Knut and Alice Wallenberg Foundation.}



\begin{document}


\maketitle

\begin{abstract}
We investigate the global boundedness of Fourier integral operators with amplitudes in the general H\"ormander classes $S^{m}_{\rho, \delta}(\Rl^n)$, $\rho, \delta\in [0,1]$ and non-degenerate phase functions of arbitrary rank $\kappa\in \{0,1,\dots, n-1\}$ on Besov-Lipschitz $B^{s}_{p,q}(\Rl^n)$ and Triebel-Lizorkin $F^{s}_{p,q}(\Rl^n)$ of order $s$ and $0<p\leq\infty$, $0<q\leq\infty$. The results that are obtained are all up to the end-point and sharp and are also applied to the regularity of Klein-Gordon-type oscillatory integrals in the aforementioned function spaces.
\end{abstract}
\vspace*{0.5cm}

\tableofcontents

\section{Introduction}
In this paper, we investigate the global regularity of oscillatory integral operators of the form
\begin{equation*}
	T_a^\varphi f(x) = \frac{1}{(2\pi)^n} \int_{\Rl^n} e^{i\varphi(x,\xi)}\,a(x,\xi)\,\widehat f (\xi) \dd\xi,
\end{equation*}
with amplitudes in the general H\"ormander class $S^m_{\rho, \delta}(\Rl^n)$ (Definition \ref{symbol class Sm}), 
on Besov-Lipschitz $B^{s}_{p,q}(\Rl^n)$ and Triebel-Lizorkin $F^{s}_{p,q}(\Rl^n)$ of order $s$ with $0<p\leq\infty$ and $0\leq q<\infty.$ Throughout the paper, we are assuming that the phase function $\varphi(x, \xi)$ is positively homogeneous of degree one in $\xi$, is strongly non-degenerate (Definition \ref{nondeg phase}) and is of rank $\kappa$ (Definition \ref{def:FIO}). The operator $T^\varphi_{a}$ is then referred to as a {\it Fourier integral operator} by L. H\"ormander \cite{Hor:acta}, for which in the case of $\rho\in(\frac{1}{2}, 1]$ and $\delta=1-\rho$, he developed a local theory and a calculus.  Note that the rank $\kappa$ of the phase here is an integer between $0$ and $n-1$. The case of $\kappa=0$ corresponds to the case of {\it pseudodifferential operators}, while the case of $\kappa=n-1$ is the relevant one for the wave operator $e^{it\sqrt{-\Delta}}$ for fixed time $t,$ which is a typical example of a Fourier integral operator.\\

Regarding the regularity of Fourier integral operators in Besov-Lipschitz and Triebel-Lizorkin spaces, in \cite{Brenner} P. Brenner showed that under suitable regularity hypotheses on the functions \(f_{0}\) and \(f_{1}\), if 
$$f_{+}:=\frac{1}{2}\big(f_{0}-i(\sqrt{-\Delta})^{-1} f_{1}\big)$$
and 
$$f_{-}:=\frac{1}{2}\big(f_{0}+i(\sqrt{-\Delta})^{-1} f_{1}\big)$$
then for a fixed time $t>0$ one has the estimate
\begin{equation}\label{brenners global besov estimate for the wave equation}
\big\Vert  e^{i t \sqrt{-\Delta}} f_{+}+e^{-i t \sqrt{-\Delta}} f_{-}\big\Vert_{B^{s}_{p,q}(\Rl^n)}\leq C_t \left ( \Vert f_0\Vert_{B^{s+m}_{p',q}(\Rl^n)}+ \Vert f_1\Vert_{B^{s+m-1}_{p',q}(\Rl^n)} \right ),
\end{equation}
where $s\in \Rl$, $p\in [2,\infty)$, $ p'=\frac{p}{p-1}$, $q\in [1,\infty]$, and $\displaystyle (n+1)\abs{\frac{1}{p}-\frac{1}{2}}\leq m\leq 2n\abs{\frac{1}{p}-\frac{1}{2}}$. 

\noindent In \cite {Kapitanski} L. V. Kapitanski\u\i, extended and improved the results of Brenner to the range $p\in [2,\infty]$ and $\displaystyle(n-1)\abs{\frac{1}{p}-\frac{1}{2}}\leq m\leq n\abs{1-\frac{2}{p}}.$\\

\noindent The pioneering results of Brenner's and Kapitanski\u\i's, and the later results of Q.J. Qiu \cite{Qiu} on Besov-Lipschitz spaces, were both concerned with estimates for Fourier integral operators with rank equal to $n-1$, where the amplitudes are in the H\"ormander classes $S^m_{1,0}(\Rn)$ for a suitable $m$ and also with the Banach-space scales of the aforementioned function spaces.\\

For operators with amplitudes in other H\"ormander classes, H. Triebel \cite{TriebelFIO} used the Frazier-Jawerth atomic decomposition of Triebel-Lizorkin spaces, to investigate the regularity of Fourier integral operators on these spaces. In detail, let \(m \in \Rl, M \in \mathbb{N}\) and \(N \in \mathbb{N}\). Then one says that \(b(x, \xi)\) belongs to the symbol
class \(S_{1}^{m\, \mathrm{loc}}(M, N)\)  if for any \(R>0\) there exists a constant \(c_{R}>0\) with

\begin{equation}
\sup _{|x| \leqq R}\left|\partial^{\alpha}_{\xi}\partial^{\beta}_x b(x, \xi)\right| \leqq c_{R}|\xi|^{m-|\alpha|+|\beta|}, \quad|\xi| \geq 1,|\alpha| \leq M,|\beta| \leq N
\end{equation}
Triebel \cite{TriebelFIO} showed the following:  Let \(0<p \leq 1<q \leq \infty\) and \(s>n\left(\frac{1}{p}-1\right)\). Let \(N \in \mathbb{N}\) with \(N>1+\frac{n}{p}\) and
 \(m \geq 0,\, m' \geq 0\) and
 \begin{equation}
(\varphi(x, \xi)-x\cdot\xi) \in S_{1}^{m' \, \mathrm{loc}}(1+[s], N), \quad a(x, \xi) \in S_{1}^{-m \, \mathrm{loc}}(1+[s], N)
\end{equation}
with \(m \geq m'(1+[s]+N)\). Assume also that
\begin{equation}\sup _{|x| \leq R}\left|\partial_{x}^{\alpha}\left(e^{i (\varphi(x, \xi)-x\cdot \xi)}\, a(x, \xi)\right)\right| \in L_{\mathrm{loc}}^{1}(\Rl^n_{\xi})
\end{equation}
for all \(|\alpha| \leq 1+[s]\) and all \(R>0\). Then the Fourier integral operator $T_a^{\varphi}$ maps $F^{s\, \mathrm{comp}}_{p,q}(\Rl^n) \to F^{s\, \mathrm{loc}}_{p,q}(\Rl^n).$ This result of Triebel's is the first known result regarding the regularity of Fourier integral operators with amplitudes in the so-called forbidden H\"ormander class $S^{m}_{1,1}(\Rl^n),$ in Triebel-Lizorkin spaces.\\

\noindent From the point of view of our work in this paper, a breakthrough regarding investigations of the regularity of integral operators on Triebel-Lizorkin (and therefore also Besov-Lipschitz) spaces was made in the paper by M. Pramanik, K. Rogers and A. Seeger \cite{PRS}. Therein, the authors proved a Calder\'on-Zygmund-type estimate with quite a few applications, including the regularity of Radon transforms and Fourier integral operators, just to name a couple. Specifically, the authors considered local Fourier integral operators $T_a^{\varphi}$ where $a(x,\xi)\in S^{m}_{1,0}(\Rn)$ is compactly supported in $x$ and $\varphi(x, \xi)$ is non-degenerate on the support of $a(x,\xi).$ Using their Calder\'on-Zygmund estimate in \cite{PRS}, they showed that if \(n \geq 2,\,2<p<\infty,\, q>0\), then  \(T_a^{\varphi}: F^{0}_{p,p}(\mathbb{R}^{n}) \rightarrow F^{0}_{p,q}(\mathbb{R}^{n})\), provided that $m=-(n-1)\big|\frac{1}{p}-\frac{1}{2}\big|.$\\

\begin{figure}[ht!]
\centering\includegraphics[scale=.75]{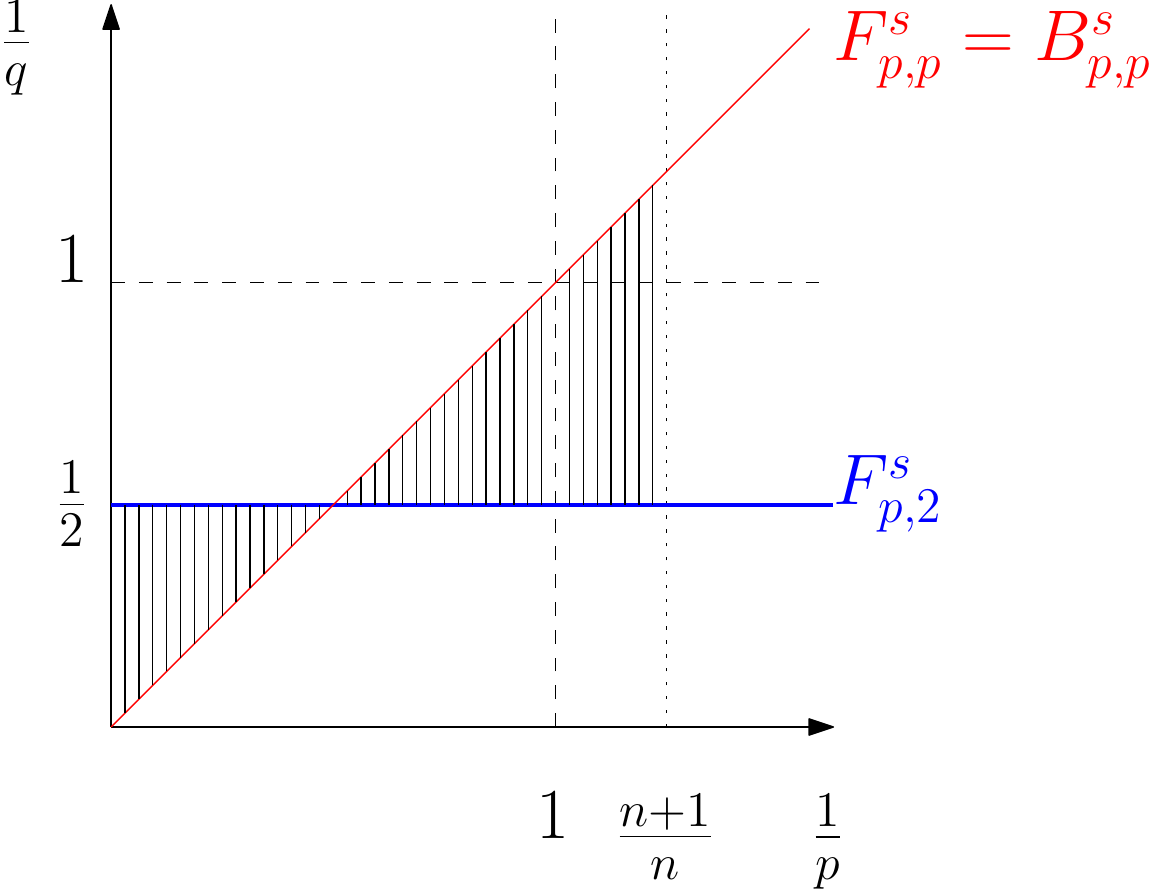}
\caption{Global boundedness in Triebel-Lizorkin scale.}\label{pic:gammalskåpmat}
\end{figure}

In recent years, the second and the third author, in collaboration with S. Rodr\'iguez-L\'opez \cite{IRS}, showed that for any $s\in \Rl$, $\displaystyle \frac{n}{n+1}<p\leq \infty$ and $0<q\leq \infty$, Fourier integral operators with amplitudes $a\in S^{m}_{1,0}(\Rl^n)$ and strongly non-degenerate phase functions in class $\Phi^2$ (see Definition \ref{def:phi2}) satisfy
\begin{equation*}
	T_a^\varphi: B_{p,q}^{s+m-m_c}(\Rl^n)\rightarrow B_{p,q}^s(\Rl^n),
\end{equation*}
 provided that $m_c=-(n-1)\big|\frac{1}{p}-\frac{1}{2}\big|.$
Moreover, under the same conditions on the phase, and the parameters $s$ and $p$, it was also shown that $T_a^\varphi$ is bounded from $F_{p,q}^{s}(\Rl^n) $ to $ F_{p,q}^s(\Rl^n),$
provided that $a\in S^{-(n-1)|\frac{1}{p}-\frac{1}{2}|}_{1,0}(\Rn)$ and $\min\, \{2,p\}\leq q\leq \max\,\{ 2,p\}.$ This is done by complex interpolation in the vertical direction between $F^s_{p,p}(\Rl^n)=B^s_{p,p}(\Rl^n)$ and $F^s_{p,2}(\Rl^n)$ (as in {\bf Figure \ref{pic:gammalskåpmat}}).\\

However for the range of $q\in [0,\infty]$, the best result so far was confined to the amplitudes $a\in S^{-(n-1)|\frac{1}{p}-\frac{1}{2}|-\varepsilon}_{1,0}(\Rn)$ for any $\varepsilon >0$ and hence the expected optimal end-point regularity result is missing. However, as it was observed by M. Christ and A. Seeger \cite{CS}, it behoves one to restrict the range of $q$'s in dealing with Triebel-Lizorkin boundedness of Fourier integral operators of order $m=-(n-1)\big|\frac{1}{p}-\frac{1}{2}\big|$ and the results for arbitrary values of $p$ and $q$ are actually false. We also mention in passing that the boundedness with $\varepsilon-$loss is also valid for operators with more general amplitudes, see e.g. Lemma \ref{lem:local and global nonendpoint TL} in this paper. \\

\vspace*{.1cm}
\vspace*{.1cm}

The discussion of the boundedness results above suggests that the regularity estimates for the Fourier integral operators could be improved in three ways. Namely with regard to the ranges of the scales of the Besov-Lipschitz and Triebel-Lizorkin spaces, i.e. $s$, $p$, and $q$, or with regard to the orders and types of the amplitude i.e. $m,$ $\rho$ and $\delta$, or finally with respect to the rank of the Fourier integral operator.\\ In this paper, we have made an attempt to make improvements and achieve optimal results in all three of these directions.\\ 

The motivation for considering Fourier integral operators of types different than $\rho=1$, $\delta=0$, and ranks different than $n-1$ comes from the theory of partial differential equations, scattering theory, inverse problems, and tomography, just to name a few. In the paper \cite{CIS} the authors divided the amplitudes of the operators according to the following scheme. One calls the amplitudes $a\in S^m_{\rho,\delta}(\Rl ^n)$ with  type $\rho\in (0,1]$ and $\delta\in [0,1)$, {\it classical}, the amplitudes with $\rho=0$ and $\delta\in [0,1)$, {\it exotic}, and amplitudes with the type $\rho\in [0,1]$ and $\delta=1$, {\it forbidden}.\\

In this paper, using certain decomposition in the frequency space, a composition theorem for the action of parameter-dependent pseudodifferential operators on Fourier integral operators (significantly more general than the local one in \cite{Hor:acta}), atomic and molecular decompositions of Triebel-Lizorkin spaces in the spirit of Frazier-Jawerth \cite{FJ:phi-transform, FJ:discrete-transform, FJ-phi-and-wavelet-transform}, and vector-valued inequalities of Pramanik-Rogers-Seeger, we manage to get a significant extension of the results  in \cite{PRS} and \cite{IRS} mentioned above. Compare for instance {\bf Figure \ref{pic:gammalskåpmat}} which illustrates previous results, to {\bf Figure \ref{pic:TLendpointresults}} (in Subsection \ref{subsec:FIO_classic_exotic}) that illustrates the extensions that are obtained here.\\ 
As an application of our results, we can also consider the case of oscillatory integral operators that are related to the Klein-Gordon equation. These operators correspond to the operators $T_a^\varphi$ where $\varphi(x, \xi)= x\cdot \xi+ \langle \xi\rangle$. Although these operators are not Fourier integral operators (due to the lack of homogeneity of their phase function), as we will see in Theorem \ref{linearhpthm}, matters could be reduced to the case of Fourier integral operators.\\

More specifically in Theorems \ref{thm:TLpgeq2FIO}, \ref{thm:Sobolev_fio0}, \ref{thm:Sobolev_fio1} (and Theorems \ref{theorem:BL-booster theorem}, \ref{thm:BS_fio0}, \ref{thm:BS_fio1}) we obtain regularity results for Fourier integral operators with amplitudes in classical, exotic and forbidden scales of the class $S^m_{\rho,\delta}(\Rn)$, and with phase functions of arbitrary rank, on Triebel-Lizorkin (and Besov-Lipschitz spaces). These results are all valid up to the endpoints. Moreover, we also obtain regularity results for operators whose amplitudes belong to the worst possible H\"ormander class $S^m_{0,1}(\Rn)$ and also operators with amplitudes in the forbidden class $S^m_{1,1}(\Rn)$ for which even the endpoint $L^2$--boundedness is in general impossible. It was shown in \cite{CIS} that for the operators with forbidden amplitudes, one could establish the analogue of the Meyer-Stein's result for pseudodifferential operators, namely that Fourier integral operators with amplitudes in $S^0_{1,1}(\Rn)$ and phase functions of rank $n-1$ are bounded on Sobolev spaces $H^s(\Rl^n)$ for any $s>0$. In this paper we extend this result (Theorem \ref{thm:Sobolev_fio1}) to operators of rank $\kappa\in \{0,1, \dots, n-1\}$, with order $m=-\kappa\big|\frac{1}{p}-\frac{1}{2}\big|$ and type $\rho=\delta=1,$ on Triebel-Lizorkin spaces $F^{s}_{p,q}(\Rl^n)$ with $s>n\big(\frac{1}{\min\{1, p, q\}}-1\big)$ and admissible ranges of $p$'s and $q$'s. This is a substantial extension of the results of Triebel \cite{TriebelFIO}, and the second and the third author with Rodr\'iguez-L\'opez \cite{CIS}, mentioned above.\\

The paper is organized as follows; in Section \ref{prelim} we provide the necessary preliminaries from both Fourier- and microlocal analysis. We also include the necessary background from the theory of function spaces. In Section \ref{subsec_composition}
we prove a parameter-dependent composition (or calculus) theorem which will be crucial in proving regularity results in Besov-Lipschitz and Triebel-Lizorkin spaces for operators with amplitudes in general H\"ormander classes. In Section \ref{SSS decomposition} we introduce the decomposition of the Fourier integral operators of constant rank, which differs from thr one given in \cite{SSS}, by the fact that it works for operators with amplitudes of type $\rho\in [0,\frac{2}{3}]$. Note that the decomposition given in \cite{SSS} works for $\rho\in [\frac{1}{2}, 1]$. This section also contains crucial estimates for kernels of the operators. In Section \ref{Sec:FIO estimates} we estimate the Fourier integral operators themselves, using atomic and molecular decompositions, and prove our main regularity theorems for these operators in Triebel-Lizorkin and Besov Lipschitz spaces. Finally, in Section \ref{Klein-Gordon} we deal with the regularity of oscillatory integral operators related to the Klein-Gordon equation and show that similar results, as in the case of the wave equation, are valid in that case.\\

{{\bf{Acknowldgements.}}
The second author is  supported by the Knut and Alice Wallenberg Foundation.
The authors are also grateful to Andreas Str\"ombergsson for his support and encouragement.}

\section{Preliminaries}\label{prelim}

As is common practice, we will denote positive constants in the inequalities by $C$, which can be determined by known parameters in a given situation but whose
value is not crucial to the problem at hand. Such parameters in this paper would be, for example, $m$, $p$, $s$, $n$,  and the constants connected to the seminorms of various amplitudes or phase functions. The value of $C$ may differ
from line to line, but in each instance could be estimated if necessary. We also write $a\lesssim b$ as shorthand for $a\leq Cb$ and moreover will use the notation $a\sim b$ if $a\lesssim b$ and $b\lesssim a$.\\

We start by recalling the definition of the Littlewood-Paley partition of unity which is the most basic tool in the frequency decomposition of the operators at hand.
\begin{Def}\label{def:LP}
 Let $\psi_0 \in \mathcal C_c^\infty(\Rl^n)$ be equal to $1$ on $B(0,1)$ and have its support in $B(0,2)$. Then let
$$\psi_j(\xi) := \psi_0 \left (2^{-j}\xi \right )-\psi_0 \left (2^{-(j-1)}\xi \right ),$$
where $j\geq 1$ is an integer and $\psi(\xi) := \psi_1(\xi)$. Then $\psi_j(\xi) = \psi\left (2^{-(j-1)}\xi \right )$ and one has the following Littlewood-Paley partition of unity

\begin{equation*}
    \sum_{j=0}^\infty \psi_j(\xi) = 1, \quad \text{\emph{for all }}\xi\in\Rl^n .
\end{equation*}

\noindent It is sometimes also useful to define a sequence of smooth and compactly supported functions $\Psi_j$ with $\Psi_j=1$ on the support of $\psi_j$ and $\Psi_j=0$ outside a slightly larger compact set. One could for instance set
\begin{equation*}
\Psi_j := \psi_{j+1}+\psi_j+\psi_{j-1},
\end{equation*}
with $\psi_{-1}:=\psi_{0}$.
\end{Def} \hspace*{1cm}\\
In what follows we define the Littlewood-Paley operators by

\begin{equation*}
    \psi_j(D)\, f(x)= \int_{\Rl^n}   \psi_j(\xi)\,\widehat{f}(\xi)\,e^{ix\cdot\xi}\, \ddd\xi,
\end{equation*}

where $\ddd \xi$ denotes the normalised Lebesgue measure ${\dd \xi}/{(2\pi)^n}$ and
\begin{equation*}
	\widehat{f}(\xi)=\int_{\Rl^n} e^{-i x\cdot\xi}\,f(x) \dd x,
\end{equation*}
is the Fourier transform of $f$.

\subsection{Brief theory of classical function spaces}
\noindent
Here we recall the definition and some basic facts about Besov-Lipschitz and Triebel-Lizorkin spaces that will be used throughout this paper. We begin by defining the mixed norms of functional sequences $\{f_k\}$, $k\in \mathbb{Z}_+$, via
\eq{
    \|\{f_k\}\|_{L^p(\ell^q)}&:=\Big\|\sum_{k} |f_k(\cdot)|^q\Big\|_{L^p(\Rn)}^{1/q},\\
    \|\{f_k\}\|_{\ell^q(L^p)}&:=\Big(\sum_{k} \|f_k\|_{L^p(\Rn)}^q\Big)^{1/q}.
}
Now using the Littlewood-Paley decomposition of Definition \ref{def:LP}, one defines the \emph{Triebel-Lizorkin space} $F^s_{p,q}(\Rl^n)$ and \emph{Besov-Lipschitz space} $B^s_{p,q}(\Rl^n)$ as follows:

\begin{Def}\label{def:TLspace}
	Let $s \in {\Rl}$ and $0< p <\infty$, $0< q \leq\infty$. The Triebel-Lizorkin space is defined by
	\[
	F^s_{p,q}(\Rl^n)
	:=
	\Big\{
	f \in {\mathscr{S}'}(\Rl^n) \,:\,
	\|f\|_{F^s_{p,q}(\Rl^n)}
	:=
	\|2^{js}\psi_j(D) f\|_{L^p(\ell^q)}<\infty
	\Big\},
	\]
where $\mathscr{S}'(\Rl^n)$ denotes the space of tempered distributions.	
\end{Def}
 
\begin{Def}\label{def:Besov}
	Let $0 < p,q \le \infty$ and $s \in {\mathbb R}$. The Besov-Lipschitz spaces are defined by
	\[
	{B}^s_{p,q}(\Rl^n)
	:=
	\Big\{
	f \in {{\SS}'(\Rl^n)} \,:\,
	\|f\|_{{B}^s_{p,q}(\Rl^n)}
	:=
	\|2^{js}\psi_j(D) f\|_{\ell^q(L^p)}<\infty
	\Big\}.
	\]
\end{Def}

\begin{Rem}\label{rem:TLspace}
Different choices of the basis $\{\psi_j\}_{j=0}^\infty$ give equivalent norms of $F^s_{p,q}(\Rl^n)$ in \emph{Definition \ref{def:TLspace},} see e.g. \cite{Triebel1}. We will use either $\{\psi_j\}_{j=0}^\infty$ or  $\{\Psi_j\}_{j=0}^\infty$ to define the norm of $F^s_{p,q}(\Rl^n)$.
\end{Rem}
\noindent We note that for $p=q=\infty$ and $0<s\leq 1$ we obtain the familiar Lipschitz space $\Lambda^s(\Rl^n)$, i.e. $B^s_{\infty,\infty}(\Rl^n)= \Lambda^s(\Rl^n)$. For $-\infty <s<\infty$ and $1\leq p<\infty,$ $F^s_{p,2}(\Rl^n)=H^{s,p}(\Rl^n)$ (various $L^p$-based Sobolev and Sobolev-Slobodeckij spaces) and for $0<p<\infty$, $F^0_{p,2}(\Rl^n)=h^p(\Rl^n)$ (the local Hardy spaces). Moreover the dual space of $F^{0}_{1,2}(\Rl^n)$ is $\mathrm{bmo}$ (the local version of $\mathrm{BMO}$).\\

Some other facts which will be useful to us are

\begin{enumerate}
    \item[$(i)$] For $-\infty <s<\infty$ and $0<p\leq \infty$ one has
    \begin{equation}\label{equality of TL and BL}
        B^s_{p,p}(\Rl^n)= F^s_{p,p}(\Rl^n),
    \end{equation}
    and
    \begin{equation}\label{embedding of TL}
        F^{s+\varepsilon}_{p,q_0}(\Rl^n)\xhookrightarrow{} F^s_{p ,q_1}(\Rl^n)\quad\text{and}\quad
        B^{s+\varepsilon}_{p,q_0}(\Rl^n)\xhookrightarrow{} B^s_{p ,q_1}(\Rl^n) 
    \end{equation}
    for $-\infty <s<\infty$, $0<p< \infty$, $0<q_0,q_1 \leq \infty$ and all $\varepsilon>0$.\\
    \item[$(ii)$] For $s'\in \Rl$, the operator $ (1-\Delta)^{s'/2}$ maps ${F}^s_{p,q}(\Rl^n)$ isomorphically into ${F}^{s-s'}_{p,q}(\Rl^n)$ and ${B}^s_{p,q}(\Rl^n)$ isomorphically into ${B}^{s-s'}_{p,q}(\Rl^n),$ see \cite[p. 58]{Triebel1}.\\
    \item[$(iii)$] Let $X^s_{p,q}(\Rn)$ be either $B^s_{p,q}(\Rn)$ or $F^s_{p,q}(\Rn)$. Then for all $-\infty <s<\infty$ and $0<p,q\leq \infty$ one has 
    \begin{equation}\label{pointwise multiplier}
        \|fg\|_{X^s_{p,q}(\Rl^n)}\lesssim \Big(\sum_{|\alpha|\leq M}\sup_{x\in\Rl^n}|\d^\alpha f(x)|\Big)\|g\|_{X^s_{p,q}(\Rn)}.
    \end{equation}
\end{enumerate}

We will also need the following result (in the proof of our Lemma \ref{Usingequivalenceofnormslemma}). The proof could be found in \cite{Triebel1}.

\begin{Th}\label{charbyapprox}
Let
\eq{
    \mathfrak{U}_p(\Rn):=\{&a;\, a=\{a_k\}_{k=0}^\infty \subset \mathscr{S}'(\Rn)\cap L^p(\Rn),\\
    &\textnormal{supp}\, \widehat{a_k}\subseteq B(0,2^{k+1})\}.
}
Such sequences $\{a_k\}_{k=0}^\infty $ are referred to as admissible sequences.

\begin{enumerate}
    \item[$(i)$] If $0<p,q< \infty$ and $s>n\Big(\frac{1}{\min\{p,q,1\}}-1\Big)$ then
\begin{align*}
    F^s_{p,q}(\Rn)=\{&f;\, f\in \mathscr{S}'(\Rn),\,\exists a=\{a_k\}_{k=0}^\infty \subset\mathfrak{U}_p(\Rn),\, f=\lim_{k}a_k,\,\\
    &\|a_0\|_{L^p(\Rn)}+\|2^{sk}(f-a_k)\|_{L^p(l^q)}<\infty\}.
\end{align*}
Furthermore, $\|f\|_{F^s_{p,q}}=\inf_{a}(\|a_0\|_{L^p(\Rn)}+\|2^{sk}(f-a_k)_{k\geq 0}\|_{L^p(l^q)})$ is an equivalent quasi-norm in $F^s_{p,q}.$
    \item[$(ii)$] If $0<p,q\leq \infty$ and $s>n\Big(\frac{1}{\min\{p,1\}}-1\Big)$ then
\begin{align*}
    B^s_{p,q}(\Rn)=\{&f;\, f\in \mathscr{S}'(\Rn),\,\exists a=\{a_k\}_{k=0}^\infty \subset\mathfrak{U}_p(\Rn),\, f=\lim_{k}a_k,\,\\
    &\|a_0\|_{L^p(\Rn)}+\|2^{sk}(f-a_k)\|_{l^q(L^p)}<\infty\}.
\end{align*}
Furthermore, $\|f\|_{B^s_{p,q}}=\inf_{a}(\|a_0\|_{L^p(\Rn)}+\|2^{sk}(f-a_k)_{k\geq 0}\|_{l^q(L^p)})$ is an equivalent quasi-norm in $B^s_{p,q}.$
\end{enumerate}
\end{Th}

We also need the following result, a proof of which can be found in \cite{RS}.

\begin{Prop}\label{prop:runst}For $\nu=\min\{1, p,q\}$ one has
\begin{equation}
\norm{\sum_{j=0}^{\infty} f_{j}}_{F_{p, q}^{s}(\Rn)} \leq\Big( \sum_{j=0}^{\infty}\Vert f_{j}\Vert_{F_{p, q}^{s}(\Rn)}^{\nu}\Big)^{1/\nu}
\end{equation}
and
\begin{equation}
\norm{\sum_{j=0}^{\infty} f_{j}}_{B_{p, q}^{s}(\Rn)} \leq\Big( \sum_{j=0}^{\infty}\Vert f_{j}\Vert_{B_{p, q}^{s}(\Rn)}^{\nu}\Big)^{1/\nu}.
\end{equation}
\end{Prop}

In connection to estimates for linear operators in Triebel-Lizorkin spaces, one often encounters various maximal functions one is the well-known Hardy-Littlewood's maximal function
 \begin{equation*}\label{HLmax}
  \mathcal{M} f(x):=\sup _{B\ni x}\frac{1}{|B|} \int_{B}|f(y)| \dd y,   
 \end{equation*} where the supremum is taken over all balls $B$ containing $x$. For \(0<p<\infty\), one also defines 
$$(\mathcal{M}_{p} f(x):=\left(\mathcal{M}\left(|f|^{p}\right)\right)^{1 / p}.$$
The other one is J. Peetre's maximal operator \cite{Triebel1}.
\begin{equation}\label{Peetremax}
  \mathfrak{M}_{a,b}(f)(x):=\Big\Vert  \frac{{f(x-\cdot)}}{(1+b|\cdot|)^a)}\Big\Vert_{L^{\infty}(\Rn) } 
\end{equation}
where $0<a, b<\infty$. For any $x\in \Rl^n$, $f\in\mathscr{S}'(\Rn)$ with $\supp \widehat{f}\subset \{\xi;\,|\xi| \leq 2b \}$ and $a\geq \frac{n}{p}$ one has that

\begin{equation}\label{hl bounds peetre}
\mathfrak{M}_{a, b} u(x) \lesssim \mathcal{M}_{p} u(x).
\end{equation}
The main tool involving maximal functions in the context of Triebel-Lizorkin spaces is the following Fefferman-Stein vector valued inequality, see \cite{FS} for the proof.
\begin{Th}\label{FSGeneralized thm}
    Let $ 0<p<\infty, \, 0<q \leq \infty, \, 0<r<\min \{p, q\}$. Then for a sequence $(g_k)_{k\geq 0}$ one has 
    \begin{equation}
        \norm{\Big(\sum_{k=0}^{\infty}|\mathcal{M}_r g_{k}(\cdot)|^{q}\Big)^{1/q}}_{L^{p}(\mathbb{R}^{n})}\lesssim \norm{\Big(\sum_{k=0}^{\infty}\left|g_{k}(\cdot)\right|^{q}\Big)^{1 / q}}_{ L^{p}\left(\mathbb{R}^{n}\right)}.
    \end{equation}
\end{Th}
We shall now present the following rather abstract lemma which will be used in connection to the boundedness of FIOs with amplitudes in the classes $S^m_{0,1}$ and $S^m_{1,1}$.

\begin{Lem}\label{Usingequivalenceofnormslemma}
    Let $X^{s}_{p,q}$ be either $F^{s}_{p,q}(\Rn)$ or $B^{s}_{p,q}(\Rn)$. Assume that $(u_k)_{k\geq 1}\in X^{s}_{p,q}$ is an arbitrary sequence such that 
    \begin{equation}\label{boundednessassumption}
        \big\| \{2^{ks}u_k\}_{k=0}^\infty\big \|_{_{L^p(l^q)}}\lesssim \Vert f\Vert_{X^{s}_{p,q}}.
    \end{equation}
    Moreover let $\{h_{k,l}(x)\}_{k,l\geq 0}$ be a sequence such that the spectrum of $h_{k,l}(x)$ $($i.e. the support of its Fourier transform$)$ is in $B(0,2^{k+l+2})$ and that there is some sufficiently large $N$ such that $h_{k,l}$ satisfies the pointwise estimate,
    \begin{equation}\label{piecewiseest}
        |h_{k,l}(x)|\lesssim 2^{-Nl}\,\mathcal{M}_{r}u_k(x),\quad\mathrm{for\,\, some}\,\, r \in (0,\min\{p,q\}).
    \end{equation}
    Let $g_l=\sum_{k\geq 0}h_{k,l}$.
    If $0<p,q<\infty$ and $s>n(\frac{1}{\min\{ p, q,1\}}-1)$ then
    \begin{equation}
        \norm{\sum_{l\geq 0} g_l }_{X^{s}_{p,q}} \lesssim \Vert f\Vert_{X^{s}_{p,q}}.
    \end{equation}
\end{Lem}

\begin{proof}
Observe first that by Proposition \ref{prop:runst} one has
\begin{equation}\label{firstprop:runst}
\norm{\sum_{l\geq 0} g_l}_{F_{p, q}^{s}(\Rn)} \leq\Big( \sum_{l\geq 0} \Vert g_{l} \Vert_{F^{s}_{p, q}}^{\nu}\Big)^{1/\nu},
\end{equation}
where $\nu=\min\{1,p,q\}$.\\

In order to approximate $\|g_l\|_{F^s_{p,q}}$ we shall apply Theorem \ref{charbyapprox}. For fixed (but otherwise arbitrary $l$) define a sequence $a=\{a_{k,l}\}_{k\geq 0}$ by
\begin{align}
    a_{k,l}=\begin{cases}
    0& k\leq l-1\\
    \sum_{j=0}^{k-1} h_{j,l}(x) & k\geq l.
    \end{cases}
\end{align}
Note that since the spectrum of $h_{k,l}(x)$ is in $B(0,2^{k+l+2})$ it follows that $(a_{k,l})_{k\geq 0}$ is an admissible sequence in the sense of Theorem \ref{charbyapprox}.\\

Thus Theorem \ref{charbyapprox} yields,
\begin{equation}\label{splitequationgl}
    \|g_l\|_{F^s_{p,q}}\lesssim \|a_{0,l}\|_{L^p(\Rn)}+\norm{2^{sk}(g_l-a_{k,l})}_{L^p(l^q)},
\end{equation}
note that $a_{0,l}=0$ for all $l\geq 0$. Now, split the last summand in \eqref{splitequationgl} into two parts,
\begin{equation*}\label{es0}
    \|2^{sk}(g_l-a_{k,l})_{k\geq 0}\|_{L^p(l^q)}\leq \norm{\{2^{ks}(g_l-a_{k,l})\}_{k=0}^{j-1}}_{L^p(l^q)}+\norm{\{2^{ks}(g_l-a_{k,l})\}_{k=j}^{\infty}}_{L^p(l^q)}
\end{equation*}

For the first part, using \eqref{piecewiseest}, \eqref{boundednessassumption} and the Fefferman-Stein Theorem \ref{FSGeneralized thm} (in index $k$) we obtain
\begin{align}\label{es1}
    &\norm{\{2^{ks}(g_l-a_{k,l})\}_{k=0}^{j-1}}_{L^p(l^q)}
    =\norm{\big(\sum_{k=0}^{l-1} 2^{ksq} \big)^{1/q} g_l}_{L^p(\Rn)}
    \lesssim 2^{ls}\|g_l\|_{L^p(\Rn)}\\
    &\nonumber= 2^{ls}\norm{\sum_{k\geq 0} h_{k,l}(x)}_{L^p(\Rn)}
    \lesssim 2^{ls}\norm{\{2^{ks}h_{k,l}(x)\}_{k=0}^{\infty}}_{L^p(l^q)}\\
    &\nonumber\lesssim  2^{l(s-N)}\norm{\{2^{ks}\mathcal{M}_r u_k(x)\}_{k=0}^{\infty}}_{L^p(l^q)}
    \lesssim  2^{l(s-N)}\norm{\{2^{ks} u_k(x)\}_{k=0}^{\infty}}_{L^p(l^q)}\\
    &\nonumber\lesssim 2^{(s-N)l}\|f\|_{F^s_{p,q}}.
\end{align}

While for the second part we have for any $\eps>0$
\begin{align}\label{es2}
    &\norm{\{2^{ks}(g_l-a_{k,l})\}_{k=l}^{\infty}}_{L^p(l^q)}
    =\norm{\{2^{ks}\Big(\sum_{j=k-(l-1)}^\infty h_{k,l}(x) \Big)\}_{k=l}^{\infty}}_{L^p(l^q)}\\
    &=\nonumber\norm{\{2^{ks}\Big(\sum_{d=1}^\infty h_{k-l+d,l} \Big)\}_{k=l}^{\infty}}_{L^p(l^q)}\\
    &=\nonumber\norm{\{2^{(l-d+\eps)s}2^{(k-l+d-\eps)s}\Big(\sum_{d=1}^\infty h_{k-l+d,l} \Big)\}_{k=l}^{\infty}}_{L^p(l^q)}\\
    &\leq\nonumber\norm{\{2^{(l-d+\eps)s}2^{(k-l+d)s}\Big(\sum_{d=1}^\infty h_{k-l+d,l} \Big)\}_{k=l}^{\infty}}_{L^p(l^q)}\\
    &\nonumber\leq\sum_{d=1}^\infty \Big\{2^{(l-d+\eps)s}\norm{\{2^{(k-l+d)s}h_{k-l+d,l} \}_{k=l}^{\infty}}_{L^p(l^q)}\Big\}\\
    &\nonumber\leq\sum_{d=1}^\infty \Big\{2^{(l-d+\eps)s}2^{-Nl}\|f\|_{F^s_{p,q}}\Big\}\\
    &\nonumber\lesssim_\eps 2^{(s-N) l}\|f\|_{F^s_{p,q}}
\end{align}
where we have used \eqref{piecewiseest}. Now \eqref{es1} and \eqref{es2} together with \eqref{splitequationgl} yield that
\begin{equation}\label{someeqaf}
    \Vert g_{l} \Vert_{F^{s}_{p, q}}\lesssim 2^{(s-N) l}\|f\|_{F^s_{p,q}}.
\end{equation}
Plugging \eqref{someeqaf} into \eqref{firstprop:runst} and taking $N>s$ one obtains
\begin{equation*}
    \norm{\sum_{l\geq 0} g_l}_{F_{p, q}^{s}(\Rn)} \leq\Big( \sum_{l\geq 0} 2^{(s-N)\nu l}\Vert f \Vert_{F^{s}_{p, q}}^{\nu}\Big)^{1/\nu}\lesssim \Vert f \Vert_{F^{s}_{p, q}}.
\end{equation*}
The proof for Besov-Lipschitz spaces follows from a similar argument. 
\end{proof}

Another important and useful fact about Besov-Lipschitz and Triebel-Lizorkin spaces is the following:

\begin{Th}\label{thm:invariance thm}
Let $\eta: \Rn \to\Rn$ with $\eta(x)=(\eta_1 (x), \dots,\eta_n (x))$ be a  diffeomorphism,  such that $|\det D\eta (x)|\geq c>0$, $\forall x\in \Rn$ \emph{(}$D\eta$ denotes the Jacobian matrix of $\eta$\emph{)}, and $\Vert\partial^{\alpha}\eta_j (x)\Vert_{L^\infty(\Rn)}\lesssim 1$ for all $j\in \{1,\dots, n\}$ and $|\alpha|\geq 1.$ Then
for $s\in \Rl$, $0<p<\infty$ and $0<q\leq \infty$ one has $$\Vert f\circ \eta\Vert_{F^{s}_{p,q}(\Rn)}\lesssim \Vert f\Vert_{F^{s}_{p,q}(\Rn)}.$$

The same invariance estimate is also true for Besov-Lipschitz spaces $B^{s}_{p,q}(\Rn)$ for  $s\in \Rl$, $0<p\leq\infty$ and $0<q\leq \infty$.
\end{Th}
For a proof see J. Johnsen, S. Munch Hansen and W. Sickel \cite[Corollary 25]{JMHS}, and  H. Triebel \cite[Theorem 4.3.2]{Triebel2}. \\

The rest of this section is dedicated to setting up the stage for atomic Hardy spaces which will be used in connection to $h^p\to L^p$ estimates for Fourier integral operators.\\

For Definitions \ref{Def:Qbe} and \ref{def:hpatom}, let $[x]$ denote the integer part of $x$.
\begin{Def}\label{Def:Qbe}
For a closed cube $Q$, define
\begin{enumerate}
    \item[$(i)$] $c_Q$ as the centre of $Q$,
    \item[$(ii)$] $l_Q$ as the side length of $Q$,
    \item[$(iii)$]     $k_Q:=[1-\log_2(l_Q)]$, 
    \item[$(iv)$] $\chi_Q$ as the characteristic function of $Q$, i.e. $\chi_Q(x):=\diffcases{1,& x\in Q\\0,& x\notin Q}$.
    \item[$(v)$] $cQ$ as the cube with centre $c_Q$ and length $l_{c\,Q}=c\,l_Q$, where $c\in \Rl_{>0}$.
\end{enumerate}

\end{Def}
Observe that $k_Q$ is the unique integer such that $$2^{-k_Q}< l_Q\leq 2^{-(k_Q-1)}.$$
\begin{Def}\label{def:hpatom}
A function $\at$ is called a $h^p$-atom if there exists a cube $Q$ such that the following three conditions are satisfied:
\begin{enumerate}
\item[$(i)$] $\supp \at\subset Q$,
\item[$(ii)$] $\displaystyle \sup_{x\in Q}|\at(x)|\leq 2^{k_Qn/p},$
\item[$(iii)$] If $k_Q\geq 1$, $\displaystyle  M_{\at}= \left[ n\brkt{\frac 1p -1} \right ]$, 
then $\displaystyle \int_{\Rn} x^{\alpha}\at(x)\dd x=0,$ for $|\alpha|\leq M_{\at}$. No further condition is assumed if $k_Q\leq 0.$ 
\end{enumerate}

It is well known $($see \cite{Triebel1}$)$ that a distribution $f\in h^p (\Rn)$ has an atomic decomposition
\nm{eq:hpsumdefinition}{
    f=\sum_{j=0}^\infty\lambda_{j}\at_{j},
}
where the $\lambda_{j}$ are constants such that $$\displaystyle \inf_{\set{\lambda_j}}\sum_{j=0}^\infty|\lambda_{j}|^{p}\sim\Vert f\Vert_{h^p(\mathbb{R}^{n})}^{p}=\Vert f\Vert_{F_{p,2}^0(\mathbb{R}^{n})}^{p}$$ and the $\at_{j}$ are $h^p$-atoms.\\

\end{Def}

\subsection{Basic definitions related to the Fourier integral operators}

Next, we define the building blocks of the FIOs and the pseudodifferential operators. These are the amplitudes (symbols in the pseudodifferential setting) and the phase functions. The class of amplitudes considered in this paper was first introduced by L. H\"ormander in \cite{Hor1}.
\begin{Def}\label{symbol class Sm}
Let $m\in \Rl$ and $\rho, \delta \in [0,1]$. An \textit{amplitude} \emph{(}symbol\emph{)} $a(x,\xi)$ in the class $S^m_{\rho,\delta}(\Rl^n)$ is a function $a\in \mathcal{C}^\infty (\Rl^n\times \Rl^n)$ that verifies the estimate
\begin{equation*}
\left|\partial_{\xi}^\alpha \partial_{x}^\beta a(x,\xi) \right| \lesssim \langle\xi\rangle ^{m-\rho|\alpha|+\delta|\beta|},
\end{equation*}
for all multi-indices $\alpha$ and $\beta$ and $(x,\xi)\in \Rl^n\times \Rl^n$, where 
$\langle\xi\rangle:= (1+|\xi|^2)^{1/2}.$\\
We shall henceforth refer to $m$ as the order of the amplitude and $\rho, \, \delta$ as its type. We will refer to the class  $S_{\rho,\delta}^m(\Rl^n)$ with $0<\rho\leq 1$, $0\leq \delta<1$ as classical, to the class $S_{0,\delta}^m(\Rl^n)$ with $0\leq\delta<1$ as the exotic class, and to $S_{\rho,1}^m(\Rl^n)$ with $0\leq\rho\leq 1$ as the forbidden class of amplitudes. \\

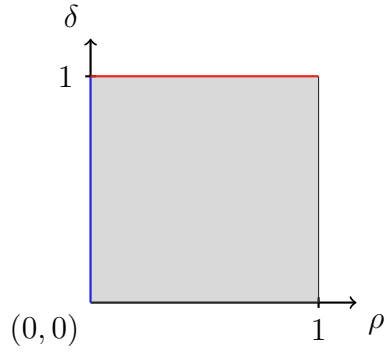
\begin{figure}
\begin{tikzpicture}

    \draw[thick,->] (0,0) -- (3.5,0) node[anchor=north west] {$\rho$};
    \draw[thick,->] (0,3) -- (0,3.5) node[anchor=south east] {$\delta$};
    \draw[thick,blue] (0,0) -- (0,3);

    \draw[thick] (3 cm,2pt) -- (3 cm,-2pt) node[anchor=north] {$1$};
    \draw[thick] (2pt,3 cm) -- (-2pt,3 cm) node[anchor=east] {$1$};
    \draw (0,0) node[anchor=north east] {$(0,0)$};

    \draw (3 cm, 0 cm) -- (3 cm,3 cm);
    \draw[red,thick] (0 cm, 3 cm) -- (3 cm,3 cm);
    
    \fill[gray, opacity = 0.3] (0,0) rectangle (3,3);
\end{tikzpicture}
\caption{The Hörmander classes. The classes on the blue line are called \textit{exotic} and the classes on the {red line} are called \textit{forbidden}.}
\end{figure}

\end{Def}

Towards the end of this paper, in connection with the forbidden amplitudes, we will use the \emph{Zygmund class} $C_{*}^{r} (\Rl^n)$ whose definition we now recall.
\begin{Def}\label{def:Zygmund}
	Let $r \in \Rl$. The Zygmund class is defined by
	\[
	C_{*}^{r} (\Rl^n)
	:=
	\Big\{
	f \in {\mathscr{S}'}(\Rl^n) \,:\,
	\|f\|_{C_{*}^{r} (\Rl^n)}
	:=
	\sup_{j \geq 0}2^{jr}\|\psi_j(D)f\|_{L^\infty(\Rl^n)}
	<\infty\Big\}.
	\]
\end{Def}
If $C^r(\Rl^n)$, $(r\in \Rl_+)$ denotes the H\"older space, and $\mathcal{C}^r(\Rl^n)$ denotes the space of continuous functions with continuous derivatives of orders up to and including $r$, then one has that 
\begin{equation}\label{Zygmundegenskap}
   C_{*}^{r}(\Rl^n) =\mathcal{C}^r(\Rl^n)  \quad \text{for}\,\,\, r\in \Rl_+ \setminus \Z_+ \quad \text{and}\quad  \mathcal{C}^r(\Rl^n) \subset C_{*}^{r}(\Rl^n) \quad \text{for }\, r\in \Z_+.
\end{equation}

\noindent In connection to the definition of the Zygmund class, there is another class of amplitudes that have low regularity in the $x$-variable, which
were considered by G. Bourdaud in \cite{Bourdaud}.

\begin{Def}\label{symbolclass ZygSm}
\noindent Let $m\in \Rl$, $0\leq \delta\leq 1$ and $r>0$. An \textit{amplitude} \emph{(}symbol\emph{)} $a(x,\xi)$ is in the class $C_{*}^{r} S_{1, \delta}^{m}(\Rl^n)$  if it is  $\mathcal{C}^\infty ( \Rl ^n)$ in the $\xi$ variable and verifies the estimates
\begin{equation*}
    \| \partial_\xi^\alpha a(\cdot,\xi) \|_{L^\infty(\Rl^n)}\lesssim \langle\xi\rangle ^{m-|\alpha|},
\end{equation*}
and 
\begin{equation*}
\| \partial_\xi^\alpha a(\cdot,\xi)\|_{C_{*}^{r}(\Rl^n)}\lesssim \langle\xi\rangle ^{m-|\alpha|+\delta r},
\end{equation*}
for all multi-indices $\alpha$ and $\xi\in \Rl^n$.

\end{Def}
It is important to note that  $S^m_{1,1} (\Rl^n)\subset C_{*}^{r} S_{1,1}^{m}(\Rl^n),$ for all $r>0$, which follows from \eqref{Zygmundegenskap}.\\

\noindent Given the symbol classes defined above, one associate to the symbol its \textit{Kohn-Nirenberg quantization} as follows:
\begin{Def}
Let $a$ be a symbol. Define a pseudodifferential operator \emph{(}$\Psi\mathrm{DO}$ for short\emph{)} as the operator
\begin{equation*}
a(x,D)f(x) := \int_{\Rl^n}e^{ix\cdot\xi}\,a(x,\xi)\,\widehat{f}(\xi) \ddd \xi,
\end{equation*}
a priori defined on the Schwartz class $\mathscr{S}(\Rl^n).$
\end{Def}

In order the define the Fourier integral operators that are studied in this paper, we consider the class of phase functions introduced by D. Dos Santos Ferreira and W. Staubach in \cite{DS}

\begin{Def}\label{def:phi2}
A \textit{phase function} $\varphi(x,\xi)$ in the class $\Phi^k$, $k\in\Z_{>0}$, is a function \linebreak$\varphi(x,\xi)\in \mathcal{C}^{\infty}(\Rl^n \times\Rl^n \setminus\{0\})$, positively homogeneous of degree one in the frequency variable $\xi$ satisfying the following estimate

\begin{equation}\label{C_alpha}
	\sup_{(x,\,\xi) \in \Rl^n \times\Rl^n \setminus\{0\}}  |\xi| ^{-1+\vert \alpha\vert}\left | \partial_{\xi}^{\alpha}\partial_{x}^{\beta}\varphi(x,\xi)\right |
	\leq C_{\alpha , \beta},
	\end{equation}
	for any pair of multi-indices $\alpha$ and $\beta$, satisfying $|\alpha|+|\beta|\geq k.$
\end{Def}\hspace*{1cm}

Following the approaches in e.g. \cite{DS, RLS}, for the global $L^p$ boundedness  results that were established in those papers, we also define the following somewhat stronger notion of non-degeneracy:
\begin{Def}\label{nondeg phase}
One says that the phase function $\varphi(x,\xi)$ satisfies the strong non-degeneracy condition \emph{(}or $\varphi$ is $\mathrm{SND}$ for short\emph{)} if
\begin{equation}\label{eq:SND}
	\big |\det (\partial^{2}_{x_{j}\xi_{k}}\varphi(x,\xi)) \big |
	\geq \delta,\qquad \mbox{for  some $\delta>0$ and all $(x,\xi)\in \mathbb{R}^{n} \times \Rl^n\setminus\{0\}$}.
\end{equation}
\end{Def}

Having the definitions of the amplitudes and the phase functions at hand, one has

\begin{Def}\label{def:FIO}
	A Fourier integral operator  \emph($\mathrm{FIO}$ for short\emph) $T_a^\varphi$ with amplitude $a$ and phase function $\varphi$, is an operator defined \emph{(}once again a-priori on $\mathscr{S}(\Rl^n)$\emph{)} by
	\begin{equation}\label{eq:FIO}
	T_a^\varphi f(x) := \int_{\Rl^n} e^{i\varphi(x,\xi)}\,a(x,\xi)\,\widehat f (\xi) \ddd\xi,
	\end{equation}
	where $\varphi(x,\xi)\in\mathcal{C}^{\infty}(\Rl^n \times \Rl^n\setminus\{0\})$ is positively homogeneous of degree one in $\xi$. We say that $\varphi$ is of rank $\kappa$ if it satisfies $\mathrm{rank}\,\partial^{2}_{\xi\xi} \varphi(x, \xi) = \kappa$ for all $(x,\xi) \in \Rl^n \times \Rl^n\setminus \{0\}$.
\end{Def}

Taking $\varphi(x,\xi)=x\cdot\xi$, one obtains the class of pseudodifferential operators associated to H\"ormanders symbol classes (Definition \ref{symbol class Sm}).

\subsection{Reduction of FIOs with phases in \texorpdfstring{$\Phi^2$}{} to FIOs with phases in \texorpdfstring{$\Phi^1$}{}}\label{sec:phase_reduction}
As was demonstrated in \cite[Section 3]{CIS}, the study of global boundedness of FIOs $T_a^{\varphi}$ with SND phase functions $\varphi\in \Phi^2$ could be reduced to the study of global boundedness of FIOs with SND phase functions $x\cdot \xi+\theta(x,\xi)$ with $\theta\in \Phi^1$. The same is valid for FIOs of the form $$ \iint_{\Rl^n\times \Rl^n} a(y,\xi)\, e^{i\varphi(y,\xi)-ix\cdot\xi}\, f(y) \ddd \xi \dd y,$$with $\varphi\in \Phi^2$, matters can be reduced to the investigation of FIOs of the form

\begin{equation}\label{the_adjoint}
    \iint_{\Rl^n\times \Rl^n} {\sigma(y,\xi)}\, e^{i\theta(y,\xi)+i(y-x)\cdot\xi}\, f(y)  \ddd\xi \dd y,
\end{equation}

where $\sigma(y,\xi)$ belongs to the same class as $a(y,\xi)$ and $\theta(y,\xi)\in \Phi^{1}.$  

\subsection{\texorpdfstring{$L^2$ boundedness of FIOs}{}}
We recall the optimal $L^2$ boundedness result for FIOs with amplitudes in general H\"ormander classes, due to Dos Santos Ferreira and Staubach 
\cite[Theorem 2.7]{DS}, which will be used in the proof of $h^p\to L^p$ boundedness of Fourier integral operators.

\begin{Th}\label{basicL2}
Let  $\rho, \delta\in [0,1]$, $\delta\neq 1$. Assume that $a(x,\xi)\in S^{m}_{\rho,\delta}(\Rn)$ and $\varphi(x,\xi)$ is in the class $\Phi^2$ and is \emph{SND}. Then the \emph{FIO} $T_a^\varphi$ is bounded on $L^2(\Rn)$ if and only if $m=-n \, \max\{0,(\delta-\rho)/2\}$. In case $\rho\in [0,1],$ $\delta=1$ then the $L^2$ boundedness is valid if and only if $m<n(\rho-1)/2$. 
\end{Th}

\subsection{Low-frequency Besov-Lipschitz and Triebel-Lizorkin estimate}

Here we recall the main boundedness result of the low-frequency part of FIOs. To this end let $a_0(x,\xi):= a(x,\xi)\psi_0(\xi)$, where $\psi_0$ is a Littlewood-Paley piece (see definition \ref{def:LP}). Then the following result, due to Israelsson, Rodr\'iguez-L\'opez and Staubach, handles the situation for FIOs with compact frequency support on all scales of Besov-Lipschitz and Triebel-Lizorkin spaces. The proof could be found in  \cite[Proposition 5.7]{IRS} and \cite[Proposition 6.4]{IRS}.

\begin{Th}\label{thm:low_freq_TL_BL_FIO}
For $\rho,\delta\in [0,1]$ let $a(x,\xi)\in S^m_{\rho,\delta}(\Rl^n)$. Suppose also that $0<q_1,q_2 \leq\infty$, $s_1,s_2\in \Rl$. Then if $\varphi (x,\xi)\in \Phi^2$ verifies the $\mathrm{SND}$ condition one has $T_{a_0}^\varphi :F_{p,q_1}^{s_1}(\Rl^n)\to F_{p,q_2}^{s_2}(\Rl^n)$, for $\frac{n}{n+1}<p\leq\infty$
Moreover, the Triebel-Lizorkin estimates above may be replaced by the corresponding Besov-Lipschitz estimate.
\end{Th}

\subsection{A Triebel-Lizorkin boundedness theorem for \texorpdfstring{$p>2$}{}}\label{subsubsec:PRS_boundedness}

The following theorem due to M. Pramnanik, K. Rogers and A. Seeger could be found in \cite[Theorem 2.1]{PRS}.
\begin{Th}\label{thm:PRS}
Let $ 2 < p < \infty$, $0<q\leq \infty$ and $ 0<b<n$. Assume that the operators $S_j$ satisfy
\begin{align}
&\sup_{
j>0}
2^{jb/p}\| S_j\|_{L^p\to L^p} \lesssim 1,\label{eq:PRS1}
\\&
\sup_{
j>0}
2^{jb/2}\|S_j\|_{L^2\to L^2} \lesssim 1.\label{eq:PRS2}
\end{align}
Furthermore, assume that for each cube $\tilde Q$ there is a measurable set $\mathcal E_{\tilde Q}$ and a constant $\Gamma\geq 1$ such that
\eq{
|\mathcal E_{\tilde Q}| \leq \Gamma  \max\{l_{\tilde Q}^{n-b},l_{\tilde Q}^n\},
}
where $l_{\tilde Q}$ is the side length of $\tilde Q$ as in \emph{Definition} \emph{\ref{Def:Qbe}}, and assume further that for every $j\in \mathbb{N}$ and every cube \(\tilde Q\) with \(2^{j} l_{\tilde Q} \geq 1\), one has
\nm{eq:PRS}{
\sup_{
x\in \tilde Q}
\int_{\Rn \setminus \mathcal E_{\tilde Q}} 
|K_j(x, y)| \dd y \lesssim  \max \{(2^{-j\eps}l_{\tilde Q}^{-\eps}, 2^{-j\eps}\},}
for some $\eps>0$.\\

Then
\eq{
\norm{ \brkt{\sum_{j=0}^\infty 
2^{jbq/p}|\Psi_j(D) S_jf_j|^q}^{1/q}}_{L^p(\Rn)} \lesssim \brkt{ \sum_{j=0}^\infty \|f_j\|_{L^p(\Rn)}^p}^{1/p}.
}
where $\Psi\in \mathscr{S}(\Rl^n)$, $\Psi_j(D):=\Psi(2^{-j}D)$ and $f_j$ is a sequence of functions.
\end{Th}

\section{Left-composition of parameter-dependent \texorpdfstring{$\Psi$}{}DOs with FIOs in general Hörmander classes}\label{subsec_composition}
In this section, we state and prove a composition result that will enable us to keep track of the parameter while a parameter-dependent $\Psi$DO acts from the left on an FIO. 
This will be crucial in the proof of the boundedness of FIOs with forbidden amplitudes on Sobolev spaces as well as in other situations.

\begin{Th}\label{thm:left composition with pseudo}
Let $m,s\in \Rl$ $\rho\in[0,1]$, $\delta\in[0,1)$. Suppose that $ a(x, \xi)\in S_{\rho,\delta}^m (\Rl^n)$, $b(x,\xi)\in S^s_{1,0}(\Rl^n)$ and $\varphi$ is a phase function that is smooth on $\supp a$ and verifies the conditions
\begin{align}\label{eq:composition_conditions}
    &|\xi| \lesssim |\nabla_x \varphi(x, \xi)| \lesssim |\xi|,  &|\partial_\xi^\alpha \partial _x^\beta \varphi (x, \xi)| \lesssim  \jap{\xi}^{1-|\alpha|}
\end{align}
for all $(x, \xi) \in \supp a$ and for all $|\alpha| \geq 0,$ and all $|\beta| \geq 1$.
\noindent For  $0<t\leq 1$ consider the parameter-dependent pseudodifferential operator 
\begin{equation*}
b(x, tD)f(x) := \int_{\Rl^n} e^{ix\cdot \xi}\,b(x,t\xi)\,\widehat f(\xi) \ddd \xi,
\end{equation*}
and the oscillatory integral operator
$$T_{a}^\varphi f(x) := \int_{\Rl^n} e^{i\varphi(x,\xi)}\,a(x, \xi)\,\widehat f(\xi) \ddd\xi.$$
Let $\sigma_t$ be the amplitude of the composition operator 
$T_{\sigma_t}^\varphi
:=b(x, tD)T_{a}^\varphi$ 
given by
\begin{equation*}
\sigma_t(x, \xi) := \iint_{\Rl^n\times \Rl^n} a(y, \xi)\, b (x,t\eta)\,e^{i(x-y)\cdot \eta+i\varphi(y,\xi)-i\varphi(x,\xi)} \ddd\eta \dd y.
\end{equation*}
Then for any $M\geq 1$ and all $0<\eps<1-\max\{\delta,1/2\}$, one can write $\sigma_t$ as

\begin{equation}\label{asymptotic expansion}
\sigma_t(x, \xi) =b(x,t\nabla_{x}\phase(x,\xi))\,a(x,\xi) + \sum_{0<|\alpha| < M}\frac{t^{\eps|\alpha|}}{\alpha!}\, \sigma_{\alpha}(t,x,\xi)+t^{\eps M} r(t,x,\xi),
\end{equation}

where, for all multi-indices $\beta, \gamma$ one has
\begin{align*}
&|\partial^{\beta}_{\xi} \partial^{\gamma}_{x}\sigma_{\alpha}(t,x,\xi)|  \lesssim  t^{\min\{s,0\}}  \bra{\xi}^{s+m-(1-\max\{\delta,1/2\}-\varepsilon)|\alpha|-\rho|\beta|+\delta|\gamma|},  \\
&|\partial^{\beta}_{\xi} \partial^{\gamma}_{x} r(t,x,\xi)|   \lesssim t^{\min\{s,0\}} \bra{\xi}^{s+m-\brkt{1-\max\{\delta,1/2\}- \varepsilon}M-\rho|\beta|+\delta|\gamma|}.
\end{align*}
\end{Th}

\begin{proof}
The idea of the proof is similar to that of the asymptotic expansion proved in \cite{RRS}, however the details are somewhat different since the calculus generalises to all amplitude classes $S^m_{\rho,\delta}(\Rn)$ with $\rho\in[0,1]$, $\delta\in[0,1)$. Observe that we can take $M$ as large as necessary since it is possible to merge terms from the main part of \eqref{asymptotic expansion} into the rest term. \\

Let $\chi(x-y)\in \mathcal{C}^{\infty}(\mathbb{R}^n \times \Rl^n)$ such that $0\leq \chi\leq 1,$ $\chi(x-y)\equiv 1$ for $|x-y|< \kappa_2/2$ and $\chi(x-y)=0$ for $|x-y|>\kappa_2$, for some small $\kappa_2$ to be specified later. We now decompose $\sigma_t(x,\xi)$ into two parts $\textbf{I}_1 (t,x,\xi)$ and $\textbf{I}_2 (t,x,\xi)$ where
\begin{equation*}
  \textbf{I}_1 (t,x,\xi)
  :=\iint_{\Rl^n\times \Rl^n} a(y,\xi)\,b(x, t\eta)\, (1-\chi(x-y))\,e^{i(x-y)\cdot\eta+i\phase(y,\xi)-i\phase(x,\xi)}\ddd\eta\dd y,
  \end{equation*}
  and
\begin{equation*}
  \textbf{I}_2 (t,x,\xi)
  :=\iint_{\Rl^n\times \Rl^n} a(y,\xi)\,b(x, t\eta)\,\chi(x-y)\,e^{i(x-y)\cdot\eta+i\phase(y,\xi)-i\phase(x,\xi)}\ddd\eta\dd y.
  \end{equation*}

\quad\\

\textbf{Step 1 -- The analysis of $\mathbf{I_1(t,x,\boldsymbol\xi)}$}\\
To this end, we introduce the differential operators
\[
L_{\eta} 
:=-i\frac{x -y}{|x-y|^2}\cdot \nabla_{\eta} \quad \text{and} \quad
L_{y} 
:=\frac{1}{\langle \nabla_{y}\phase(y,\xi)\rangle ^2 -i\Delta_{y}\phase(y,\xi)}(1-\Delta_{y}).
\]
Because of \eqref{eq:composition_conditions}, one has 
$$|\langle \nabla_{y}\phase(y,\xi) \rangle^2 -i\Delta_{y}\phase(y,\xi)|\geq \langle \nabla_{y}\phase(y,\xi) \rangle^2 \gtrsim \xxi^2.$$
Now integration by parts yields
\eq{
  \textbf{I}_1 (t,x,\xi)&=\iint_{\Rl^n\times \Rl^n} (L_{y}^*)^{N_{2}} \{e^{-iy\cdot\eta} \,a(y,\xi)\, (L_{\eta}^*)^{N_1}[(1-\chi(x-y))\,b(x, t\eta)]\}\\ &\qquad \qquad \times e^{ix\cdot\eta+i\phase(y,\xi)-i\phase(x,\xi)} \ddd\eta \dd y.
}

Now since $0<t\leq 1$, provided $0<N_3<N_1-s$, we have
\begin{align*}
\abs{\partial^{N_1}_{\eta_j}b(x, t\eta)} & \lesssim t^{N_1} \langle t\eta\rangle^{s-N_1}= t^{N_1} \langle t\eta\rangle^{-N_3}\langle t\eta\rangle^{s-(N_1 -N_3)} \\
& \lesssim t^{N_{1}} (t^2+|t\eta|^2)^{-N_3/2} \langle t\eta\rangle^{s-(N_{1}-N_3)}\lesssim t^{N_{1}-N_3}\langle \eta\rangle^{-N_3}.
\end{align*}
Therefore, choosing $N_1 >n$ and $2N_2<N_3 -n$
\eq{
  &| \textbf{I}_1 (t,x,\xi)| \lesssim t^{N_1 -N_3}\xxi^{-2N_2 +m}\iint_{|x-y|>\kappa} \langle \eta\rangle^{2N_2} |x-y|^{-N_1} \langle\eta\rangle^{-N_3} \ddd\eta \dd y \\ &\lesssim  t^{N_1 -N_3} \xxi^{-2N_2 +m}.
}
Estimating derivatives of $\textbf{I}_1(t,x,\xi)$ with respect to $x$ and $\xi$ may introduce factors estimated by powers of $\xxi$, $\langle \eta\rangle$, and $|x-y|$, which can all be handled by choosing $N_1$ and $N_2$ appropriately. Therefore, for all $N$ and any $\nu>0$
\begin{equation*}
  \abs{ \partial^{\alpha}_{\xi} \partial^{\beta}_{x}\textbf{I}_1 (t,x,\xi)} \lesssim t^{\nu}\xxi^{-N},
  \end{equation*}
and so $\textbf{I}_1 (t,x,\xi)$ forms part of the error term $t^{\eps M}r(t,x,\xi)$ in \eqref{asymptotic expansion}.\\

\textbf{Step 2 -- The analysis of $\mathbf{I_2(t,x,\boldsymbol\xi)}$}\\
First, we make the change of variables $\eta=\nabla_x\phase(x,\xi)+\zeta$ in the integral defining $\textbf{I}_2 (t,x,\xi)$ and then expand $b(x, t\eta)$ in a Taylor series to obtain
\eq{
    &b(x,t\nabla_x\phase(x,\xi)+t\zeta)  = \sum_{0 \leq |\alpha|<M} t^{|\alpha|}\frac{\zeta^\alpha}{\alpha !} (\partial_\eta^\alpha b))(x,t\nabla_x\phase(x,\xi)) \\ & \qquad+ t^{M} \sum_{|\alpha|=M} C_\alpha {\zeta^\alpha} r_\alpha(t, x,\xi,\zeta),
}
 where 
\begin{equation}\label{eq:ralapha}
    r_\alpha(t, x,\xi,\zeta) 
 := \int_0^1 (1-\tau)^{M-1} (\partial_\eta^{\alpha} b)(x, t\nabla_x\phase(x,\xi)+\tau t\zeta)\dd \tau.
\end{equation}
If we set
\[
\Phi(x,y,\xi)
:=\phase(y,\xi)-\phase(x,\xi)+(x-y)\cdot\nabla_x\phase(x,\xi),
\]
we obtain
\[
\textbf{I}_2 (t,x,\xi)= \sum_{|\alpha|<M} \frac{t^{\eps|\alpha|}}{\alpha!}\, \sigma_{\alpha}(t,x,\xi) + t^{\eps M}\, \sum_{|\alpha|=M} C_\alpha\, R_{\alpha}(t,x,\xi),
\]
where using integration by parts, we have
\begin{align*}
&\sigma_{\alpha}(t,x,\xi) 
 := t^{(1-\eps)|\alpha|}\iint_{\Rl^n\times \Rl^n} e^{i(x-y)\cdot\zeta+i\Phi(x,y,\xi)} \,\zeta^{\alpha}\, a(y,\xi) \\
 &\qquad\qquad\qquad\times\chi(x-y)\, (\partial_\eta^\alpha b)(x, t\nabla_x\phase(x,\xi))\dd y\ddd\zeta \\
& =t^{(1-\eps)|\alpha|}(\partial_\eta^\alpha b)(x,t\nabla_x\phase(x,\xi)) (i)^{-|\alpha|}\partial_y^{\alpha}\left[ e^{i\Phi(x,y,\xi)}\, a(y,\xi)\,\chi(x-y) \right]_{|_{_{y=x}}},
\end{align*}
and
\[
R_{\alpha}(t,x,\xi) 
:= t^{(1-\eps)|\alpha|}\iint_{\Rl^n\times \Rl^n} e^{i(x-y)\cdot\zeta} e^{i\Phi(x,y,\xi)} \zeta^{\alpha}\,a(y,\xi)\, \chi(x-y) \, r_\alpha(t, x,\xi,\zeta)\, \dd y\ddd\zeta.
\]

\quad\\

\textbf{Step 2.1 -- The analysis of $\mathbf{\boldsymbol\sigma_{\boldsymbol\alpha}(t,x,\boldsymbol\xi)}$}\\
We now claim that
\begin{equation} \label{half}
\left|\partial_y^{\gamma} e^{i\Phi(x,y,\xi)}{}_{|_{y=x}}\right|
\lesssim \bra{\xi}^{|\gamma|/2}.
\end{equation}
We first observe that when $\gamma = 0$, \eqref{half} is obvious. To obtain \eqref{half} for $\gamma \neq 0$ we recall Fa\`a di Bruno's formulae
\[
\partial_y^{\gamma} e^{i\Phi(x,y,\xi)}=\sum_{\gamma_1 + \cdots+ \gamma_k =\gamma} C_\gamma \brkt{\partial^{\gamma_{1}}_{y}\Phi(x,y,\xi)}\cdots \brkt{\partial^{\gamma_{k}}_{y}\Phi(x,y,\xi)}\,e^{i\Phi(x,y,\xi)},
\]
where the sum ranges of $\gamma_j$ such that $|\gamma_{j}|\geq 1$ for $j=1,2,\dots, k$ and $\gamma_1 + \cdots+ \gamma_k =\gamma$ for some $k \in \Z_+$. Since $\Phi(x,x,\xi)=0$ and $\left.\partial_y\Phi(x,y,\xi)\right|_{y=x}=0$, setting $y=x$ in the expansion above leaves only terms in which $|\gamma_j|\geq 2$ for all $j = 1,2,\dots,k$. But $\sum_{j=1}^{k} |\gamma_j |\leq |\gamma|,$ so we actually have $2k\leq |\gamma|$, that is $k\leq  |\gamma|/2$. Estimate \eqref{eq:composition_conditions} on the phase tells us that $|\partial^{\gamma_{j}}_{y}\Phi(x,y,\xi)| \lesssim \xxi$, so
\[
\left|\partial_y^{\gamma} e^{i\Phi(x,y,\xi)}{}_{|_{y=x}}\right| \lesssim \xxi \cdots \xxi \lesssim \xxi ^{k} \lesssim \bra{\xi}^{|\gamma|/2},
\]
which is \eqref{half}.\\

If we use the fact that $t\leq 1$ and the assumption \eqref{eq:composition_conditions} on the phase function $\phase$, then we have
\eq{
&|\sigma_{\alpha}(t,x,\xi)| \lesssim t^{(1-\eps)|\alpha|}\bra{t\nabla_x\phase(x,\xi)}^{s-|\alpha|} \bra{\xi}^{|\alpha|\max\{\delta,1/2\}} \bra{\xi}^m
\\ &\lesssim  t^{(1-\eps)|\alpha|} \,\bra{t\xi}^{s-(1-\varepsilon)|\alpha|}\,\bra{t\xi}^{-\varepsilon|\alpha|}\, \bra{\xi}^{m+|\alpha|\max\{\delta,1/2\}}
\\ &\lesssim t^{\min\{s,0\}}  \bra{\xi}^{s+m-(1-\max\{\delta,1/2\}-\varepsilon)|\alpha|},
}
when $|\alpha| > 0$.\\ 

By the assumptions of the theorem, the derivatives of $\sigma_{\alpha}$ with respect to $x$ or $\xi$ do not change the estimates when applied to $b$, and the same is true when derivatives are applied to $\partial_y^{\alpha} e^{i\Phi(x,y,\xi)}|_{y=x}$. Therefore, for all multi-indices $\beta$, $\gamma\in \mathbb{Z}_{+}$,
\[
\abs{\partial^{\beta}_{\xi} \partial^{\gamma}_{x} \sigma_{\alpha}(t,x,\xi)} \lesssim  t^{\min\{s,0\}}  \bra{\xi}^{s+m-(1-\max\{\delta,1/2\}-\varepsilon)|\alpha|-\rho|\beta|+\delta|\gamma|},
\]
as required.\\

\textbf{Step 2.2 -- The analysis of $\mathbf{R_{\boldsymbol\alpha}(t,x,\boldsymbol\xi)}$}\\
Take $g\in \mathcal{C}_c^\infty(\Rl^n)$ such that $g(x)=1$ for $|x|<\kappa_1/2$ and $g(x)=0$ for $|x|>\kappa_1$, where $\kappa_1>0$ will be chosen later. We then decompose
\begin{equation*}\label{EQ:Ralphas}
 \begin{aligned} 
& R_{\alpha}(t,x,\xi) 
= t^{(1-\eps)|\alpha|}\iint_{\Rl^n\times \Rl^n} e^{i(x-y)\cdot\zeta} g\Big(\frac{\zeta}{\bra{\xi}}\Big) \\
&\qquad\qquad\qquad\times \partial_y^{\alpha} \left[ e^{i\Phi(x,y,\xi)}\, \chi(x-y)\,\,a(y,\xi)\, r_\alpha(t,x,\xi,\zeta) \right] \dd y\ddd\zeta \\ 
& \qquad + t^{(1-\eps)|\alpha|} \iint_{\Rl^n\times \Rl^n} e^{i(x-y)\cdot\zeta} \Big(1-g\Big( \frac{\zeta}{\bra{\xi}}\Big)\Big) \\
&\qquad\qquad\qquad\times \partial_y^{\alpha}\left[  e^{i\Phi(x,y,\xi)}\,\chi(x-y)\,\,a(y,\xi)\, r_\alpha(t,x,\xi,\zeta)\right]\dd y \ddd\zeta \\
&=: R_\alpha^I(t,x,\xi) + R_\alpha^{I\!\!I}(t,x,\xi). 
\end{aligned}
\end{equation*}

\quad\\

\textbf{Step 2.2.1 -- The analysis of $\mathbf{R^I_{\boldsymbol\alpha}(t,x,\boldsymbol\xi)}$}\\
Note that the inequality
\[
\bra{\xi}\leq 1+|\xi|\leq \sqrt{2}\bra{\xi},
\]
and \eqref{eq:composition_conditions} yield
$$
\bra{t\nabla_x\phase(x,\xi)+t\tau\zeta} 
\leq  (C_2\sqrt{2}+\kappa_1)\bra{t\xi},  $$
and
\begin{align*}
& \sqrt{2}\bra{t\nabla_x\phase(x,\xi)+t\tau\zeta} \geq 1+|t\nabla_x\phase|-|t\zeta|  \\
& \geq 1+C_1|t\xi|-t\kappa_1\bra{\xi} \\
& \geq (1-\kappa_1)+(C_1-\kappa_1)|t\xi| \geq (\min\{1,C_1\}-\kappa_1)\bra{t\xi}.
\end{align*}

Therefore, if we choose $\kappa_1<\min\{1,C_1\}$, then for any $\tau\in (0,1)$, $\bra{t\nabla_x\phase(x,\xi)+t\tau\zeta}$ and $\bra{t\xi}$ are equivalent.\\

This yields that for $|\zeta|\leq r\bra{\xi}$, $\partial^\beta_\zeta r_\alpha(t, x,\xi,\zeta)$  are dominated by $ t^{|\beta|}\bra{t\xi}^{s-|\alpha|-|\beta|}.$ Furthermore, for $t\leq 1$, it follows from the representation \eqref{eq:ralapha} for $r_\alpha$ that
\begin{equation}\label{eq:estr}
\begin{aligned}
& \Big|\partial_\zeta^{\beta}\Big( g\Big( \frac{\zeta}{\bra{\xi}}\Big) r_\alpha(t,x,\xi,\zeta)\Big)\Big| 
 \lesssim\sum_{\gamma\leq\beta} \Big|\partial_\zeta^\gamma g\Big(\frac{\zeta}{\bra{\xi}}\Big) \partial_\zeta^{\beta-\gamma} r_\alpha(t,x,\xi,\zeta)\Big| \\
& \leq C_{\alpha,\beta} \sum_{\gamma\leq\beta}t^{|\beta|-|\gamma|}\bra{\xi}^{-|\gamma|} \bra{t\xi}^{s-|\alpha|-|\beta|+|\gamma|} \\
&  \lesssim \sum_{\gamma\leq\beta}t^{\min\{s,0\}+|\beta|-|\gamma|-(1-\varepsilon)|\alpha|}\bra{\xi}^{-|\gamma|}\\ 
&\qquad \qquad \qquad \times \bra{\xi}^{s-(1-\varepsilon)|\alpha|}   t^{-(|\beta|-|\gamma|)}\bra{\xi}^{-(|\beta|-|\gamma|)}\\
&  \lesssim t^{\min\{s,0\}-(1-\varepsilon)|\alpha|} \,\bra{\xi}^{s-(1-\varepsilon)|\alpha|-|\beta|}.
\end{aligned}
\end{equation}

At this point we also need estimates for $\partial_y^\alpha e^{i\Phi(x,y,\xi)}$ off the diagonal, that is, when $x\neq y.$ This derivative has at most $|\alpha|$ powers of terms $\nabla_y\phase(y,\xi)-\nabla_x\phase(x,\xi)$, possibly also multiplied by at most $|\alpha|$ higher order derivatives $\partial_y^\beta\phase(y,\xi)$, which can be estimated by $(|y-x|\bra{\xi})^{|\alpha|}$ using \eqref{eq:composition_conditions}. The term containing the difference $\nabla_y\phase(y,\xi)-\nabla_x\phase(x,\xi)$ is the product of at most $|\alpha|/2$ terms of the type $\partial_y^\beta\phase(y,\xi)$, which can be estimated by $\bra{\xi}^{|\alpha|/2}$ in view of \eqref{eq:composition_conditions}. These observations yield
\[
\abs{\partial_y^\alpha e^{i\Phi(x,y,\xi)}}\lesssim (1+|x-y|\bra{\xi})^{|\alpha|} \bra{\xi}^{|\alpha|/2},
\]
and therefore we also have
\begin{equation}\label{eq:ests}
\left|\partial_y^{\alpha}\left[ e^{i\Phi(x,y,\xi)} \chi(x-y) \right] \right|\lesssim (1+|x-y|\bra{\xi})^{|\alpha|}\bra{\xi}^{|\alpha|/2}.
\end{equation}

Let
\[
L_\zeta
:=\frac{(1-\bra{\xi}^2\Delta_\zeta)}{1+\bra{\xi}^2|x-y|^2}, \quad \text{so} \quad L_\zeta^N e^{i(x-y)\cdot\zeta}=e^{i(x-y)\cdot\zeta}.
\]
Integration by parts with $L_\zeta$ yields
\[
 \begin{aligned}
& R^I_{\alpha}(t,x,\xi) 
= t^{(1-\eps)|\alpha|}\iint_{\Rl^n\times \Rl^n} \frac{e^{i(x-y)\cdot\zeta}\,\partial_y^{\alpha}\left[ \chi(x-y)\, a(y,\xi)\, e^{i\Phi(x,y,\xi)}\right]}{(1+\bra{\xi}^2 |x-y|^2)^N} \\
&\qquad \qquad \times (1-\bra{\xi}^{2}\Delta_\zeta)^N \Big\{ g\Big(\frac{\zeta}{\bra{\xi}}\Big) \, r_\alpha(t,x,\xi,\zeta) \Big\} \dd y \ddd\zeta \\
&  = t^{(1-\eps)|\alpha|}\iint_{\Rl^n\times \Rl^n} \frac{e^{i(x-y)\cdot\zeta}\,\partial_y^{\alpha}\left[\chi(x-y)\,a(y,\xi)\,e^{i\Phi(x,y,\xi)}\right]}{(1+\bra{\xi}^2 |x-y|^2)^N}\\
&\qquad \qquad \times \sum_{|\beta|\leq 2N} c_{\beta}\bra{\xi}^{|\beta|} \Big\{ \partial_\zeta^{\beta}\Big( g\Big(\frac{\zeta}{\bra{\xi}}\Big) r_\alpha(t,x,\xi,\zeta)\Big) \Big\} \dd y \ddd\zeta.
\end{aligned}
\]
Using estimates (\ref{eq:estr}), (\ref{eq:ests}) and that the size of the support of $g(\zeta/\bra{\xi})$ in $\zeta$ is bounded by $(\kappa_1\bra{\xi})^n$, we obtain
\[
\begin{aligned}
&|R^I_{\alpha}(t,x,\xi)|  \lesssim t^{\min\{s,0\}}\,\sum_{|\beta|\leq 2N} \bra{\xi}^{n+|\beta|} \bra{\xi}^{-(1-\varepsilon)|\alpha|-|\beta|} \bra{\xi}^{\max\{\delta,1/2\}+s+m} \\
&\qquad \qquad \qquad\times\int_{|x-y|<\kappa}\frac{(1+|x-y|\bra{\xi})^{|\alpha|}} {(1+\bra{\xi}^2 |x-y|^2)^N} \dd y \\ 
& \lesssim t^{\min\{s,0\}}\sum_{|\beta|\leq 2N} \bra{\xi}^{n+|\beta|} \bra{\xi}^{s-(1-\varepsilon)|\alpha|-|\beta|} \bra{\xi}^{\max\{\delta,1/2\}+m} \\
&\qquad \qquad \qquad\times \langle \xi\rangle^{-n} \int_{0}^{\infty}\frac{\tau^{n-1}(1+\tau)^{|\alpha|}} {(1+\tau^2)^N} \dd \tau\\
& \lesssim t^{\min\{s,0\}} \bra{\xi}^{s+m-(1-\max\{\delta,1/2\}- \varepsilon)|\alpha|},
\end{aligned}
\]
if we choose $N>(n+|\alpha|)/2$,  and the hidden constants in the estimates are independent of $t$ (because of (\ref{eq:estr})). The derivatives of $R_{\alpha}^I(t,x,\xi)$ with respect to $x$ and $\xi$ give an extra power of $\zeta$ under the integral. This amounts to taking more $y$-derivatives, yielding a higher power of $\bra{\xi}.$ However, for a given number of derivatives of the remainder $R_{\alpha}^I(t,x,\xi)$, we are free to choose $M=|\alpha|$ as large as we like and therefore the higher power of $\bra{\xi}$ will not cause a problem. Thus for all multi-indices $\beta$, $\gamma,$ and $|\alpha|$ large enough we have 
\begin{equation*}
|\partial^{\beta}_{\xi} \partial^{\gamma}_{x} R_{\alpha}^{I}(t,x,\xi)| \lesssim  t^{\min\{s,0\}} \bra{\xi}^{s+m-(1-\max\{\delta,1/2\}- \varepsilon}|\alpha|-\rho|\beta|),
\end{equation*}
where the hidden constant in the estimate does not depend on $t$.\\

\textbf{Step 2.2.2 -- The analysis of $\mathbf{R_{\boldsymbol\alpha}^{I\!I}(t,x,\boldsymbol\xi)}$}\\
Define
\[
\Psi(x,y,\xi,\zeta)
:=(x-y)\cdot\zeta+\Phi(x,y,\xi)= (x-y)\cdot(\nabla_x\phase(x,\xi)+\zeta)+\phase(y,\xi)-\phase(x,\xi).
\]
It follows from \eqref{eq:composition_conditions} that if we choose $\kappa_2 <\kappa_1/8C_0,$ then since $|x-y|<\kappa_2$ on the support of $\chi$, one has (using that we are in the region $|\zeta|\geq \kappa_1\bra{\xi}/2$)
\begin{equation*}\label{eq:rho}
\begin{aligned}
|\nabla_y\Psi| & =|-\zeta+\nabla_y\phase-\nabla_x\phase|\leq 2C_2(|\zeta|+\bra{\xi}), \quad \text{and} \\
|\nabla_y\Psi| & \geq |\zeta|-|\nabla_y\phase-\nabla_x\phase| \geq \frac{1}{2}|\zeta|+\Big(\frac{\kappa_1}{4}-C_0|x-y| \Big)\bra{\xi}\geq C(|\zeta|+\bra{\xi}).
\end{aligned}
\end{equation*}
Now, using \eqref{eq:composition_conditions}, for any $\beta$ we have the estimate
\begin{equation}\label{no idea what to call 1}
\abs{\partial_y^\beta\brkt{e^{-i\Phi(x,y,\xi)}\,\d_y^{\gamma} e^{i\Phi(x,y,\xi)}}}\lesssim\bra{\xi}^{|\gamma|}.
\end{equation}
For $M=|\alpha|>s$ we also observe that
\begin{equation} \label{eq:rs}
|r_\alpha(t, x,\xi,\zeta)|\lesssim 1.
\end{equation}
For the differential operator defined to be 
$$L_y:=i|\nabla_y\Psi|^{-2}\sum_{j=1}^n (\partial_{y_j}\Psi) \,\partial_{y_j},$$ 
induction shows that $L_y^N$ has the form
\begin{equation*} (L_y^*)^N=\frac{1}{|\nabla_y\Psi|^{4N}}\sum_{|\beta|\leq N} P_{\beta,N}\,\partial_y^\beta, 
\end{equation*}
where
$$ P_{\beta,N}
:=\sum_{|\mu|=2N} c_{\beta\mu\gamma_j}\,(\nabla_y\Psi)^\mu\, \partial_y^{\gamma_1}\Psi\cdots \partial_y^{\gamma_N}\Psi,$$
$|\gamma_j|\geq 1$ and $\sum_{j=M}^{N}|\gamma_j|+|\beta|=2N$. It follows from \eqref{eq:composition_conditions} that $|P_{\beta,N}|\leq C(|\zeta|+\bra{\xi})^{3N}.$ Now Leibniz's rule yields
\begin{equation*}
\begin{aligned}
& R^{I\!\!I}_{\alpha}(t,x,\xi) 
= t^{(1-\eps)|\alpha|}\iint_{\Rl^n\times \Rl^n} e^{i(x-y)\cdot\zeta} \,
\Big(1-g\Big( \frac{\zeta}{\bra{\xi}}\Big)\Big)\, r_\alpha(x,\xi,\zeta) \\
& \qquad \qquad \qquad\times\partial_y^{\alpha}\left[ e^{i\Phi(x,y,\xi)}\, a(y,\xi)\,\chi(x-y) \right] \dd y \ddd\zeta \\
& = t^{(1-\eps)|\alpha|}\iint_{\Rl^n\times \Rl^n} e^{i\Psi(x,y,\xi,\zeta)} \,
\Big(1-g\Big(\frac{\zeta}{\bra{\xi}}\Big)\Big)\, r_\alpha(t,x,\xi,\zeta) \\
& \qquad\qquad \qquad\times \sum_{\gamma_1+\gamma_2 +\gamma_3=\alpha} \brkt{e^{-i\Phi(x,y,\xi)}\partial_y^{\gamma_1} e^{i\Phi(x,y,\xi)}}\, \partial^{\gamma_2}_{y}\chi(x-y) \,\partial^{\gamma_3}_{y} a(y,\xi)  \dd y \ddd\zeta \\
& = t^{(1-\eps)|\alpha|}\iint_{\Rl^n\times \Rl^n} e^{i\Psi(x,y,\xi,\zeta)} |\nabla_y\Psi|^{-4N}\sum_{|\beta|\leq N} P_{\beta,N}(x,y,\xi,\zeta)\\
& \qquad\qquad\qquad \times \Big(1-g\Big(\frac{\zeta}{\bra{\xi}} \Big)\Big) r_\alpha(t,x,\xi,\zeta)  \sum_{\gamma_1+\gamma_2 +\gamma_3=\alpha} \partial_y^\beta \big [ \brkt{e^{-i\Phi(x,y,\xi)}}\\
& \qquad\qquad \qquad\times \partial_y^{\gamma_1} e^{i\Phi(x,y,\xi)}\, \partial^{\gamma_2}_{y}\chi(x-y)\, \partial^{\gamma_3}_{y} a(y,\xi)\big ] \dd y \ddd\zeta.
\end{aligned}
\end{equation*}
It follows now from \eqref{no idea what to call 1} and (\ref{eq:rs}) that
\[
\begin{aligned}
& |R^{I\!\!I}_{\alpha}(t,x,\xi)|  \lesssim t^{(1-\eps)|\alpha|}\int_{|\zeta|\geq \kappa_1\bra{\xi}/2} \int_{|x-y|<\kappa} (|\zeta|+\bra{\xi})^{-N} \bra{\xi}^{|\alpha|+m}\dd y \ddd\zeta\\ 
&\lesssim t^{(1-\eps)|\alpha|}\bra{\xi}^{|\alpha|+m} 
\int_{|\zeta|\geq \kappa_1\bra{\xi}/2}|\zeta|^{-N}\ddd\zeta  \leq C \bra{\xi}^{|\alpha|+n+m-N},
\end{aligned}
\]
which yields the desired estimate when $N>|\alpha|+n$. For the derivatives of $R_{\alpha}^{I\!\!I}(t,x,\xi)$, we can get, in a similar way to the case for $R_{\alpha}^I$, an extra power of $\zeta$, which can be taken care of by choosing $N$ large and using the fact that $|x-y|<\kappa_2.$  Therefore for all multi-indices $\beta$, $\gamma\in \mathbb{Z}_{+},$
\begin{equation*}
            \abs{\partial^{\beta}_{\xi} \partial^{\gamma}_{x} R_{\alpha}^{I\!\!I}(t,x,\xi)} \lesssim \bra{\xi}^{|\alpha|+n+m-N},
\end{equation*}
where the constant hidden in the estimate does not depend on $t$. The proof of Theorem \ref{thm:left composition with pseudo} is now complete.
\end{proof}

\section{Decomposition of FIOs of constant rank}\label{SSS decomposition}

In connection to the study of the $L^p$-regularity of FIOs, based on an idea of C. Fefferman \cite{Feff}, Seeger, Sogge and Stein \cite{SSS} introduced a second dyadic decomposition superimposed on a preliminary Littlewood-Paley decomposition.\\

Note that, since we are dealing with FIOs with amplitudes in general H\"ormander classes, the tools at hand have to be adapted to include the treatments of amplitudes of type $\rho\in[0,\frac{1}{2})$ which was missing in the analysis given in \cite{SSS}. Here we follow the decomposition of the frequency space, as in Dos Santos Ferreira-Staubach \cite{DS}.\\

To start, one considers an FIO $T^{\varphi}_{a}$ with amplitude $a(x, \xi)\in S^{m}_{\rho, \delta}(\Rl^n)$, $0<\rho\leq \frac23$, $0\leq \delta<1$, $m\in\Rl$ and that for $0 \leq \kappa\leq n-1$, $\varphi\in\Phi^2$,  $\mathrm{rank}\,\partial^{2}_{\xi\xi} \varphi(x, \xi) = \kappa$, on the support of $a(x, \xi)$ and its Littlewood-Paley decomposition
\begin{equation}\label{eq:LPdecomp}
T^{\varphi}_{a}= \sum_{j=0}^\infty  T^{\varphi}_{a}\psi_j(D) =: \sum_{j=0}^\infty  T_j, 
\end{equation}
where the kernel $K_j$ of $T_j$ is given by
\nm{eq:LPdecompkernel}{
    &K_j(x,y):= \int_{\Rn} e^{i\varphi(x,\xi)-iy\cdot\xi}\,\psi_j(\xi)\,a(x,\xi)\, \ddd \xi \\
    &=2^{jn\rho}\int_{\Rn} e^{i2^{j\rho}\varphi(x,\xi)-i2^{j\rho}y\cdot\xi}\,\psi_j(2^{j\rho}\xi)\,a(x,2^{j\rho}\xi) \ddd \xi \nonumber
}
Here each $\psi_j$ is supported in a dyadic shell $\left \{2^{j-1}\leq \vert \xi\vert\leq 2^{j+1}\right \}$ (as in Definition \ref{def:LP}) so the $\xi$-support of the integrand in the second integral is 
\nm{eq:Aj}{A_j:=\{\xi\in\Rn:2^{j(1-\rho)-1}\leq|\xi|\leq 2^{j(1-\rho)+1}\}.}

The shells $A_j$ will in turn be decomposed into truncated cones using the following construction:\\

Since $\varphi$ has constant rank $\kappa$ on the support of $a(x, \xi)$ we may assume that there exists some $\kappa$-dimensional submanifold $S_\kappa(x)$ of $\mathbb S^{n-1}\cap \Gamma$ for some sufficiently narrow cone $\Gamma$, such that $\mathbb S^{n-1}\cap \Gamma$ is parameterised by $\bar{\xi}=\bar{\xi}_x(u,v)$, for $(u,v)$ in a bounded open set $U\times V$ near $(0,0)\in \Rl^{\kappa}\times\Rl^{n-\kappa-1}$, and such that $\bar{\xi}_x(u,v)\in S_\kappa(x)$ if and only if $v=0$, and $\nabla_\xi\varphi(x,\bar{\xi}_x(u,v))=\nabla_\xi\varphi(x,\bar{\xi}_x(u,0))$.

\begin{Def}\label{def:LP2}
For each $j\in\mathbb Z_{>0}$, $0<\rho\leq 1$ and $0\leq \kappa\leq n-1$, let $\{u^{\nu}_{j} \} $ be a collection of points in $U$ such that,
\begin{enumerate}
    \item [$(i)$]$  | u^{\nu}_{j}-u^{\nu'}_{j} |\geq 2^{-j\rho/2},$ whenever $\nu\neq \nu '$.
    \item [$(ii)$] If $u\in U$, then there exists a $
    u^{\nu}_{j}$ so that $\vert u -u^{\nu}_{j}  \vert
    <2^{-j\rho/2}$.
\end{enumerate}
Set $\xi^{\nu}_{j}=\bar{\xi}_x(u^{\nu}_{j},0)$.
One may take such a sequence by choosing a maximal collection $\{\xi_{j}^{\nu}\}$ for which (i) holds, then (ii) follows. Now, denote the number of cones needed by $\mathscr N_j$. See $\mathrm{Figure}$ $\ref{fig:submanisphere}$ for an illustration of the case $\mathbb S^2$ and $\kappa=1$.\\

Let $\Gamma^{\nu}_{j}$ denote the cone in the $\xi$-space, with the apex at the origin, whose
 central direction is $\xi^{\nu}_{j}$, i.e.
\begin{equation}\label{eq:gammajnu}
\Gamma^{\nu}_{j}
:= \set{ \xi \in \Rl^n ;\, \xi=s\bar{\xi}_x(u,v), \abs{ \frac{u}{|u|} - u^{\nu}_{j} }
    <2^{-j\rho/2}, v\in V, s>0}
\end{equation}
One also defines the following partition of unity on $\Gamma$. To construct this, take
\eq{\tilde \eta_j^\nu (u):= \phi\brkt{2^{j\rho/2}\brkt{\frac{u}{|u|} - u^{\nu}_{j}}}}
where $\phi$ is a non-negative function in $ \mathcal C_c^\infty(\Rl^n)$ with $\phi(u)=1$ for $|u|\leq 1$ and $\phi(u)=0$ for $|u|\geq 2$, in such a way that 
$$\tilde\chi_j^\nu(u) := \frac{\tilde\eta_j^\nu(u)}{ \sum_{\nu=1}^{\mathscr N_j} \tilde\eta_j^\nu(u)}.$$
The corresponding homogeneous partition of unity $\chi_j^\nu$ is defined such that
\nm{eq:defchijnu}{\chi_j^\nu (s\bar{\xi}_x(u,v))=\tilde\chi_j^\nu (u),}
for $v\in V$ and $s>0$.
\end{Def}
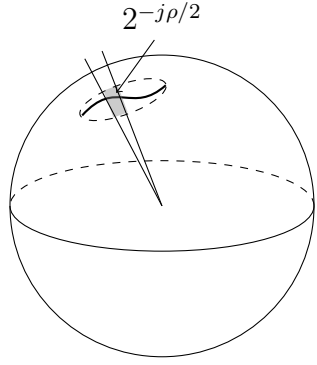
\begin{figure}
\begin{tikzpicture}
    \draw (0,0) circle (2cm);
    \draw (-2,0) arc (180:360:2 and 0.6);
    \draw[dashed] (2,0) arc (0:180:2 and 0.6);
    \draw[dashed,rotate=20] (0.0,1.5) ellipse (0.6 and 0.2);
    \draw[thick] (-1.05,1.2) .. controls (-0.6,1.7) and (-0.4,1.2) .. (0.05, 1.6);
    \fill[gray, opacity = 0.4, rotate = 20] (-0.2,1.32) rectangle (-0.0,1.7);
    \draw (0,0) -- (111:2.2);
    \draw (0,0) -- (117.5:2.2);
    \draw (0,2.2) node[anchor=south] 
    {$2^{-j\rho/2}$};
    \draw[->] (-0.1,2.2) -- (-0.6,1.5);
\end{tikzpicture}
\caption{Illustration of the case $\mathbb S^2$ and $\kappa=1$. The thick line illustrates the submanifold $S_1(x)$ in the open set $\bar{\xi}_x(U\times V)$. The length of the curve $S_1(x)$ within the grey area doesn't exceed $2^{-j\rho/2}$. }\label{fig:submanisphere}
\end{figure}

Using Definition \ref{def:LP2}, we can make a second dyadic decomposition of the kernel \eqref{eq:LPdecompkernel} as 
\nm{eq:LPdecompkernel2}{
    K_j^\nu(x,y):= 2^{jn\rho}\int_{\Rn} e^{i2^{j\rho}\varphi(x,\xi)-i2^{j\rho}y\cdot\xi}\,\psi_j(2^{j\rho}\xi)\,\chi_j^\nu(\xi)\,a(x,2^{j\rho}\xi) \ddd \xi 
}
and observe that 
$$
    K_j(x,y) = \sum_{\nu=1}^{\mathscr N_j} K_j^\nu(x,y).
$$

Now turning back to the shells in \eqref{eq:Aj}, we decompose each $A_j$ into truncated cones $\Gamma^{\nu}_{j}\cap A_j$.  We claim that $\mathscr N_j= O\big(2^{j\rho\kappa/2}\big)$
and that each truncated cone has volume of the size $O\big(2^{jn(1-\rho)-j\rho\kappa/2})$ 
\\

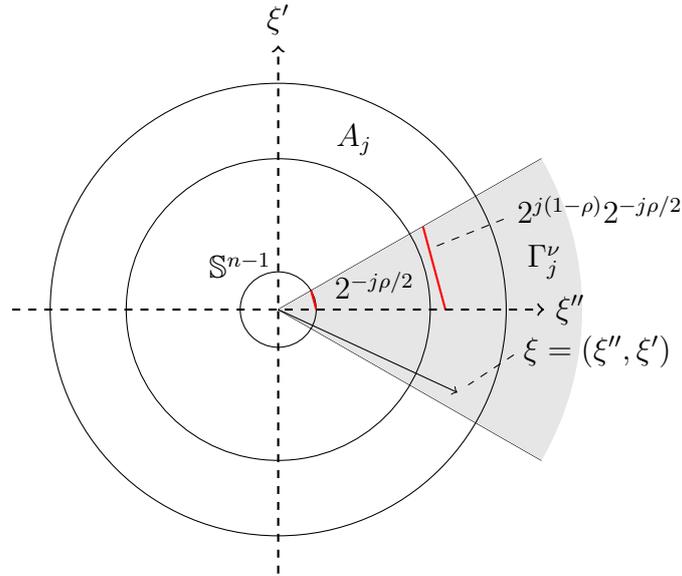
\begin{figure}
    \begin{tikzpicture}

        \draw[rotate = 30] (0,0) -- (4, 0); 
        \draw[rotate = -30] (0,0) -- (4, 0);
        \fill[gray!20] (0,0) -- (30:4) arc(30:-30:4) -- cycle ;

        \draw (0,0) circle (0.5);
        \draw (0,0) circle (2);
        \draw (0,0) circle (3);
    
        \draw (-0.5,0.9) node[anchor=north] {$\mathbb S^{n-1}$};
        \draw (1, 1.9) node[anchor=south] {$A_j$};
        \draw (3.5, 1) node[anchor=north] {$\Gamma^\nu_j$};
        \draw (0.6,0) node[anchor=south west] {$2^{-j\rho/2}$};
        \draw (3,1) node[anchor=south west] {$2^{j(1-\rho)}2^{-j\rho/2}$};
        \draw (4.2,-0.2) node[anchor=north] {${\xi= (\xi'',\xi')}$};
        
        \draw[red,thick] (30:0.5) -- (0.5,0);
        \draw[red,thick] (30:2.2) -- (2.2,0);
    
        \draw[->,rotate = -25] (0,0) -- (2.6,0);
        \draw[dashed] (3.1,-0.6) -- (2.5,-1);
        \draw[dashed] (22:2.25) -- (22:3.25);
        \draw[->,dashed,thick] (-3.5,0)--(3.5,0) node[right]{$\xi''$};
        \draw[->,dashed,thick] (0,-3.5)--(0,3.5) node[above]{$\xi'$};

    \end{tikzpicture}
\caption{An illustration of the intersection between $A_j$ and $\Gamma_j^\nu$, containing frequencies $\xi$ of size $2^{j(1-\rho)}$.}\label{fig:conesize}
\end{figure}

Indeed observe that by the definition \eqref{eq:gammajnu} we have that for each truncated cone $\Gamma^{\nu}_{j}\cap A_j$ there are $\kappa$ directions with length (roughly) equal to $2^{j(1-\rho)} 2^{-j \rho / 2} $ and $n-\kappa-1$ directions with length roughly equal to $2^{j(1-\rho)}$ (which is the thickness of $A_j$). Hence we infer that 
\begin{equation}\label{eq:measureintersect}
    |\Gamma^{\nu}_{j}\cap A_j|
\sim 
    \left(2^{-j \rho / 2} 2^{j(1-\rho)}\right)^{\kappa} 2^{j(n-\kappa-1)(1-\rho)} 2^{j(1-\rho)}
=
    2^{jn(1-\rho)-j\rho\kappa/2}.
\end{equation}
Using this, it also follows that there are roughly
$$\frac{2^{jn(1-\rho)}}{2^{jn(1-\rho)-j\rho\kappa/2}} = 2^{j\rho\kappa/2}$$
such truncated cones needed to cover one shell $A_j$, which proves the claim. See {\bf Figure \ref{fig:conesize}} for an illustration.\\ 

For the following lemma and throughout the rest of the  paper, we choose the coordinate axes in $\xi$-space such that $\xi''\in T_{\xi^{\nu}_j}S_\kappa(x)$ and $\xi'\in\big(T_{\xi^{\nu}_j}S_\kappa(x)\big)^\perp$, so that $\xi=(\xi'',\xi')$. Observe that since $\frac{|\xi'|}{|\xi|}\lesssim 2^{-j\rho/2}$ one has that
\nm{eq:xiprimestimate}{
    |\xi'|\lesssim 2^{j(1-\rho)}2^{-j \rho / 2}=2^{j(1-3\rho/2)},
}
for $\xi \in \Gamma_j^\nu\cap A_j.$

\begin{Lem}\label{lem:chijnu}
The functions $\chi_j^\nu$ belong to $\mathcal C^\infty(\Rl^n\setminus \{0\} )$ and are supported in the cones $\Gamma_j^\nu$. They sum to $1$ in $\nu$\emph:
\begin{equation*}
\sum_{\nu=1}^{\mathscr N_j}\chi^{\nu}_{j}(\xi) =1,\quad j\in \mathbb{N}\,\, \text { and } \xi\neq 0 
\end{equation*}
and moreover they satisfy the estimates 
\nm{eq:quadrseconddyadcond}{|\d_\xi^\alpha \chi_j^\nu (\xi) |\lesssim 
2^{j\rho|\alpha|/2}
|\xi|^{-|\alpha|},}
for all multi-indices $\alpha$ and
\nm{eq:quadrseconddyadcond2}{ 
    \vert \partial^{\beta}_{\xi''}\chi^{\nu}_{j}(\xi)  \vert\lesssim\vert \xi\vert ^{-|\beta|},
}
for all $\beta \in \mathbb Z^{n-\kappa}_{\geq0}$.
\end{Lem}

\begin{proof}
Note that by \eqref{eq:defchijnu} it is enough to show the result for $\tilde\chi_j^\nu$, since $\xi_x$ is a diffeomorphism. In proving \eqref{eq:quadrseconddyadcond} we note that the argument of $\eta_j^\nu$ contains a factor of $2^{j\rho/2}$ followed by a factor that is homogeneous of degree zero. Hence $\alpha$ derivatives yield a factor of $2^{j\rho|\alpha|/2}$ and a function that is homogeneous of degree $-|\alpha|$. 
To prove \eqref{eq:quadrseconddyadcond2} one observes that on the support of $\chi^\nu_j$ one can write $$\partial_{\xi''_\ell}= \partial_{r}+O(2^{-j\rho/2}) \cdot\nabla_\xi,$$  for any $1\leq\ell\leq n-\kappa,$ where $\partial_{r}$ is the radial derivative and $\partial_{r} \tilde\chi^\nu_j=0$  since $\tilde\chi^\nu_j$ is homogeneous of degree zero.
\end{proof}

We will split the phase $\varphi(x,\xi)-y\cdot\xi$ into two different pieces,  $(\nabla_\xi\varphi(x,\xi_j^\nu)-y)\cdot \xi$ (which is linear in $\xi$), and $\varphi(x,\xi) - \nabla_\xi \varphi(x,\xi_j^\nu) \cdot \xi$. The following lemma yields an estimate for the nonlinear second piece.

\begin{Lem}\label{lem:h}
Let $\varphi(x,\xi)\in\mathcal{C}^{\infty}(\Rl^n \times \Rl^n\setminus\{0\})$ be positively homogeneous of degree one in $\xi$. For $j\in \Z_{\geq0}$, $1\leq\nu \leq \mathscr N_j,$ define
$$h_j^\nu(x,\xi):= \varphi(x,\xi) -\nabla_\xi \varphi(x,\xi_j^\nu)\cdot \xi.$$
Then for $0<\rho\leq \frac23$ and for  $\xi$ in $A_j \cap \Gamma^{\nu}_{j}$, see \eqref{eq:Aj} and \eqref{eq:gammajnu}, one has that 
$$
    |\partial_{\xi'}^\alpha h_j^\nu(x,\xi)| \lesssim \begin{cases}
    |\xi'| |\xi|^{-|\alpha|} & \text{ if } |\alpha|= 1  \\ 
    |\xi| ^{1-|\alpha|} & \text{ if } |\alpha| \geq 2
    \end{cases} \,\lesssim 2^{-j\rho/2}.
    $$
Moreover, for $\beta\in \mathbb Z_{> 0}^{n-\kappa}$ one has that
$$
    |\partial_{\xi''}^\beta h_j^\nu(x,\xi)|\lesssim |\xi'|^2 |\xi| ^{-|\beta|-1}\lesssim 2^{-j\rho}.
$$
\end{Lem}

\begin{proof}
The proof is based on \eqref{eq:xiprimestimate} simple Taylor expansions and homogeneity considerations, see \cite[p. 407]{Stein}. 
\end{proof}

\begin{Rem}\label{rem:alpha=1}
    Note that by the previous lemma, one has the estimate
    \begin{equation}\label{eq:trunkerad0}
        |\nabla_{\xi} h^{\nu}_j(x,\xi)|\lesssim 2^{-j\rho /2},
    \end{equation}
    for all $\xi\in  \Gamma^{\nu}_j,$ due to the fact that $\nabla_{\xi} h^{\nu}_j (x, \xi)$ is homogeneous of degree zero in $\xi.$ 
\end{Rem}

In \cite{SSS} the authors define an \lq\lq influence set\rq\rq\ associated to the SND phase function $\varphi$. We have to make a similar definition but it has to be fitted to the more general classes of amplitudes that we are considering here. To this end we have

\begin{Def}\label{def:influenceset}
Assume that $\varphi$ is an \emph{SND} phase function of rank $\kappa$ in the class $\Phi^2$ and let $Q$ be a cube. Define, for $j\in \Z_{\geq0}$, $1\leq\nu \leq \mathscr N_j,$ the set $R_j^\nu$ as
\eq{
    R_j^\nu := \set{x\in \Rn : |(\nabla_\xi \varphi(x,\xi_j^\nu)-\bar y)''|\leq c 2^{-\rho j}, \quad | (\nabla_\xi \varphi(x,\xi_j^\nu)-\bar y)'| \leq c 2^{-\rho j/2}}.
}

Here, $c$ is a large constant depending on the size of the Hessian matrix of $\varphi$ but independent of $j$ and $\nu$, to be specified later. Now set
\nm{eq:Bstar}{
    Q^*_\rho := 
        \diffcases{
            \displaystyle\bigcup_{j= k_Q }^\infty\bigcup_{\nu=1}^{\mathscr N_j} R_j^\nu,& \rho\in(0,1]\text{ and } k_Q\geq 1
        \\
            2\sqrt n Q,&  \rho= 0 \text{ or } k_Q\leq 0.
        }
}

\end{Def}
In this connection, we have the following estimates.

\begin{Lem}\label{lem:sizeofQstar}
Assume that $\varphi$ is an \emph{SND} phase function of rank $\kappa$ in the class $\Phi^2$ and let $Q$ be a cube with $k_Q\geq 1$. Moreover, let $Q_\rho^*$ be defined as in \eqref{eq:Bstar}. Then for $\rho\in (0,1]$ we have the following properties.
\begin{enumerate}
    \item[$(i)$]  The measure of $Q_\rho^*$ satisfies
        \eq{|Q^*_\rho| \lesssim 2^{- k_Q \rho (n-\kappa)}.}
    \item[$(ii)$] If $x\in \Rl^n \setminus Q_\rho^*$ and $y\in Q$, 
    then for $c$ in \emph{Definition \ref{def:influenceset}} large enough, we have
    \begin{equation}\label{eq:keypoint}
        2^{2j\rho}|(\nabla_\xi \varphi(x,\xi_j^\nu)-y)'|^2+ 2^{j\rho}|(\nabla_\xi \varphi(x,\xi_j^\nu)-y)''|^2  \gtrsim 2^{(j-k_Q)\rho }, \quad j \geq k_Q
    \end{equation}
\end{enumerate}

\end{Lem}

\begin{Rem}
By a similar argument, one can prove the following version of \eqref{eq:keypoint}.
\begin{equation}\label{eq:keypoint2}
    2^{2j\rho}|(x-\nabla_\xi \varphi(y,\xi_j^\nu))'|+ 2^{j\rho}|(x-\nabla_\xi \varphi(y,\xi_j^\nu))''|  \gtrsim 2^{(j-k_Q)\rho }, \quad j \geq k_Q.
\end{equation}
This will be used for proving $h^p\to L^p$ boundedness for the adjoint.
\end{Rem}
\begin{proof} Recall that $\mathscr N_j = O\big(2^{j\rho\kappa/2}\big)$, $k_Q :=[1-\log_2(l_Q)]$ and that $c_Q$ is the centre of $Q.$\\

\textbf{Step 1 -- Proof of (\textit{i})}\\
Since $R_j^\nu$ is of size $O(2^{-j\rho/2})$ in $\kappa$ directions and $O(2^{- j\rho})$ in the other $n-\kappa$ directions, we have for $0<\rho\leq 1$ that
\eq{|Q^*_\rho| \leq \sum_{j=k_Q}^\infty \sum_{\nu=1}^{\mathscr N_j} |R_j^\nu| \lesssim \sum_{j=k_Q}^\infty 2^{j\rho \kappa/2}\,2^{-j(\rho (n-\kappa) + \rho \kappa/2)}\lesssim 2^{-k_Q\rho (n-\kappa)}.}

\textbf{Step 2 -- Proof of (\textit{ii})}\\
Observe that \eqref{eq:keypoint} is equivalent to
\begin{equation*}
2^{(j+ k_Q )\rho/2} |(\nabla_\xi \varphi(x,\xi_j^\nu)-  y)''|
+2^{ k_Q \rho/2}|(\nabla_\xi \varphi(x,\xi_j^\nu)-  y)'|  
\gtrsim 1, \quad j \geq k_Q.
\end{equation*}
Moreover, it is enough to show that
\begin{equation}\label{lower bound estimate3} 
2^{ k_Q \rho} |(\nabla_\xi \varphi(x,\xi_j^\nu)-  y)''|
+2^{{k_Q} \rho/2}|(\nabla_\xi \varphi(x,\xi_j^\nu)-  y)'|  
\gtrsim 1.
\end{equation}
Since  $\xi^{\nu}_j \in \mathbb{S}^{n-1},$  by Definition \ref{def:LP2} part ($ii$), there exists a unit vector $\xi^{\nu'}_{k_Q}$  such that
\begin{equation}\label{close points on sphere}
    |\xi^{\nu}_j-\xi^{\nu'}_{k_Q}|<  2^{-{k_Q}\rho /2}.  
\end{equation}

We now state and prove the following three assertions in order to conclude the proof. Observe that $c$ in \eqref{important_estimate} comes from Definition \ref{def:influenceset}.
\begin{align}
&
    2^{{k_Q}\rho}|(\nabla_\xi \varphi(x,\xi_{k_Q}^{\nu'})-c_Q)''| + 2^{{k_Q}\rho/2}| (\nabla_\xi \varphi(x,\xi_{k_Q}^{\nu'})-c_Q)'|\geq c,\label{important_estimate}
\\&
    |(\nabla_\xi \varphi (x,\xi^{\nu}_j) - \nabla_\xi \varphi (x,\xi^{\nu'}_{k_Q}))'|\lesssim {2^{-{k_Q}\rho /2}},\label{pain3}
\\&
    |(\nabla_\xi \varphi (x,\xi^{\nu}_j) - \nabla_\xi \varphi (x,\xi^{\nu'}_{k_Q}))''|\lesssim {2^{-{k_Q}\rho }}\label{pain2}
\end{align}

Before we prove the above assertions we start with proving that
\begin{equation}\label{pain1}
    |\nabla_\xi \varphi (x,\xi^{\nu}_j) - \nabla_\xi \varphi (x,\xi^{\nu'}_{k_Q})|\lesssim 2^{-{k_Q}\rho/2}.
\end{equation}
Indeed, we observe that since $\varphi\in \Phi^2$ one has that $|\partial^2_{\xi\xi}\varphi(x,\xi)|\lesssim |\xi|^{-1},$ for all $x$ and all $\xi\neq 0.$ Therefore the mean-value theorem and \eqref{close points on sphere} yield that for some $t\in [0,1]$ and all ${k_Q}\geq 1$ one has
\begin{align}
&|\nabla_\xi \varphi (x,\xi^{\nu}_j) - \nabla_\xi \varphi (x,\xi^{\nu'}_{k_Q})|
 \leq | \xi^{\nu}_j- \xi^{\nu'}_{k_Q}||\xi^{\nu}_j+t(\xi^{\nu'}_{k_Q}-\xi^{\nu}_j)|^{-1} \nonumber \\
& \leq {2^{-{k_Q}\rho/2}} (1-t2^{-{k_Q}\rho/2})^{-1}
\leq {2^{-{k_Q}\rho /2}} (1-2^{-\rho/2})^{-1} \lesssim 2^{-{k_Q}\rho /2}.\nonumber
\end{align}
which proves \eqref{pain1}.\\

\makeatletter 
\renewcommand{\eqref}[1]{\tagform@{\ref{#1}}}
\def\maketag@@@#1{\hbox{#1}}
\textbf{Step 2.1 -- Proof of \eqref{important_estimate} and \eqref{pain3}}\\
\makeatother
Firstly, \eqref{important_estimate} follows directly from the fact that $x\in \Rn \setminus Q^*_\rho\subset \Rn \setminus R_{k_Q}^{\nu'}$. And \eqref{pain3} trivially follows from \eqref{pain1}.\\

\makeatletter 
\renewcommand{\eqref}[1]{\tagform@{\ref{#1}}}
\def\maketag@@@#1{\hbox{#1}}
\textbf{Step 2.3 -- Proof of \eqref{pain2}}\\
\makeatother
We first claim that
\begin{equation}\label{eq:SSStrick}
|(\nabla_\xi \varphi (x,\xi^{\nu}_j) - \nabla_\xi \varphi (x,\xi^{\nu'}_{k_Q}))\cdot \xi^{\nu'}_{k_Q}|\lesssim 2^{-{k_Q}\rho }.
\end{equation}
To see this, using the homogeneity of $\varphi$ (dictated by the $\Phi^2$-condition) we have 
\begin{equation*}
\begin{split}
&(\nabla_\xi \varphi (x,\xi^{\nu}_j) - \nabla_\xi \varphi (x,\xi^{\nu'}_{k_Q}))\cdot \xi^{\nu'}_{k_Q}\\& = (\nabla_\xi \varphi (x,\xi^{\nu}_j) - \nabla_\xi \varphi (x,\xi^{\nu'}_{k_Q}))\cdot (\xi^{\nu'}_{k_Q}- \xi^{\nu}_j)+ \varphi (x,\xi^{\nu}_j) - \nabla_\xi \varphi (x,\xi^{\nu'}_{k_Q})\cdot \xi^{\nu}_j\\& = (\nabla_\xi \varphi (x,\xi^{\nu}_j) - \nabla_\xi \varphi (x,\xi^{\nu'}_{k_Q}))\cdot (\xi^{\nu'}_{k_Q}- \xi^{\nu}_j)+ h^{\nu'}_{k_Q}(x,\xi^{\nu}_{j}),
\end{split}
\end{equation*}
where 
$$h^{\nu'}_{k_Q}(x,\xi) :=\varphi (x,\xi) - \nabla_\xi \varphi (x,\xi^{\nu'}_{k_Q})\cdot \xi.$$
 Moreover, by \eqref{eq:trunkerad0} one has the estimate
\begin{equation}\label{eq:trunkerad}
    |(\nabla_{\xi} h^{\nu'}_{k_Q})(x,\xi)|\leq 2^{-{k_Q}\rho /2},
\end{equation} for all $\xi \in \Gamma^{\nu'}_{k_Q}$
. 
Recalling that $\xi^{\nu}_j$ is in $\mathbb{S}^{n-1}$, \eqref{close points on sphere} shows that $\xi^{\nu}_{j}$ also belongs to the cone $\Gamma^{\nu'}_{{k_Q}}$, and so for all $t\in [0,1]$, the expression $\xi^{\nu}_j+t(\xi^{\nu'}_{k_Q}-\xi^{\nu}_j)$ which represents the line segment joining  $\xi^{\nu}_j$ and $\xi^{\nu'}_{k_Q}$  belongs to $\Gamma^{\nu'}_{{k_Q}},$ due to the convexity of the cone.
Therefore, \eqref{eq:trunkerad} yields that for all $t\in [0,1]$
$$|(\nabla_{\xi} h^{\nu'}_{k_Q})(x,\xi^{\nu}_j+t(\xi^{\nu'}_{k_Q}-\xi^{\nu}_j))|\lesssim 2^{-{k_Q}\rho /2}.$$ Now we also observe that $h^{\nu'}_{k_Q}(x,\xi^{\nu'}_{k_Q})=0$. Hence, using the mean-value theorem and \eqref{close points on sphere}, one readily sees that $h^{\nu'}_{k_Q}(x,\xi^\nu_j)= O(2^{-{k_Q}\rho }).$\\

For the term $ (\nabla_\xi \varphi (x,\xi^{\nu}_j) - \nabla_\xi \varphi (x,\xi^{\nu'}_{k_Q}))\cdot (\xi^{\nu'}_{k_Q}- \xi^{\nu}_j)$ we just use the Cauchy-Schwarz inequality, \eqref{close points on sphere} and \eqref{pain1}, which concludes the proof of \eqref{eq:SSStrick}.\\

Now we turn to the proof of  \eqref{pain2}. Our initial convention that $\xi^{\nu}_j$ lies along the $\xi_1$-axis, the triangle inequality, the Cauchy-Schwarz inequality, \eqref{close points on sphere}, \eqref{pain1} and \eqref{eq:SSStrick}, yield that
\begin{align*}
& |(\nabla_\xi \varphi (x,\xi^{\nu}_j) - \nabla_\xi \varphi (x,\xi^{\nu'}_{k_Q}))''|
= 
|(\nabla_\xi \varphi (x,\xi^{\nu}_j) - \nabla_\xi \varphi (x,\xi^{\nu'}_{k_Q}))\cdot\xi^{\nu}_j| \nonumber \\
&  
= |(\nabla_\xi \varphi (x,\xi^{\nu}_j) - \nabla_\xi \varphi (x,\xi^{\nu'}_{k_Q}))(\xi^{\nu}_j-\xi^{\nu'}_{k_Q})+ (\nabla_\xi \varphi (x,\xi^{\nu}_j) - \nabla_\xi \varphi (x,\xi^{\nu'}_{k_Q}))\cdot\xi^{\nu'}_{k_Q}| \nonumber \\
&  \leq
|\nabla_\xi \varphi (x,\xi^{\nu}_j) - \nabla_\xi \varphi (x,\xi^{\nu'}_{k_Q})| \, 
|\xi^{\nu}_j-\xi^{\nu'}_{k_Q}|
+
|(\nabla_\xi \varphi (x,\xi^{\nu}_j) - \nabla_\xi \varphi (x,\xi^{\nu'}_{k_Q}))\cdot\xi^{\nu'}_{k_Q}|
\nonumber \\
&  \lesssim 
2^{-{k_Q}\rho }.
\end{align*}

\makeatletter 
\renewcommand{\eqref}[1]{\tagform@{\ref{#1}}}
\def\maketag@@@#1{\hbox{#1}}
\textbf{Step 3 -- Proof of \eqref{lower bound estimate3}}\\
\makeatother
Finally, to show \eqref{lower bound estimate3} (which as we mentioned above implies the desired estimate \eqref{eq:keypoint}, we use the triangle inequality, \eqref{important_estimate},
\eqref{pain3} and \eqref{pain2} to obtain 
\begin{align*}
    & 2^{{k_Q}\rho}|(\nabla_\xi \varphi(x,\xi_j^\nu)-  y)''|
    + 2^{{k_Q} \rho/2 }|(\nabla_\xi \varphi(x,\xi_j^\nu)-  y)'| 
\\& =  
    2^{{k_Q}\rho}|\big(\nabla_\xi \varphi(x,\xi_j^\nu)-c_Q - (y-c_Q)\big)''|
    + 2^{{k_Q}\rho/2 }\big|\big(\nabla_\xi \varphi(x,\xi_j^\nu)-c_Q - (y-c_Q)\big)'\big| 
\\& \geq 
    2^{{k_Q}\rho}|\big(\nabla_\xi \varphi(x,\xi_j^\nu)-c_Q \big)''| - 2^{{k_Q}\rho}|y''-c_Q''|
    + 2^{{k_Q}\rho/2 }\big|\big(\nabla_\xi \varphi(x,\xi_j^\nu)-c_Q \big)'\big|
\\ &\qquad -
         2^{{k_Q}\rho/2 }|y'-c_Q'| 
\\& \geq 
    2^{{k_Q}\rho}\big|\big(\nabla_\xi \varphi(x,\xi_j^\nu)-c_Q \big)''\big| 
    + 2^{{k_Q}\rho/2 }\big|\big(\nabla_\xi \varphi(x,\xi_j^\nu)-c_Q \big)'\big|
\\ &\qquad -
        A2^{{k_Q}\rho} 2^{-k_Q} -A2^{{k_Q}\rho/2 } 2^{-k_Q} 
\\& = 
    2^{{k_Q}\rho}\big|\big(\nabla_\xi \varphi(x,\xi_j^\nu) - \nabla_\xi \varphi(x,\xi_{k_Q}^{\nu'}) + \nabla_\xi \varphi(x,\xi_{k_Q}^{\nu'})-c_Q \big)''\big| 
\\& \qquad +
         2^{{k_Q}\rho/2 }\big|\big(\nabla_\xi \varphi(x,\xi_j^\nu) - \nabla_\xi \varphi(x,\xi_{k_Q}^{\nu'}) + \nabla_\xi \varphi(x,\xi_{k_Q}^{\nu'})-c_Q \big)'\big|
\\ &\qquad -
        A2^{{k_Q}\rho} 2^{-k_Q}-A 2^{{k_Q}\rho/2 } 2^{-k_Q}
\\ & \geq 
    2^{{k_Q}\rho}\big|\big( \nabla_\xi \varphi(x,\xi_{k_Q}^{\nu'})-c_Q \big)''\big| - 2^{{k_Q}\rho} |\big(\nabla_\xi \varphi(x,\xi_j^\nu) - \nabla_\xi \varphi(x,\xi_{k_Q}^{\nu'}) \big)''|
\\& \qquad + 
        2^{{k_Q}\rho/2 }\big|\big( \nabla_\xi \varphi(x,\xi_{k_Q}^{\nu'})-c_Q \big)'\big| - 2^{{k_Q}\rho/2 }\big|\big(\nabla_\xi \varphi(x,\xi_j^\nu) - \nabla_\xi \varphi(x,\xi_{k_Q}^{\nu'}) \big)'\big| 
\\& \qquad  -
        A 2^{{k_Q}\rho} 2^{-k_Q}-A2^{{k_Q}\rho/2 } 2^{-k_Q}
\\& \geq 
     2^{{k_Q}\rho}\big|\big( \nabla_\xi \varphi(x,\xi_{k_Q}^{\nu'})-c_Q \big)''\big| +  2^{{k_Q}\rho/2 }\big|\big( \nabla_\xi \varphi(x,\xi_{k_Q}^{\nu'})-c_Q \big)'\big| 
\\& \qquad -
        2^{{k_Q}\rho}B 2^{-{k_Q}} - 2^{{k_Q}\rho/2 }C2^{-{k_Q}/2} -A2^{{k_Q}\rho} 2^{-k_Q}-A2^{{k_Q}\rho/2 } 2^{-k_Q}
\\& \geq 
     c - 2^{{k_Q}\rho}B  2^{-{k_Q}} - C 2^{{k_Q}\rho/2 }2^{-{k_Q}/2} 
\\& = 
    c - B 2^{-{k_Q}(1-\rho)} - C 2^{-{k_Q}(1-\rho )/2} -A2^{-{k_Q}(1-\rho)} -A2^{-{k_Q}(1-\rho/2 )}
\\& \geq 
    c - A - B - C,
\end{align*}
where the constants $A=2\sqrt n$, $B$ stems from estimate \eqref{pain2} and $C$ from \eqref{pain3}. Therefore picking $c$ large enough we obtain \eqref{eq:keypoint}. \\
\begin{Rem}
    Note that for $\rho>\frac{2}{3}$ we do not use the scaling in \emph{\eqref{eq:LPdecompkernel}} and follow the same lines of reasoning as in \cite{SSS}. This also means that the discussions starting from \emph{Definition \ref{def:LP2}} up-to and including \emph{Lemma \ref{lem:sizeofQstar}} have to be modified according to the decomposition given on pages \emph{247-248} in \cite{SSS}, which indeed works for $\rho\geq\frac{1}{2}$.
\end{Rem}

\end{proof}
\begin{Lem}\label{lemma:HLS-lemma}
Let $\tilde m\leq 0$, $n\geq 1,$ $\rho\in [0,1]$, $\delta\in [0,1)$ and 
\begin{equation*}
m:=\tilde m-n \max\set{0,\frac{\delta-\rho}2}.
\end{equation*} 
Suppose that $a\in S^m_{\rho,\delta}(\Rn)$ and that $a(x,\xi)$ vanishes in a neighborhood of $\xi=0$. Also, let $\varphi$ be an \emph{SND} phase function in the class $\Phi^2$ with rank $\kappa$.
Then $T_a^\varphi$, see \emph{Definition \ref{def:FIO},} satisfies
\nm{eq:L2toLq}{
\|T_a^\varphi f\|_{L^2(\Rn)} 
\lesssim 
\|f\|_{h^{2n/(n-2 \tilde m)}(\Rn)}.}
Also for the adjoint operator one has
\nm{eq:L2toLqadj}{
\|(T_a^{\varphi})^* f\|_{L^2(\Rn)} \lesssim 
\|f\|_{h^{2n/(n-2 \tilde m)}(\Rn)}.}
\end{Lem}

\begin{proof}
Since the operator $T_a^\varphi (1-\Delta)^{-\tilde m/2}$ is an FIO with the phase $\varphi$ and an amplitude in $S^{-n \max\{(\delta-\rho)/2,0\}}_{\rho,\delta}(\Rn)$ it is $L^2$-bounded by Theorem  \ref{basicL2}. This $L^2$ boundedness together with the estimates for the Bessel potential operators reformulated in terms of embedding of Triebel-Lizorkin spaces (see \cite[Corollary 2.7]{Triebel4}) yield
\eq{
&\|T_a^\varphi f\|_{L^2(\Rn)}  = \|T_a^\varphi (1-\Delta)^{-\tilde m /2}(1-\Delta)^{\tilde m /2}f\|_{L^2(\Rn)} \\ 
&\lesssim \|(1-\Delta)^{\tilde m /2 }f\|_{L^2(\Rn)} 
\lesssim \|f\|_{h^{2n/(n-2 \tilde m)}(\Rn)},}
which proves \eqref{eq:L2toLq}. Observe that the choice of the range of $\tilde m$ implies that $0<  2n/(n-2 \tilde m)\leq 2$.\\
Next we prove \eqref{eq:L2toLqadj}. By Theorem \ref{thm:left composition with pseudo} the composition  $ (1-\Delta)^{-m' /2} T_a^\varphi$ is an FIO with the phase $\varphi$ and an amplitude in $S^{-n \max\{(\delta-\rho)/2,0\}}_{\rho,\delta}(\Rn)$, and therefore $L^2$-bounded. Finally, observing that $$\Vert (T_a^\varphi)^*(1-\Delta)^{-m' /2}\Vert_{L^2 \to L^2} = \Vert  (1-\Delta)^{-m' /2}T_a^\varphi\Vert_{L^2 \to L^2} ,$$ one can proceed as above. 
\end{proof}

\begin{figure}
\begin{tikzpicture}

    \draw[thick,->] (0,0) -- (6.5,0) node[anchor=north west] {$p$};
    \draw[thick,->] (0,-1) -- (0,1) node[anchor=south east] {$m_c$};

    \draw[thick] (6 cm,-2pt) -- (6 cm,2pt) node[anchor=south] {$\infty$};
    \draw[thick] (0.5,-2pt) -- (0.5,2pt) node[anchor=south] {$1$};
    \draw[thick] (3.25,-2pt) -- (3.25,2pt) node[anchor=south] {$2$};
    \draw[thick] (2pt,-0.7) -- (-2pt,-0.7) node[anchor=east] {$-\frac{\kappa+(n-\kappa)(1-\rho)}2$};

    \draw (0.5,-0.7) -- (3.25,0);
    \draw (6,-0.7) -- (3.25,0);
    
    \filldraw[gray, opacity = 0.6, path fading = south] (0.5,-0.7) -- (3.25,0) -- (6,-0.7) -- (6,-1) -- (0.5,-1) -- cycle;
\end{tikzpicture}
\caption{Critical order when $\rho\geq\delta$. }\label{pic:criticaldecay}
\end{figure}
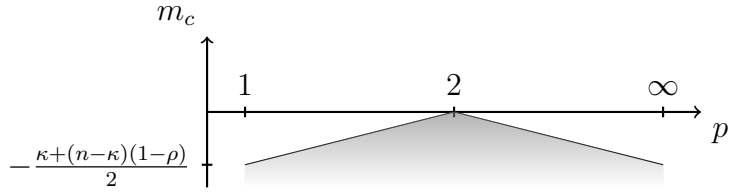

{\bf{Notation.}} {We shall henceforth set
\begin{equation}\label{eq:criticaldecay}
    m_c(p) := -\big(\kappa+(n-\kappa)(1-\rho)\big)\Big|\frac{1}{p}-\frac{1}{2}\Big| - n\max\set{0,\frac{\delta-\rho}2},
\end{equation}
which as will turn out, be the critical order of the amplitudes for bounded FIOs (see \ref{pic:criticaldecay}).}

\begin{Lem}\label{lem:bstarcestimate}
Let $T_a^\varphi$ be an $\mathrm{FIO}$ with an amplitude $a\in S^{m_c(p)}_{\rho,\delta}(\Rn)$ for $\rho\in\big[0,\frac23\big],\,\delta\in[0,1)$ and $m_c(p)$ in \eqref{eq:criticaldecay}. Let $\varphi \in \Phi^2$ be an $\mathrm{SND}$ phase function of rank $\kappa$ satisfying the conditions in  \emph{Definition \ref{def:FIO}}. Moreover $T_j:= T_a^\varphi\psi_j(D)$, where $\psi_j(D)$ is a Littlewood-Paley piece as in \emph{Definition \ref{def:LP}}.  Furthermore, suppose that $f$ is supported in a cube $Q$ with
\eq{
\|f\|_{L^{ 1}(\Rn)}\lesssim 2^{k_Qn(1/p-1)}.
}
Also let $Q^*_\rho$ be as in \eqref{eq:Bstar}.
\begin{enumerate}
    \item[$(i)$] If $0<p<\frac{2n-2n\rho+2\kappa\rho}{2n-n\rho+\kappa\rho}$, then
        \nm{eq:LPoperator2}{
            \sum_{j=\max\{k_Q+1,1\}}^\infty\| T_j  f\|_{L^p(\Rn \setminus Q^*_\rho)}^p\lesssim 1.
        }
    \item[$(ii)$] If $0<p\leq 1$ and $f$ is a $h^p$-atom, see \emph{Definition \ref{def:hpatom}}, then
        \nm{eq:LPoperator1}{
            \sum_{j=1}^{\max\{k_Q,0\}}\| T_j  f\|_{L^p(\Rn \setminus Q^*_\rho)}^p\lesssim 1.
        }
    
\end{enumerate}
Moreover, the same estimates hold true for the adjoint operator $(T_j)^*.$

\end{Lem}
\begin{Rem}\label{rem:pseudoTLremark}
Observe that the condition on $p$ in \emph{Lemma \ref{lem:bstarcestimate}} $(i)$ excludes the case when $\kappa=0$ and $\rho=1$. This case is treated in e.g. \cite{Paivarinta:pdo-TL}, \cite{Triebel:pdo-TL}.
\end{Rem}
\begin{proof}[Proof of \emph{Lemma \ref{lem:bstarcestimate}}]

Observe that it is enough to prove \eqref{eq:LPoperator1} when $k_Q\geq 1$, otherwise the sum is trivially $0$. Moreover, it is enough to show that for some $\epsilon>0$,
\begin{align*}
    &\| T_j  f\|_{L^p(\Rn \setminus Q^*_\rho)} \lesssim 2^{(\max\{k_Q,0\}-j)\varepsilon},\quad j\geq \max\{k_Q+1,1\},\\
    &\| T_j  f\|_{L^p(\Rn \setminus Q^*_\rho)} \lesssim 2^{(j-k_Q)\varepsilon},\quad \quad \quad \;\;1\leq j\leq \max\{k_Q,0\}.
\end{align*}

We start the proof of this lemma with some useful prerequisites. First of all, we will assume that $x\in \Rn\setminus Q^*_\rho$ and $y\in Q$ throughout the proof. Next, we define two types of metrics that will simplify the exposition of the proof.
\begin{align*}
    g(x,y)&:=1+ |x-y|^2,\\
    g_j^\nu(x,y)&=\big(1+2^{2j\rho}|( \nabla_\xi \varphi(x,\xi_j^\nu)-y)''|^2\big)\big(1+2^{j\rho}|( \nabla_\xi \varphi(x,\xi_j^\nu)-y)'|^2\big).
\end{align*}

Also, we define some variables that will be used throughout the proof. Let $p':=\frac{2p}{2-p}$ or equivalently $\frac 1{p'}=\frac1p-\frac12$, $m':=-\big(\kappa+(n-\kappa)(1-\rho)\big)\big|\frac{1}{p}-\frac{1}{2}\big|$ and $\mathfrak M_{\at}:=\big[n\big(\frac 1p-1\big)\big]$ as in Definition \ref{def:hpatom}. We also need $N_1>\frac{n}{2p'}$,  $N_2>8N_1$ and $N_3 > 2N_1+\frac{N_2}\rho+n$. \\

Next, we prove the following two assertions. For $t\in [0,1],$ $x\in  \Rn\setminus 2\sqrt n Q$ and $y\in Q$, one has 
\begin{equation}\label{eq:estimateybar}
    \vert x-c_Q\vert \lesssim \vert x-c_Q-t(y-c_Q)\vert.
\end{equation}
For all integers $N \geq 0$ and the kernel $K_j(x,y)$, defined in \eqref{eq:LPdecompkernel},
\nm{eq:newone}{
    \big\| g(x,y)^N\, K_{j}(x,y)\big\|_{L^2_x( \Rn)}&\lesssim 2^{j(n/2+ m')}.
}

We divide this proof into different steps, starting with proving the assertions above.\\

\makeatletter 
\renewcommand{\eqref}[1]{\tagform@{\ref{#1}}}
\def\maketag@@@#1{\hbox{#1}}
\textbf{Step 1 -- Proof of \eqref{eq:estimateybar}}\\
\makeatother
We have
\eq{
    &|x-c_Q| \leq |x-c_Q-t(y-c_Q)|+t|y-c_Q| \leq |x-c_Q-t(y-c_Q)|+tl_Q
\\&\leq 
    |x-c_Q-t(y-c_Q)|+\frac{t}{2\sqrt{n}}|x-c_Q|.
}
The assertion follows by rearrangements in the inequality above.\\

\makeatletter 
\renewcommand{\eqref}[1]{\tagform@{\ref{#1}}}
\def\maketag@@@#1{\hbox{#1}}
\textbf{Step 2 -- Proof of \eqref{eq:newone}}\\
\makeatother
In this step we will use the arguments in Section \ref{sec:phase_reduction}, which show that it is enough to consider phase functions of the form $\varphi(x,\xi)= x\cdot\xi+\theta(x,\xi)$, where $\theta(x,\xi)\in\Phi^1$. We claim that  
\nm{eq:kernelestimate1}{
\Vert (x-y)^\alpha \,\partial^{\beta} _{y}K_{j}(x,y)\Vert_{L^2_x(\Rn)}\lesssim 2^{j(  \vert \beta\vert+ m'+n/2)},
}
for all multi-indices $\alpha$ and $\beta$. Since differentiating $\beta$ times in $y$ will only introduce factors of the size $2^{j|\beta|},$ it is enough to establish \eqref{eq:kernelestimate1} for $\beta=0$. Now the global $L^2$ boundedness \eqref{eq:kernelestimate1} of the kernel can be formulated as the $L^2$ boundedness of a kernel of the form
\begin{equation*}
\tilde K^\alpha_j(x,x-y):= \int_{\Rn} a_j(x,\xi)\,(x-y)^\alpha\, e^{i\theta(x,\xi)+i(x-y)\cdot\xi}  \ddd\xi.
\end{equation*}
To this end, take ${\Psi}_j$ as in Definition \ref{def:LP}, integrate by parts, and rewrite 
\eq{
    &\tilde K^\alpha_j(x,x-y)
    =\int_{\Rn} a_j(x,\xi)\,e^{i\theta(x,\xi)}\,
    (-i)^{|\alpha|}\,\partial_\xi^\alpha e^{i(x-y)\cdot\xi} \ddd\xi\\
& 
    = i^{|\alpha|} \int_{\Rn} \partial_\xi^\alpha \Big[a_j(x,\xi)e^{i\theta(x,\xi)}\Big]\, e^{i(x-y)\cdot\xi} \,{\Psi}_j (\xi)\ddd\xi \\
& 
    = \sum_{\alpha_1 + \alpha_2=\alpha} \!\!C_{\alpha_1, \alpha_2}\int_{\Rn}
    \partial_\xi^{\alpha_1} a_j(x,\xi) \,
    \partial_\xi^{\alpha_2} e^{i\theta(x,\xi)} \,e^{i(x-y)\cdot\xi}\, {\Psi}_j (\xi)\ddd\xi \\
&  =
    \!\!\sum_{\substack{\alpha_1 + \alpha_2=\alpha \\ \lambda_1 + \dots + \lambda_r = \alpha_2} }\!\!\!\!
    C_{\alpha_1, \alpha_2, \lambda_1, \dots \lambda_r} \int_{\Rn} 
    \partial_\xi^{\alpha_1} a_j(x,\xi) \\
&  \qquad \qquad   \times
    \partial_\xi^{\lambda_1}\theta(x,\xi)
    \cdots
    \partial_\xi^{\lambda_r}\theta(x,\xi)
    \,e^{i\theta(x,\xi)} \,e^{i(x-y)\cdot\xi}\,{\Psi}_j (\xi) \ddd\xi \\
&  =\!\!
    \sum_{\substack{\alpha_1 + \alpha_2=\alpha \\ \lambda_1 + \dots + \lambda_r = \alpha_2} }\!\!\!\!
    C_{\alpha_1, \alpha_2, \lambda_1, \dots \lambda_r} 
    2^{jm'}  \!\!
    \int_{\Rn} b_j^{\alpha_1,\alpha_2, \lambda_1, \dots, \lambda_r}(x,\xi)\,
     e^{i\theta(x,\xi)+ix\cdot\xi}\, e^{-iy\cdot\xi}\,{\Psi}_j (\xi) \ddd\xi \\
&  =:\!\!
    \sum_{\substack{\alpha_1 + \alpha_2=\alpha \\ \lambda_1 + \dots + \lambda_r = \alpha_2} }
    \!\!\!\!C_{\alpha_1, \alpha_2, \lambda_1, \dots, \lambda_r} 
    2^{jm'} \,  
    S_{j}^{\alpha_1, \alpha_2, \lambda_1, \dots, \lambda_r}(\tau_{-y}\Psi_j^\vee)(x),
}
where $S_{j}^{\alpha_1, \alpha_2, \lambda_1, \dots, \lambda_r}$ is an FIO with  the phase function $\theta(x,\xi)+x\cdot\xi$ and amplitude $b_j^{\alpha_1,\alpha_2, \lambda_1, \dots, \lambda_r}(x,\xi)$ given by
\begin{equation*}
b_j^{\alpha_1,\alpha_2, \lambda_1, \dots, \lambda_r}(x,\xi)
 := 2^{-jm'}\, \partial_\xi^{\alpha_1} a_j(x,\xi) \,
\partial_\xi^{\lambda_1}\theta(x,\xi)
\dots
\partial_\xi^{\lambda_r}\theta(x,\xi).   
\end{equation*}
Moreover  $|\lambda_j| \geq 1$ and $\tau_{-y}$ is a translation by $-y$.\\

We observe that 
$b_j^{\alpha_1,\alpha_2, \lambda_1, \dots, \lambda_r}(x,\xi) \in S^{-n \max\{(\delta-\rho)/2,0\}}_{\rho,\delta}(\Rn)$ uniformly in $j$, since $a \in S^{m}_{\rho,\delta}(\Rn)$ and $\theta\in \Phi^1$.\\

Therefore by Theorem \ref{basicL2}, $S_{j}^{\alpha_1, \alpha_2, \lambda_1, \dots \lambda_r}$ is an $L^2$-bounded FIO, so
\begin{align*}
    &\Vert \tilde K^\alpha_j(x,x-y)\Vert_{L^2_x(\Rn)}
    \lesssim\sum_{\substack{\alpha_1 + \alpha_2=\alpha \\ \lambda_1 + \dots + \lambda_r = \alpha_2} }\!\! 2^{jm'} 
    \|S_{j}^{\alpha_1, \alpha_2, \lambda_1, \dots \lambda_r}(\tau_{y}\Psi_j^\vee)\|_{L^2(\Rn)} \\
& \lesssim 
    2^{jm'} \| \Psi_j\|_{L^2(\Rn)} 
    \lesssim 2^{j(   m'+n/2)},
\end{align*}
which proves \eqref{eq:kernelestimate1}.\\

Now, \eqref{eq:kernelestimate1} yields \eqref{eq:newone} for any integer $N\geq 0$, if one sums over $|\alpha|\leq 2N$.\\

\textbf{Step 3 -- The case when  $\boldsymbol{\rho}\,\mathbf{=0}$}\\
The H\"older and the Minkowski inequalities together with \eqref{eq:newone} and \eqref{eq:estimateybar} (with $t=1$) yield
\begin{align*}
    &\|T_jf \|_{L^p( \Rn\setminus Q_0^*)}= \Big\| \int_{Q} K_{j}(x,y)\,f(y)\dd y \Big\| _{L^p_x(\Rn\setminus Q_0^*)} 
\\ 
    & \leq \Big\|   g(x,c_Q)^{N_1} \int_{Q}K_{j}(x,y)\,f(y)\dd y \Big\| _{L^2_x( \Rn\setminus Q_0^*)} \,\norm{ \frac{1}{g(x,c_Q)^{N_1}} } _{L^{p'}( \Rn)}  \nonumber 
\\ 
    &\lesssim  \int_{Q} \big\|g(x,c_Q)^{N_1}\,K_{j}(x,y)\,f(y)\big\|_{L^2_x( \Rn\setminus Q_0^*)}  \dd y \nonumber 
\\
    &\lesssim   \int_{Q} \vert f(y) \vert \, \big\| g(x,y)^{N_1} \,   K_{j}(x,y)\big\| _{L^2_x( \Rn\setminus Q_0^*)} \!\dd y \nonumber 
\\
    &\lesssim 2^{-k_Q(n-n/p)}\,   2^{j(n/2+ m')}= 2^{-k_Q(n-n/p)}\,2^{j(n-n/p)}, \nonumber 
\end{align*}
since $ m'=-n\big(\frac 1p-\frac 12\big)$. Hence \eqref{eq:LPoperator2} is proved when $\rho=0$.\\

To prove \eqref{eq:LPoperator1} we expand $K_{j}(x,y)$ in a Taylor polynomial of order $\mathfrak M_{\at}$ at the point $y=c_Q$, which yields that
 \eq{
& K_{j}(x,y)
 = \sum_{ |\beta| \leq \mathfrak M_{\at}} \frac{(y-c_Q)^\beta}{\beta!}\,  \partial^\beta_y K_{j}(x, y)\Big|_{y=c_Q} \\
& \qquad + (\mathfrak M_{\at}+1) \sum_{ |\beta| = \mathfrak M_{\at}+1}  \frac{(y-c_Q)^\beta}{\beta!}  \int_0^1 (1-t)^{\mathfrak M_{\at}}\,\partial^\beta_y K_{j}(x,y)\Big|_{y=c_Q+t(y-c_Q)} \dd t
}
and due to vanishing moments of the atom in Definition \ref{def:hpatom}, $(iii)$, we may express the operator as
\begin{align*}
T_jf(x)   \sim  \sum_{ |\beta| = \mathfrak M_{\at}+1} \int_{Q}  \int_0 ^1 \frac{(y-c_Q)^\beta}{\beta!}\, (1-t)^{\mathfrak M_{\at}}\,
\partial^\beta_y K_{j}(x,y)\Big|_{y=c_Q+t(y-c_Q)}
f(y) \dd t\dd y.
\end{align*}

Noting that $| (y-c_Q)^\beta |\lesssim 2^{-k_Q(\mathfrak M_{\at}+1)}$ and applying the same procedure as above together with estimates \eqref{eq:kernelestimate1} and \eqref{eq:estimateybar}, we obtain
\eq{
    &\Vert T_j f \Vert _{L^p(\Rn\setminus Q_0^*)} 
    \lesssim 2^{-k_Q(\mathfrak M_{\at}+1-n/p+n)}\,2^{j(\mathfrak M_{\at}+1+ m'+n/2)}\\& = 2^{-k_Q(\mathfrak M_{\at} +1+n-n/p )}\,2^{j(\mathfrak M_{\at}+1+n-n/p)}.
}
Using this together with the fact that $\mathfrak M_{\at} > n/p-n-1$, we see that \eqref{eq:LPoperator1} is proved for $\rho=0$.\\

\textbf{Step 4 -- The case when $\boldsymbol{0<\rho\leq \frac23}$ and $\boldsymbol{k_Q\geq 1}$}\\
First, we observe that Peetre's inequality yields that,
\begin{align*} 
    &\nonumber\frac{ \,(1+2^{2j\rho}|( \nabla_\xi \varphi(x,\xi_j^\nu)-c_Q)''|^2)^{2N_1}}
    {(1+2^{2j\rho}|(\nabla \varphi(x,\xi_j^\nu)-y)''|^2)^{2N_1}}
     \frac{(1+2^{j\rho}|( \nabla_\xi \varphi(x,\xi_j^\nu)-c_Q)'|^2)^{2N_1}}{(1+2^{j\rho}|( \nabla \varphi(x,\xi_j^\nu)-y)'|^2)^{2N_1}}
\\&\leq 
    (1+2^{2j\rho}|(c_Q-y)''|^2)^{2N_1}(1+2^{j\rho}|(c_Q-y)'|^2)^{2N_1}, 
\end{align*}
which yields that, for $y\in\supp f$ and $j\geq k_Q$, 
\nm{eq:equation2}{
    g_j^\nu(x,c_Q)^{2N_1} 
\lesssim 
    g_j^\nu(x,y)^{2N_1} \sum_{\ell=0}^{8N_1} 2^{(j-k_Q)\ell}
\lesssim
    \diffcases{
        g_j^\nu(x,y)^{2N_1}\,2^{8N_1(j-k_Q)}, & j\geq k_Q,\\
        g_j^\nu(x,y)^{2N_1}, & j < k_Q.
    }
}
Let $K_j^\nu(x,y)$ be as in \eqref{eq:LPdecompkernel2}.
Then by H\"older's inequality and Minkowski's integral inequality, one has that 
\begin{align} \label{eq:111}
    &\nonumber\brkt{\int_{\Rn \setminus Q^*_\rho}\brkt{\int_{\Rn} |K_j(x,y)\,f(y)|\dd y }^p\dd x}^{1/p} 
\\&\nonumber=
    \brkt{\int_{\Rn \setminus Q^*_\rho}\brkt{\frac1{g_j^\nu(x,c_Q)^{N_1}}\int_{\Rn}\sum_{\nu=1}^{\mathscr N_j} |g_j^\nu(x,c_Q)^{N_1}\, K_j^\nu(x,y)\,f(y)|\dd y }^p\dd x}^{1/p}
\\&\nonumber\leq
    \sum_{\nu=1}^{\mathscr N_j}\norm{\frac1{g_j^\nu(\cdot,c_Q)^{N_1}}}_{L^{p'}(\Rn \setminus Q^*_\rho)}\norm{\int_{\Rn} |g_j^\nu(\cdot,c_Q)^{N_1}\, K_j^\nu(\cdot,y)\,f(y)|\dd y}_{L^2(\Rn \setminus Q^*_\rho)} 
\\&\lesssim
    \sum_{\nu=1}^{\mathscr N_j} 2^{ -j\rho (n-\kappa/2)/p'}\int_{\Rn}\norm{\int_{\Rn} |g_j^\nu(\cdot,c_Q)^{N_1}\, K_j^\nu(\cdot,y)\,f(y)|\dd y}_{L^2(\Rn \setminus Q^*_\rho)}
\end{align}
for $N_1>n/2p'$.
Here we have used that
\nm{eq:bloodyhell1}{
\int_{\Rn }\frac1{g_j^\nu(x,c_Q)^{ N_1p'}}\dd x \lesssim 2^{ -j\rho (n-\kappa/2)}
}
uniformly in $c_Q$. \\

Take $h^{\nu}_j(x,\xi)= \varphi(x,\xi) -\nabla_\xi \varphi(x,\xi_j^\nu)\cdot \xi$ as in Lemma \ref{lem:h} and define
$$b_j^{\nu,\gamma} (x,\xi) :=  2^{j\rho|\gamma|}\xi^\gamma\,a(x,2^{j\rho}\xi)\,\chi_j^\nu(\xi)\,\psi_j(2^{j\rho}\xi)\,e^{i2^{j\rho}h^{\nu}_j(x,\xi)},$$
$\mathbf t_j^\nu(x):= \nabla_\xi \varphi(x,\xi_j^\nu)$ and
\eq{
    L:= \big(1- \Delta_{\xi''}^2\big)\big(1-{2^{-j\rho}\Delta_{\xi'}}\big).
}

Note that $L $ satisfies the property
\nm{eq:LMclaim}{
    g_j^\nu(x,y)^N\, e^{i2^{j\rho}(\nabla_\xi \varphi(x,\xi_j^\nu)-y)\cdot \xi}=L^N e^{i\nabla_\xi \varphi(x,\xi_j^\nu)\cdot \xi}
}
for all $N\in \Z_{\geq0}.$\\

Using Lemma \ref{lem:chijnu} and Lemma \ref{lem:h}, Table \ref{tab:table1} shows how $\xi$-derivatives act on each term of $b_j^{\nu,\gamma}$. Using these facts we see that, for all $N\geq0$,
\eq{
    L^M b_j^{\nu,\gamma}(\mathbf t^{-1}(x+y),\xi) \in S^{m+|\gamma|}_{0,\delta}(\Rn)
}
and thus
\nm{eq:claim1}{
    \big(L^N b_j^{\nu,\gamma}\big)(\mathbf t^{-1}(x+y),2^{-j\rho}\xi) \in S^{m+|\gamma|}_{\rho,\delta}(\Rn)
}
uniformly in $y$.\\

\begin{table}
\begingroup
\setlength{\tabcolsep}{10pt} 
\renewcommand{\arraystretch}{1.5} 
\centering
\begin{tabular}{| c  |c|c |}
    \hline
    Factor &$\partial_{\xi'}^{\alpha}$ & $\partial_{\xi''}^{\beta}$ \\
    \hline\hline 
    $a(x,2^{j\rho}\xi)$&$2^{j\rho|\alpha|}2^{j(m-\rho|\alpha|)}$ &$2^{j\rho|\beta|}2^{j(m-\rho|\beta|)}$ \\
    \hline 
    $2^{j\rho|\gamma|}\xi^\gamma$ & $2^{j\rho|\gamma|}2^{j|\gamma-\alpha|(1-\rho)}$  & $2^{j\rho|\gamma|}2^{j|\gamma-\beta|(1-\rho)}$\\
    \hline 
    $\chi_j^\nu(\xi)$&$2^{j |\alpha|(3\rho/2-1)}$ & $2^{j|\beta|(\rho-1)}$\\
    \hline 
    $\psi_j(2^{j\rho}\xi)$&$2^{j|\alpha|(\rho-1)}$ &$2^{j|\beta|(\rho-1)}$ \\
    \hline 
     \raisebox{-0.7mm}{$e^{i2^{j\rho}h^{\nu}_j(x,\xi)}$}  &$\displaystyle\sum_{\gamma_1+\dots+\gamma_k =\alpha}\brkt{\prod_{\ell=1}^k 2^{j\gamma_\ell\rho}\, 2^{-j\gamma_\ell\rho/2}}$ & $\displaystyle\sum_{\gamma_1+\dots+\gamma_k =\beta} \brkt{\prod_{\ell=1}^k 2^{j\gamma_\ell\rho}\, 2^{-j\gamma_\ell\rho}}$ \\
    \hline 
\end{tabular}
\vspace{0.2cm} 
\caption{\label{tab:table1} The estimates on the size of the factors of $b_j^{\nu,\gamma} (x,\xi)$ after being acted upon by differential operators. The sums are taken over all partitions of the multi-indices $\alpha\in\mathbb Z_{\geq 0}^{\kappa}$ and
$\beta\in\mathbb Z_{\geq 0}^{n-\kappa}$. The action of $\partial_{\xi'}^{\alpha}$ adds a factor $2^{j|\alpha|/2}$ to the boundedness, while the action of $\partial_{\xi''}^{\beta}$ does not have an impact on the boundedness. 
}
\endgroup
\end{table}

Let $\widehat{f}_j^\nu$ be a smooth cut-off function such that $\widehat{f}_j^\nu(2^{j\rho}\cdot)$ is constantly equal to one on the $\xi$-support of $b_j^{\nu,\gamma}$ {and vanishes outside a compact set which is slightly larger than the aforementioned $\xi$-support}. We then have

\nm{eq:fjnuL2estimate}{
    &\|f_j^\nu\|^2_{L^2(\Rn)} 
\sim 
    \|\widehat{f}_j^\nu\|^2_{L^2(\Rn)} = 2^{jn\rho}\|\widehat{f}_j^\nu(2^{j\rho}\cdot)\|^2_{L^2(\Rn)}\lesssim 2^{jn\rho} |\Gamma^{\nu}_{j}\cap A_j|
\\&\lesssim 
    2^{jn\rho}\, 2^{jn(1-\rho)-j\rho\kappa/2}
=
    2^{j(n-\rho\kappa/2)}\nonumber
}

Now we define the $L^2$-bounded operator $S_{j,y}^{\nu,M,\gamma}$ via

\eq{
    S_{j,y}^{\nu,N,\gamma}  {f}_j^\nu (x) := 2^{-j (m'+|\gamma|)}\int_{\Rn} e^{ix\cdot\xi  }\,\big(L^N b_j^{\nu,\gamma}\big)(\mathbf t^{-1}(x+y),2^{-j\rho}\xi)\,\widehat{f}_j^\nu(\xi)\ddd\xi.
}

Now we are ready to prove the assertions of the lemma.\\

\makeatletter 
\renewcommand{\eqref}[1]{\tagform@{\ref{#1}}}
\def\maketag@@@#1{\hbox{#1}}
\textbf{Step 4.1 -- The proof of \eqref{eq:LPoperator2} when $\boldsymbol{k_Q\geq 1}$}\\
\makeatother
Note that  \eqref{eq:111} is bounded by
\eq{
    \sum_{\nu=1}^{\mathscr N_j} 2^{ -j\rho (n-\kappa/2)/p'}\int_{\Rn} \|g_j^\nu(\cdot,c_Q)^{N_1}\, K_j^\nu(\cdot,y)\|_{L^2(\Rn \setminus Q^*_\rho)}\,|f(y)|\dd y
}

Recalling that $N_2>8N_1$ and $N_3>2N_1+N_2/\rho+n$,
we use \eqref{eq:LMclaim}, \eqref{eq:keypoint}, \eqref{eq:equation2}, \eqref{eq:bloodyhell1}, \eqref{eq:fjnuL2estimate} and the fact that $j\geq k_Q,$ to deduce that
\begin{align} 
    & \int_{\Rn \setminus Q^*_\rho}  g_j^\nu(x,c_Q)^{2N_1}\,| K_j^\nu(x,y)|^2\dd x   \label{eq:equation6}
\\ & \lesssim
     \int_{\Rn \setminus Q^*_\rho} \frac{g_j^\nu(x,c_Q)^{2N_1} \,2^{2jm'}\,\|{S}_{j,y}^{\nu,N_3,0}
    {f}_j^\nu\|_{L^2(\Rn)}^2}{g_j^\nu(x,y)^{N_3-N_2/\rho+N_2/\rho}} \dd x  \nonumber 
\\ &  \lesssim
     2^{-8k_QN_1}2^{8N_1j}2^{-N_2(j-k_Q)}\,    2^{2jm'}\,\|    f_j^\nu\|_{L^2(\Rn)}^2 \int_{\Rn}   \frac{1}{g_j^\nu(x,y)^{N_3-2N_1-N_2/\rho}} \dd x \nonumber
\\&  \lesssim 
    2^{-8k_QN_1}2^{8N_1j}2^{-N_2(j-k_Q)}\,    2^{2jm'}\,    2^{j(n-\rho\kappa/2)}\,2^{ -j\rho (n-\kappa/2)} \nonumber
\\ &  =
     2^{-8k_QN_1}2^{8N_1j}2^{-N_2(j-k_Q)}\,2^{j(\kappa\rho-2\rho n+2n)}2^{-2j(n-\rho n+\kappa \rho)/p}.\nonumber
\end{align}
Using this and the fact that $\mathscr N_j= O\big(2^{j\rho\kappa/2}\big)$, \eqref{eq:111} is bounded by 
\eq{
    &2^{-4k_QN_1}2^{4N_1j}2^{-N_2(j-k_Q)}\,2^{j(\kappa\rho-2\rho n+2n)/2}2^{-j(n-\rho n+\kappa \rho)/p}2^{j\rho \kappa/2} 2^{ -j\rho (n-\kappa/2)/p'}2^{- n(1/p-1)}
\\&\leq 
    2^{-4k_QN_1- n(1/p-1)}2^{-N_2(j-k_Q)}\,2^{4N_1j-jn(1/p-1)}
\\&\lesssim
    2^{j({-4k_QN_1-N_2-n(1/p-1)})}2^{-k_Q({4N_1-N_2-n(1/p-1)})}
}
Thus proves \eqref{eq:LPoperator2} when $\rho\in(0,1]$ and $k_Q\geq 1$.\\

\makeatletter 
\renewcommand{\eqref}[1]{\tagform@{\ref{#1}}}
\def\maketag@@@#1{\hbox{#1}}
\textbf{Step 4.2 -- The proof of \eqref{eq:LPoperator1}}
\makeatother

Observe that  
\begin{align*}
    &\int_{\Rn}\norm{\int_{\Rn} |g_j^\nu(\cdot,c_Q)^{N_1} \,K_j^\nu(\cdot,y)\,f(y)|\dd y}_{L^2(\Rn \setminus Q^*_\rho)}
\\
    &=\norm{\int_{\Rn} |g_j^\nu(\cdot,c_Q)^{N_1}\, K_j^\nu(\cdot,y)-p_j^\nu(x,y-c_Q) \,f(y)|\dd y}_{L^2(\Rn \setminus Q^*_\rho)}
\\
    &\leq \int_{\Rn}\| g_j^\nu(\cdot,c_Q)^{N_1}\, K_j^\nu(\cdot,y)-p_j^\nu(x,y-c_Q)\, f(y)\|_{L^2(\Rn \setminus Q^*_\rho)}\dd y
\\
    & \lesssim \sum_{|\alpha|=\mathfrak M_{\at}+1} \int_{\Rn}|y-c_Q|^{\mathfrak M_{\at}+1}\,\| g_j^\nu(\cdot,c_Q)^{N_1} \,\partial_y^\alpha K_j^\nu(\cdot,\tilde y)\|_{L^2(\Rn \setminus Q^*_\rho)}\, |f(y)|\dd y
\end{align*}
where $p_j^\nu$ is the $\mathfrak M_{\at}$:th order Taylor polynomial of $K_j^\nu$ around $c_Q$ and $\tilde y$ is some point on the line segment between $y$ and $c_Q$. Now using \eqref{eq:LMclaim},  \eqref{eq:equation2}, \eqref{eq:bloodyhell1}, \eqref{eq:fjnuL2estimate} and $|\alpha|=\mathfrak M_{\at}+1$, we have 
\begin{align*}
    &\| g_j^\nu(\cdot,c_Q)^{N_1} \,\partial_y^\alpha K_j^\nu(\cdot,\tilde y) \|_{L^2(\Rn \setminus Q^*_\rho)}^2 
\\&=
    \int_{\Rn \setminus Q^*_\rho}  g_j^\nu(x,c_Q)^{2N_1}\,| \partial_y^\alpha K_j^\nu(x,\tilde y)|^2\dd x
\\&\lesssim 
    2^{2j (m'+\mathfrak M_{\at}+1)}\int_{\Rn \setminus Q^*_\rho}  \frac{g_j^\nu(x,c_Q)^{2N_1}\,\|{S}_{j,\tilde y}^{\nu,3N_1,\alpha}
    {f}_j^\nu\|_{L^2(\Rn)}^2}{g_j^\nu(x,y)^{3N_1}}\dd x
\\&\lesssim 
    2^{2j (m'+\mathfrak M_{\at}+1)}\|\,    f_j^\nu\|_{L^2(\Rn)}^2\int_{\Rn}  \frac1 {g_j^\nu(x,y)^{3N_1-2N_1}}\dd x
\\&\lesssim 
    2^{2j (m'+\mathfrak M_{\at}+1)}\,2^{j(n-\rho\kappa/2)}\,2^{ -j\rho (n-\kappa/2)}
\\&\lesssim 
    2^{2j(\mathfrak M_{\at}+1)}\,2^{j(\kappa\rho-2\rho n+2n)}2^{-2j(n-\rho n+\kappa \rho)/p}
\end{align*}
Using this and the fact that $\mathscr N_j= O\big(2^{j\rho\kappa/2}\big)$, \eqref{eq:111} is bounded by 
\begin{align*}
    &\sum_{\nu=1}^{\mathscr N_j} 2^{ -j\rho (n-\kappa/2)/p'}\int_{\Rn} \|g_j^\nu(\cdot,c_Q)^{N_1} \,K_j^\nu(\cdot,y)\|_{L^2(\Rn \setminus Q^*_\rho)}\,|f(y)|\dd y
\\&\lesssim 
    \sum_{\nu=1}^{\mathscr N_j} 2^{ -j\rho (n-\kappa/2)/p'}\sum_{|\alpha|=\mathfrak M_{\at}+1} \int_{\Rn}|y-c_Q|^{\mathfrak M_{\at}+1}
\\&\qquad\qquad\qquad\times
    \| g_j^\nu(\cdot,c_Q)^{N_1}\, \partial_y^\alpha K_j^\nu(\cdot,\tilde y)\|_{L^2(\Rn \setminus Q^*_\rho)}\, |f(y)|\dd y
\\&\lesssim 
    2^{j\rho\kappa/2}2^{ -j\rho (n-\kappa/2)/p'} 2^{-k_Q(\mathfrak M_{\at}+1)} 2^{j(\mathfrak M_{\at}+1)}\,2^{j(\kappa\rho-2\rho n+2n)/2}2^{-j(n-\rho n+\kappa \rho)/p}
\\&\leq
    2^{(\mathfrak M_{\at}+1)(j-k_Q)}2^{-jn(1/p-1)}2^{k n(1/p-1)}
\end{align*}
Observe that $\mathfrak M_{\at}+1>n(1/p-1)$, so \eqref{eq:LPoperator1} is proved when $\rho\in(0,1]$ and $k_Q\geq 1$.\\

\makeatletter 
\renewcommand{\eqref}[1]{\tagform@{\ref{#1}}}
\def\maketag@@@#1{\hbox{#1}}
\textbf{Step 5 -- The proof of \eqref{eq:LPoperator2} when $\boldsymbol{k_Q\leq 0}$}\\
\makeatother
Take $\varepsilon:= -m'-n/2$. Observe that the condition on $p$ yields that $\varepsilon>0$. Then we have
\eq{
    &\| T_j  f\|_{L^p(\Rn \setminus Q^*_\rho)}
=
    \norm{\frac{g(x,c_Q)^{N_1}}{g(x,c_Q)^{N_1}}\, T_j  f}_{L^p(\Rn \setminus Q^*_\rho)} 
\\& \lesssim
    2^{-j\varepsilon}\norm{\frac{1}{g(x,c_Q)^{N_1}}}_{L^{p'}(\Rn \setminus Q^*_\rho)} \|g(x,c_Q)^{N_1}\, 2^{j\varepsilon}T_j  f\|_{L^2(\Rn \setminus Q^*_\rho)} 
}
Now
\nm{eq:equation7}{
    &\|g(x,c_Q)^{N_1}\, 2^{j\varepsilon}T_j  f\|_{L^2(\Rn \setminus Q^*_\rho)}
\leq
    \int_{\Rn}\|g(x,c_Q)^{N_1}\, 2^{j\varepsilon}K_j(x,y)\|_{L^2(\Rn \setminus Q^*_\rho)} f(y)\dd y 
\\& \lesssim
    2^{k_Q(1/p-1)}\|g(x,y)^{N_1}\, 2^{j\varepsilon}K_j(x,y)\|_{L^2(\Rn \setminus Q^*_\rho)}\nonumber
}
where we have used \eqref{eq:estimateybar}. Now since $k_Q\leq 0$ by assumption, \eqref{eq:newone} yields that \eqref{eq:equation7} is uniformly bounded. Hence
\eq{
    \| T_j  f\|_{L^p(\Rn \setminus Q^*_\rho)}\lesssim 2^{-j\varepsilon}
}
so \eqref{eq:LPoperator2} is proved when $k_Q\leq 0$.\\

\textbf{Step 6 -- The adjoint operator}\\
The corresponding proof for the adjoint $(T_j)^*$ is similar to the one above with few modifications. First, regarding the $L^2$ boundedness in Lemma \ref{lemma:HLS-lemma}, estimate \eqref{eq:L2toLq} in that Lemma has to be replaced by \eqref{eq:L2toLqadj}. Second, the $x$- and $y$-dependencies of the kernel are reversed. This means the following replacements:
\eq{
    &\varphi(x,\xi) \longrightarrow x,\\
    &y \longrightarrow \varphi(y,\xi),\\
    &c_Q \longrightarrow \varphi(c_Q,\xi).
}
This implies that \eqref{eq:keypoint2} is used instead of \eqref{eq:keypoint} and when applying $\beta$ derivatives in the $y$-variable the $y$-dependence in both arguments has to be taken into consideration. Otherwise, the proof remains the same. 
\end{proof}

\section{Fourier integral operator estimates}\label{Sec:FIO estimates}

This section is devoted to the estimates for Fourier integral operators on various function spaces. It begins with showing a global $h^p\to L^p$ result for FIOs with amplitudes in $S^m_{\rho,\delta}$ for $\rho\in[0,1],\,\delta\in[0,1)$ and $0<p<\infty$. We then lift this result to both Triebel-Lizorkin and Besov-Lipschitz spaces in the forthcoming sections.

\begin{Prop}\label{prop:hptoLpboundednessFIO}
Let $T_a^\varphi$ be an $\mathrm{FIO}$ with an amplitude $a\in S^{m_c(p)}_{\rho,\delta}(\Rn)$ for $\rho\in[0,1],\,\delta\in[0,1)$ and $m_c(p)$ in \eqref{eq:criticaldecay}. Let $\varphi \in \Phi^2$ be an $\mathrm{SND}$ phase function of rank $\kappa$ satisfying the conditions in \emph{Definition \ref{def:FIO}}. Then  $T_{a}^\varphi:h^p(\Rn)\rightarrow L^p(\Rn)$ for $0<p<\infty$ and $T_{a}^\varphi:L^\infty(\Rn)\rightarrow \mathrm{bmo}(\Rn)$ when $p=\infty$.
\end{Prop}
\begin{proof}
The case $\frac23<\rho\leq 1$ follows by first using the decomposition of the FIO as outlined in \cite{SSS}, and then following the method below for operators for which $0\leq \rho\leq \frac{2}{3}$ . Thus we confine ourselves to the proof in this case.\\

The proof is carried out as follows. We start with proving the $h^p\to L^p$ boundedness for $T_a^\varphi$ when $0<p<\frac{2n-2n\rho+2\kappa\rho}{2n-n\rho+\kappa\rho}$. To do this it is enough to show that $\|T_a^\varphi \at\|_{L^p(\Rn)} \lesssim 1$ uniformly in the $h^p$-atom $\at$. Then we proceed with proving $h^p\to L^p$ boundedness for $(T_a^\varphi)^*$. To finish the proof, we do Riesz-Thorin interpolation with the $L^2$ boundedness of $T_a^\varphi.$\\

The special case when $\rho=1$ and $\kappa=0$ is treated later in Remark \ref{rem:pseudoTLremark}.\\

In what follows, let $\at(x)$ be an $h^p$-atom supported in a cube $Q$. Moreover let $$m'= -\big(\kappa+(n-\kappa)(1-\rho)\big)\Big|\frac{1}{p}-\frac{1}{2}\Big|$$ as before. We split up $\|T_a^\varphi \at\|_{L^p(\Rn)}$ into two pieces.
$$\|T_a^\varphi \at\|_{L^p(\Rn)} = \|T_a^\varphi \at\|_{L^p(Q^*_\rho)}+\|T_a^\varphi \at\|_{L^p(\Rn\setminus Q^*_\rho)}$$
where $Q^*_\rho$ is given in Definition \ref{def:influenceset}. \\

We begin by estimating $\|T_a^\varphi \at\|_{L^p(Q^*_\rho)}$. We split up this process into two substeps.\\

\textbf{Step 1 -- The case when $\boldsymbol{k_Q\geq 1}$ and $\boldsymbol{0<\rho\leq \frac23}$ }\\

In this step we take $\mathfrak b := 2^{-k_Qn(1/p+m'/n-1/2)}\at$. We claim that 
\eq{
    \|\mathfrak b\|_{h^{2n/(n-2 m')}} \lesssim 1
}
Observe that $$ 0<p<\frac{2n-2n\rho+2\kappa\rho}{2n-n\rho+\kappa\rho}\Rightarrow \frac{2n}{n-2m'}<1.$$

The claim follows if we show that $\mathfrak b$ is an $h^{2n/(n-2 m')}$-atom, since then the sum in \eqref{eq:hpsumdefinition} contains one term only. Hence it is enough to show that $\mathfrak b$ satisfies the three conditions in Definition \ref{def:hpatom} for $h^{2n/(n-2 m')}$. $(i)$ is obvious. $(ii)$ follows since 
\eq{
    &\sup_{x\in Q}|\mathfrak b| = 2^{-k_Qn(1/p+m'/n-1/2)}\,\sup_{x\in Q}|\at|\lesssim 2^{-k_Qn(1/p+m'/n-1/2)} 2^{k_Qn/p}
\\&= 
    2^{-k_Qn(m'/n-1/2)} = 2^{k_Qn/(2n/(n-2 m'))}.
}
$(iii)$ follows since $\frac{2n}{n-2m'}>p$ which implies that $\mathfrak M_{\mathfrak b}<\mathfrak M_\at$.\\

Now to finish Step 1, H\"older's inequality and the assertion in Lemma \ref{lem:sizeofQstar} $(i)$ and Lemma \ref{lemma:HLS-lemma}, yield that
\eq{
    &\int_{Q^*_\rho}|T_a^\varphi \at(x)|^{p}\dd x
\leq 
    |Q^*_\rho|^{1-\frac p2}\brkt{\int_{Q^*_\rho}|T_a^\varphi  \at(x)|^{2}\dd x}^{\frac p2}
\lesssim 
    2^{-k_Q\rho(n-\kappa)(1-\frac p2)}\Vert T_a^\varphi  \at\Vert_{L^{2}(\Rn)}^{p}
\\&\lesssim 
    2^{-k_Q\rho(n-\kappa)(1-\frac p2)}\|\at\|_{h^{2n/(n-2 m')}(\Rn)}^{p} = \|\mathfrak b\|_{h^{2n/(n-2 m')}(\Rn)}^{p} 
\lesssim
    1
}

\textbf{Step 2 -- The case when $\boldsymbol{k_Q\leq 0}$ or $\boldsymbol{\rho=0}$ }\\
In this case $Q^*_\rho = 2\sqrt{n}Q$, hence
\eq{
    &\int_{Q^*_\rho}|T_a^\varphi \at(x)|^{p}\dd x\leq |Q^*_\rho|^{1-\frac p2}\brkt{\int_{Q^*_\rho}|T_a^\varphi  \at(x)|^{2}\dd x}^{\frac p2} 
\\&\lesssim
    2^{-k_Qn(1-p/2)} \| \at \|_{L^2}^p \lesssim 2^{-k_Qn(1-p/2)} 2^{nk_Q} 2^{-k_Qp/2}=1.
}\hspace*{1cm}

\textbf{Step 3 -- Estimates of $\boldsymbol{\|T_a^\varphi \at\|_{L^p(\Rn\setminus Q^*_\rho)}}$}\\

Using the Littlewood-Paley partition of unity that was introduced in Definition \ref{def:LP} we can write
\begin{equation}\label{eq:defn of sss pieces of T}
 T_a^\varphi  =\sum_{j=1}^{\infty}T_a^\varphi\,\psi_j(D)=:\sum_{j=1}^{\infty}T_{j}.
\end{equation}

Now to deal with the integral $\displaystyle \int_{\Rn\setminus Q^*_\rho}|T_a^\varphi \at(x)|^{p}\dd x,$ using the notation in \eqref{eq:defn of sss pieces of T}, we observe that
\eq{
    &\int_{\Rn\setminus Q^*_\rho}|T_a^\varphi  \at(x)|^p \, \dd x\leq \sum_{j=1}^{\max\{k_Q,0\}}\int_{\Rn\setminus Q^*_\rho}|T_{j}\at(x)|^p \, \dd x
\\&\qquad+
    \sum_{j=\max\{k_Q+1,1\}}^\infty\int_{\Rn\setminus Q^*_\rho}|T_{j}\at(x)|^p \,\dd x.
}

We use Lemma \ref{lem:bstarcestimate} and the properties of the atom $\at$ to conclude Step 3.\\

\textbf{Step 4 -- Estimates of the adjoint $\boldsymbol{(T_a^\varphi)^*}$}\\
The proof of the $h^p\to L^p$ boundedness for $(T_a^\varphi)^*$ is similar to the one of $T_a^\varphi$ itself. The only difference is that Lemma \ref{lem:sizeofQstar} $(ii)$ is used instead of $(i).$ This in particular implies the $h^p\to L^p$ boundedness of $(T_a^\varphi)^*$. \\

\textbf{Step 5 -- Interpolation with $\boldsymbol{L^2(\Rn)}$}\\
Now that we have boundedness from $h^p(\Rn)$ to $L^p(\Rn)$ for both $T_a^\varphi$ itself and its adjoint when $0<p<\frac{2n-2n\rho+2\kappa\rho}{2n-n\rho+\kappa\rho}$, as well as $L^2$ boundedness, we can use a standard Riesz-Thorin interpolation argument to conclude that $T_a^\varphi$ is bounded from $h^p(\Rn)$ to $L^p(\Rn)$ for all $0<p<\infty$ and from $L^\infty(\Rn)$ to $\mathrm{bmo}(\Rn)$ when $p=\infty$.
\end{proof}

\begin{Rem}
    Note that using \emph{Proposition \ref{prop:hptoLpboundednessFIO}} together with \emph{Theorem \ref{thm:left composition with pseudo}} one could also lift the $h^p\to L^p$ boundedness result to a $h^p\to h^p$ boundedness, see the proof of \cite[Proposition 6.2]{IRS} for details. Note also that this result will also yield the $h^{1}-h^{1}$ boundedness of $(T_a^\varphi)^*$ which by duality yields the boundedness of $T_a^\varphi$ from $\mathrm{bmo}(\Rl^n) \to \mathrm{bmo}(\Rl^n).$
\end{Rem}

\subsection{Triebel-Lizorkin estimates related to classical and exotic amplitudes} \label{subsec:FIO_classic_exotic}

In this section, we state and prove the global boundedness of Fourier integral operators with classical and exotic amplitudes on Triebel-Lizorkin spaces. To this end, we define a molecular representation of the Triebel-Lizorkin spaces. Similar to the atomic representation of the Hardy spaces, this can be used to prove boundedness results.

\begin{Def}[Notation]
Let $\mathcal D$ be the set of all dyadic cubes in $\Rn$ and define the following sets:
\begin{enumerate}
    \item[$(i)$] $\mathcal D_j :=\set{Q\in \mathcal D: k_Q= j} $
    \item[$(ii)$] $\mathcal D_+:= \set{Q\in \mathcal D: k_Q\geq0}$
    \item[$(iii)$] $\mathcal D_j(Q):= \{\tilde  Q\in \mathcal D_{j}:\tilde Q\subseteq Q\} $.
\end{enumerate}
\end{Def}

Observe that for $j<k_Q$, $D_j(Q)=\emptyset.$\\

We start with defining a space of sequences that is easier to handle than the Triebel-Lizorkin spaces themselves.

\begin{Def}\label{def:TLatomseq1}For a sequence of complex numbers $b = \set{b_Q}_{\substack{ Q\in\mathcal D\\
l(Q)\leq 1}}$
we define
\eq{
g^{s,q}(b)(x) := \brkt{\sum_{ \tilde Q\in\mathcal D_+}
\big(2^{k_{\tilde Q}(s+n/2)}|b_{\tilde Q}|\chi_{\tilde Q}(x)\big)^q}^{1/q}.
}
We say that $b\in f^s_{p,q}$ is
\eq{
\| b\|_{f^s_{p,q}} := \|g^{s,q}(b)\|_{L^p(\Rn)}<\infty.
}
\end{Def}
In what follows, take $\check\Psi^{\tilde Q}(x) := 2^{-nk_{\tilde Q}/2} \check\Psi_{k_{\tilde Q}}(x- c_{\tilde Q})$, where $\Psi_{k_{\tilde Q}}$ as given in Definition \ref{def:LP} is the usual Littlewood-Paley piece. The following Lemma is a corollary of \cite[Theorem II B]{FJ:phi-transform}.
\begin{Lem}\label{lem:TLatomseq1}
Suppose $0 < p < \infty,$ $ 0 < q \leq \infty$, $s \in \Rl.$ For any sequence $b= \set{b_{\tilde Q}}_{\tilde Q\in\mathcal D}$ of complex numbers satisfying
$\|b\|_{
f^s_{p,q}} < \infty,$ one has
\eq{
    f(x) := \sum_{\tilde Q\in \mathcal D_+}b_{\tilde Q}\check\Psi^{\tilde Q}(x)
}
belongs to $F^s_{p,q}(\Rn)$ and
\eq{
\|f\|_{F^s_{p,q}(\Rn)}\lesssim 
\|b\|_{f^s_{p,q}}.
}
\end{Lem}
\begin{proof}
    The functions $\Psi^{\tilde Q}$ satisfy the conditions for being "smooth molecules", i.e. there exists integers $M,N,K>0$ such that the functions satisfy
    \eq{
        \diffcases{\displaystyle\int_{\Rl^n} x^\gamma \check{\Psi}^{\tilde Q} \dd x= 0,& |\gamma|\leq N\\
        \displaystyle|\partial^\gamma \check{\Psi}^{\tilde Q}(x)| \lesssim \frac{|\tilde Q|^{-1/2-|\gamma|/n}}{1+l_{\tilde Q}^{-1}|x-c_{\tilde Q}|^{M+|\gamma|}},  & |\gamma|\leq K.}
    }
    The rest follows from \cite[Theorem II B]{FJ:phi-transform}. Observe that the proof of \cite[Theorem II B]{FJ:phi-transform} can easily be modified to inhomogeneous Triebel-Lizorkin spaces, see \cite[Chapter 12]{FJ:discrete-transform} for details.  
\end{proof}

We now discuss the converse of Lemma \ref{lem:TLatomseq1}, namely when a Triebel-Lizorkin function can be expressed in terms of molecules and so-called "$\infty$-atoms", here denoted $b_{\iota,\tilde  Q}$.
\begin{Lem}\label{lem:TLatomseq}
Suppose $0 < p \leq 1,$ $p \leq q \leq \infty$. Every $f\in F^0_{p,q}(\Rn)$ has an atomic decomposition 
\eq{
    f(x)=\sum_{\iota=1}^\infty \lambda_\iota  \sum_{ \tilde Q\in \mathcal D_+}b_{\iota,\tilde  Q}\check\Psi^{\tilde  Q}(x), \qquad b_{\iota,\tilde  Q}\in 
    f^0_{\infty,q},
}
where $ \set{b_{\iota,\tilde  Q}}_{\tilde Q\in\mathcal D_+}=:b_\iota$ satisfies 
\eq{
\|b_\iota\|_{f^0_{\infty,q}} \leq 2^{k_{\tilde Q}n/p}. 
}
Moreover
\eq{
\|f\|_{F^0_{p,q}(\Rn)}\sim\inf_{\set{\lambda_\iota}}\brkt{\sum_{\iota=1}^\infty |\lambda_\iota|^p}^{1/p}
}

\end{Lem}

\begin{proof}
This follows by combining \cite[Theorem II A i]{FJ:phi-transform} and \cite[Theorem 1]{HanPalWeiss}. Observe that once again this can easily be modified to inhomogeneous Triebel-Lizorkin spaces, see \cite[Chapter 12]{FJ:discrete-transform} for details. 
\end{proof}
Next, we define the analog of the Hardy space atoms, which will be used to prove the boundedness results of our FIO's.
\begin{Lem}\label{lem:TLatomseq4}
Let $Q \in \mathcal D$, $0<p\leq 1$, $0<q\leq\infty$, $j\geq 1$ and $b = \{b_{\tilde Q}\}_{\tilde Q\in\mathcal D_+}\in f^0_{\infty,q}$ with $\|b\|_{f^0_{\infty,q}} \leq 2^{k_Qn/p}$. Define
\nm{eq:RQatom}{
    R_{Q,j}(x):=\sum_{\tilde Q\in\mathcal D_j( Q)}b_{\tilde Q}\check\Psi^{\tilde Q}(x).
}
Then $R_{Q,j}(x)=0$ for $j< k_{\tilde Q}$. Moreover, for $0<p<\infty$,
\eq{
\|R_{Q,j}\|_{L^{2}(\Rn)}\lesssim 2^{k_{ Q}n(1/p-1/2)}.
}

\end{Lem}
\begin{proof}
We first claim that
\nm{eq:TLseqclaim}{
    |b_{\tilde Q}|\leq 2^{k_Qn/p-nk_{\tilde Q}/2}.
}
Indeed, using that $b \in f^0_{\infty,q}$ we have
\eq{
    2^{nk_{\tilde Q}/2}  |b_{\tilde Q}|\chi_{\tilde Q}(x) \leq \sup_{x\in\Rn} \brkt{\sum_{\tilde Q\in \mathcal D_+} (2^{nk_{\tilde Q}/2}\,  |b_{\tilde Q}|\,\chi_{\tilde Q}(x))^q}^{1/q} = \|b\|_{f^0_{\infty,q}} \leq 2^{k_Qn/p}.
}
Now taking supremum in $x$ on the left-hand side yields the claim.\\

Using Lemma \ref{lem:TLatomseq1} we have
\eq{
    &\|R_{Q,j}\|_{L^{2}(\Rn)} 
=
    \|R_{Q,j}\|_{F^0_{2,2}}
\lesssim
    \norm{\brkt{\sum_{\tilde Q\in\mathcal D_j(Q)} \brkt{2^{nk_{\tilde Q}/2}\,|b_{\tilde Q}|\,\chi_{\tilde Q}(x)}^2}^{1/2}}_{L^{2}(\Rn)}
\\&\lesssim 
    2^{k_Qn/p}\norm{\brkt{\sum_{\tilde Q\in\mathcal D_j(Q)} \chi_{\tilde Q}(x)^2}^{1/2}}_{L^{2}(\Rn)} 
=
    2^{k_Qn/p}\|\chi_Q\|_{L^{2}(\Rn)} = 2^{k_{Q} n(1/p-1/{2})}
}
\end{proof}

Now we start with the lifting results for FIO's to Triebel-Lizorkin spaces. In order to do that, we need to introduce a partition of unity to $R_{Q,j}$. We estimate the pieces separately. The first estimate is given in Lemma \ref{lem:TLatomseq5}.

\begin{Lem}\label{lem:TLatomseq5}
Suppose that $0<p<\frac{2n-2n\rho+2\kappa\rho}{2n-n\rho+\kappa\rho}$ and $q\geq p$. Let $a (x,\xi)\in S^{m_c(p)}_{\rho, \delta}(\Rn)$ be supported outside a neighborhood of the origin and $\varphi \in \Phi^2$ be \emph{SND} of rank $\kappa$. Then
\nm{eq:equation3}{
\| T_a^\varphi \psi_j(D)  \chi_{\Rn \setminus2\sqrt n  Q}R_{Q,j}\|_{L^p(\Rn )} \lesssim 2^{n(k_{Q}-j)}.
}
for $j\geq k_{ Q}$.
\end{Lem}

\begin{proof}
We make two claims: for $N_1$ and $N_2$ positive integers
\begin{align}
    &|\check\Psi^{\tilde Q}(y)|\lesssim 2^{-nk_{\tilde Q}/2}2^{jn-jN_1}|y-c_{ Q}|^{-N_1}\label{eq:equation4}
\\
    &|\check\Psi_j(x-y)|\lesssim  \frac{2^{jn+jN_2}|y-c_{ Q}|^{N_2}}{\langle2^j( x-c_{ Q})\rangle^{N_2}} 
    \label{eq:equation5},\qquad y\in \Rn \setminus2\sqrt n Q
\end{align}

We prove \eqref{eq:equation4}. Observe that
$$
|y-c_{ Q}|\leq |y-c_{\tilde Q}|+|c_{\tilde Q}-c_{ Q}| \leq |y-c_{\tilde Q}|+\sqrt n l_{ Q} \leq |y-c_{\tilde Q}|+\frac12 |y-c_{ Q}|
$$
so 
\eq{
    |\check\Psi^{\tilde Q}(y)| = 2^{-nk_{\tilde Q}/2} |\check\Psi_j(x-c_{\tilde Q})|
\lesssim 
    \frac{2^{-nk_{\tilde Q}/2}2^{jn-jN_1}}{|y-c_{\tilde Q}|^{N_1}} 
\lesssim
    \frac{2^{-nk_{\tilde Q}/2}2^{jn-jN_1}}{|y-c_{ Q}|^{N_1}}
}
We turn to \eqref{eq:equation5}. By Peetre's inequality we obtain
\eq{
    &\langle 2^j(x-y)\rangle^{ -1 } = \langle 2^j(x-c_{ Q})+2^j(c_{ Q}-y)\rangle^{ -1 } 
\leq 
    \langle 2^j(x-c_{ Q})\rangle^{ -1 }(1+2^{j}|y-c_{ Q}|)
\\&\lesssim
    \langle 2^j(x-c_{ Q})\rangle^{ -1 }2^{j}|y-c_{ Q}|^{ } (2^{k_Q-j}+1) 
\lesssim 
    \langle 2^j(x-c_{ Q})\rangle^{ -1}2^{j}|y-c_{ Q}| 
}
Hence we have
\eq{
    &|\check\Psi_j(x-y)|\lesssim \frac{2^{jn}}{\langle 2^j(x-y)\rangle^{N_2}} \lesssim \frac{2^{jn+jN_2}|y-c_{ Q}|^{N_2}}{\langle2^j( x-c_{ Q})\rangle^{N_2}},
}
which is \eqref{eq:equation5}.\\

Now using Proposition \ref{prop:hptoLpboundednessFIO}, \eqref{eq:equation4}, \eqref{eq:equation5} and \eqref{eq:TLseqclaim} we obtain 
\begin{align*}
    &\| T_j  \chi_{\Rn \setminus2\sqrt n Q}R_{Q,j}\|_{L^p(\Rn )}^p\lesssim     \| \Psi_j(D) \chi_{\Rn \setminus2\sqrt n Q}R_{Q,j}\|_{L^p(\Rn)}^p 
\\&\leq
    \int_{\Rn} \brkt{ \int_{\Rn \setminus2\sqrt n Q} |\check\Psi_j(x-y)| \sum_{\tilde Q\in\mathcal D_j(Q)} |b_{\tilde Q}|\,|\check\Psi^{\tilde Q}(y)| \dd y }^p \dd x
\\&=
    2^{jp(N_2-N_1+2n)}\int_{\Rn} \brkt{ \int_{\Rn \setminus2\sqrt n Q}  \sum_{\tilde Q\in\mathcal D_j(Q)}|b_{\tilde Q}|2^{-nk_{\tilde Q}/2} \frac{|y-c_{ Q}|^{N_2-N_1}}{\langle2^j( x-c_{ Q})\rangle^{N_2}}\dd y }^p \dd x
\\&\lesssim
    2^{jp(N_2-N_1+2n)}\brkt{\sum_{\tilde Q\in\mathcal D_j(Q)}|b_{\tilde Q}|2^{-nk_{\tilde Q}/2}}^p \int_{\Rn} \brkt{ \int_{\Rn \setminus2\sqrt n Q}  \frac{|y-c_{ Q}|^{N_2-N_1}}{\langle2^j( x-c_{ Q})\rangle^{N_2}}\dd y }^p \dd x
\\&\lesssim
    2^{jp(N_2-N_1+2n)-jn}2^{-k_{ Q}p(n+N_2-N_1)}\brkt{\sum_{\tilde Q\in\mathcal D_j(Q)}2^{k_Qn/p-k_{\tilde Q}n}}^p
\\&\lesssim
    2^{jp(N_2-N_1+2n)-jn}2^{-k_{ Q}p(n+N_2-N_1)} 2^{np(j-k_{ Q})} 2^{-jpn} 2^{k_{ Q}n} = 2^{-n(j-k_{ Q})}
\end{align*}
for $N_2>n/p$ and $N_1=n+N_2+1$, since there are $O\big(2^{n(j-k_{ Q})}\big)$ elements in ${Q\in\mathcal D_j( Q)}$.
\end{proof}

\begin{Lem}\label{lem:local and global nonendpoint TL}
Let $a(x,\xi)\in S_{\rho,\delta}^{m}(\Rn)$ and that $T_a^\varphi$ is an oscillatory integral operator that is bounded from $F_{p,p}^s(\Rn)$ to $F_{p,p}^s(\Rn)$. Then $T_\sigma^\varphi$ is bounded from $F_{p,q}^{s}(\Rn)$ to $F_{p,q}^s(\Rn)$ for $\sigma(x,\xi)\in S_{\rho,\delta}^{m-\varepsilon}(\Rn)$ where $\varepsilon>0$ is arbitrary.
\end{Lem}
\begin{proof}
	Using the embedding \eqref{embedding of TL}, equality \eqref{equality of TL and BL}, and finally the fact that $(1-\Delta)^{- \varepsilon/4}$ is an isomorphism from $F_{p,q}^{s}(\Rn)$  to $F_{p,q}^{s+\varepsilon/2}(\Rn)$,
we have that
\eq{
&\Vert T_\sigma^\varphi f\Vert_{F^{s}_{p,q}(\Rn)}=\Vert T_\sigma^\varphi (1-\Delta)^{\varepsilon/2}(1-\Delta)^{-\varepsilon/2}f\Vert_{F^{s}_{p,q}(\Rn)}\\ 
&\lesssim \Vert T_\sigma^\varphi (1-\Delta)^{\varepsilon/2}(1-\Delta)^{-\varepsilon/2}f\Vert_{F^{s+\varepsilon/2}_{p,p}(\Rn)}\lesssim \Vert (1-\Delta)^{-\varepsilon/2}f\Vert_{F^{s+\varepsilon/2}_{p,p}(\Rn)} \\ 
& \lesssim \Vert (1-\Delta)^{-\varepsilon/2}f\Vert_{F^{s+\varepsilon}_{p,q}(\Rn)}\simeq \Vert f\Vert_{F^{s}_{p,q}(\Rn)}.
}
\end{proof}

\begin{Prop}\label{prop:TLp<1}
Suppose that $0<p<\frac{2n-2n\rho+2\kappa\rho}{2n-n\rho+\kappa\rho}$ and $q\geq p$. Let $a (x,\xi)\in S^{m_c(p)}_{\rho, \delta}(\Rn)$ be supported outside a neighborhood of the origin and $\varphi \in \Phi^2$ be \emph{SND} of rank $\kappa$ for all $(x, \xi) \in \supp a$. Then the \emph{FIO} $T_a^\varphi$ is bounded from $F_{p,q}^{s}(\Rn)$ to $F_{p,q}^s(\Rn)$.
\end{Prop}

\begin{proof}
Note that for $\frac23<\rho\leq 1$, the proof has to be modified in accordance with the decomposition in \cite{SSS}, and thereafter the same method as in the case of $0<\rho\leq \frac23$ can be applied. Thus, we confine ourselves to the proof in the case of $0<\rho\leq \frac23$, which is the new contribution here, and leave the details of the case $\rho> \frac23$ to the interested reader.\\

Observe that it is enough to show the result for $s=0$ and by the inclusion $F^0_{p,p}\xhookrightarrow{}F^0_{p,q}$ it is enough to show that $T_a^\varphi$ is bounded from $F_{p,q}^{0}(\Rn)$ to $F_{p,p}^0(\Rn).$ \\

Compose a Littlewood-Paley piece $\psi_j(D)$ with $T_a^\varphi$ and apply Theorem \ref{thm:left composition with pseudo}. This yields 
\eq{
    \psi_j(D) T_a^\varphi f(x) = T_a^\varphi \psi_j(\nabla_x\varphi(x,D))f(x) +\sum_{\alpha\leq M} 2^{-j\varepsilon}Rf(x),
}
where $\varepsilon>0$ stems from Theorem \ref{thm:left composition with pseudo} and $R$ is an operator of better decay, therefore by Lemma \ref{lem:local and global nonendpoint TL} we have
\eq{
    \norm{\sum_{j=0}^\infty 2^{-j\varepsilon}|Rf(x)|}_{L^p(\Rn)} \sim \|Rf\|_{L^p(\Rn)} \lesssim \| f\|_{F^{-(1-\max\{\delta,1/2\}-\varepsilon)}_{p,2}(\Rn)}\lesssim \| f\|_{F^0_{p,q}(\Rn)}.
}
So from now on we will only consider the $F^0_{p,q}\to L^p(\ell^q)$ boundedness of the first term.\\

Let $\mathbf t(\xi)=\nabla_x\varphi(x,\xi)$ and $\eta:\Rn\to\Rn $ be diffeomorphisms such that $\eta(y)\cdot \mathbf t^{-1}(\xi) = y\cdot \xi$. Then we have
\eq{
    &T_a^\varphi \psi_j(\nabla_x\varphi(x,D))f(x) = \iint_{\Rn\times\Rn} e^{i\varphi(x,\xi)-iy\cdot\xi} \,\psi_j(\mathbf t(\xi))\, a(x,\xi) f(y)\dd y\ddd \xi  
\\&=\frac1{|\det(\nabla t)|}
    \iint_{\Rn\times\Rn} e^{i\varphi(x,\mathbf t^{-1}(\xi))-iy\cdot t^{-1}(\xi)}\, \psi_j(\xi)\, a(x,\mathbf t^{-1}(\xi)) f(y)\dd y\ddd \xi  
\\&= 
    \frac{|\det(\nabla \eta)|}{|\det(\nabla t)|}\iint_{\Rn\times\Rn} e^{i\varphi(x,\mathbf t^{-1}(\xi))-iy\cdot \xi} \,\psi_j(\xi)\, a(x,\mathbf t^{-1}(\xi)) f(\eta(y))\dd y \ddd \xi 
\\&=:
    T_j (f\circ\eta)(x)
}
Observe that by Theorem \ref{thm:invariance thm} it is enough to consider $T_j f$ from now on.\\

By Lemma \ref{lem:TLatomseq},  $f\in F^0_{p,q}$ has an atomic decomposition 
\eq{
    f(x)=\sum_{\iota=1}^\infty \lambda_\iota  \sum_{ \tilde Q\in \mathcal D_+}b_{\iota,\tilde  Q}\,\check\Psi^{\tilde Q}(x), \qquad b_{\iota,\tilde Q}\in f^0_{p,q}\cap f^0_{\infty,q},
}
such that
\eq{
    \|f\|_{F^0_{p,q}(\Rn)}\sim\inf_{\set{\lambda_\iota}}\brkt{\sum_{\iota=1}^\infty |\lambda_\iota|^p}^{1/p}.
}

Since $\supp \psi_j\subset \supp \Psi^{\tilde Q}$ it is enough to consider $\tilde Q\in \mathcal D_j$ and hence
\eq{
&\norm{\brkt{\sum_{j=1}^\infty|T_jf|^p}^{1/p}}_{L^p(\Rn)} = \norm{\brkt{\sum_{j=1}^\infty \abs{\sum_{\iota=1}^\infty\lambda_\iota T_j \sum_{\tilde Q\in\mathcal D_j}b_{\iota,\tilde Q}\,\check\Psi^{\tilde Q}(x)}^p}^{1/p}}_{L^p(\Rn)} 
\\&\lesssim
\brkt{\sum_{\iota=1}^\infty|\lambda_\iota|^p\int_{\Rn}  \sum_{j=1}^\infty\abs{ T_j \sum_{Q\in\mathcal D_j}b_{\iota,\tilde Q}\,\check\Psi^{\tilde Q}(x)}^p\dd x}^{1/p} 
\\&\lesssim
\brkt{\sum_{\iota=1}^\infty|\lambda_\iota|^p}^{1/p}\sup_{\iota\in \mathbb Z_{>0}}\brkt{ \sum_{j=1}^\infty\norm{ T_j \sum_{\tilde Q\in\mathcal D_j}b_{\iota,\tilde Q}\,\check\Psi^{\tilde Q}(x)}_{L^p(\Rn)}^p}^{1/p}
}

Therefore it is enough to show that one has an expression of the form
\nm{eq:TLgoal}{
    \sum_{j=1}^\infty\| T_j R_{Q,j}\|_{L^p(\Rn)}^p\lesssim 1
}
uniformly in $ Q\in\mathcal D$, where \eq{
    R_{Q,j}(x):=\sum_{\tilde Q\in\mathcal D_j( Q)}b_{\tilde Q}\,\check\Psi^{\tilde Q}(x).
}
(Recall from Lemma \ref{lem:TLatomseq4} that $R_{Q,j}(x)=0$ for $j< k_{\tilde Q}$.)\\

To show \eqref{eq:TLgoal} we need to show the following three estimates:
\begin{align}
& \sum_{j=\max\{1,k_Q\}}^\infty\| T_j R_{Q,j}\|_{L^p(Q^*_\rho)}^p\lesssim 1,\label{eq:TLnewgoal1}\\
& \sum_{j=\max\{1,k_Q\}}^\infty\| T_j  \chi_{2\sqrt n Q} R_{Q,j}\|_{L^p(\Rn \setminus Q^*_\rho)}^p\lesssim 1,\label{eq:TLnewgoal2}\\
& \sum_{j=\max\{1,k_Q\}}^\infty\| T_j  \chi_{\Rn \setminus2\sqrt n Q}R_{Q,j}\|_{L^p(\Rn \setminus Q^*_\rho)}^p\lesssim 1,\label{eq:TLnewgoal3}
\end{align}

\makeatletter 
\renewcommand{\eqref}[1]{\tagform@{\ref{#1}}}
\def\maketag@@@#1{\hbox{#1}}
\textbf{Step 1 -- Proof of \eqref{eq:TLnewgoal1} when $\boldsymbol{k_Q\geq 1}$ and $\boldsymbol{0<\rho\leq 1}$.}
\makeatother

Using Lemma \ref{lem:TLatomseq4} and Theorem \ref{basicL2}, we have
\eq{
&\| T_j R_{Q,j}\|_{L^p(Q^*_\rho)} \leq |Q^*_\rho|^{1/p-1/2} \| T_j R_{Q,j}\|_{L^{2}(\Rn)} \\&\lesssim 2^{-k_{  Q} (n-\kappa)\rho(1/p-1/2)}2^{-j(\kappa+(n-\kappa)(1-\rho))(1/p-1/2)} \| R_{Q,j}\|_{L^{2}(\Rn)} \\&\lesssim 2^{k_{  Q}(\kappa+(n-\kappa)(1-\rho))(1/p-1/2)}2^{-j(\kappa+(n-\kappa)(1-\rho))(1/p-1/2)}
}
Hence LHS of \eqref{eq:TLnewgoal1} is bounded by
\eq{
\sum_{j=\max\{1,k_Q\}}^\infty 2^{k_{  Q}(\kappa+(n-\kappa)(1-\rho))(1-p/2)}2^{-j(\kappa+(n-\kappa)(1-\rho))(1-p/2)} \lesssim 1.
}
\hspace*{1cm}\\

\textbf{Step 2 -- Proof of \eqref{eq:TLnewgoal1} when $\boldsymbol{k_Q\leq 0}$ or $\boldsymbol{\rho=0}$.}\\
In this case $Q^*_\rho = 2\sqrt n Q$. We have, once again using Theorem \ref{basicL2},
\eq{
    &\| T_j R_{Q,j}\|_{L^p(Q^*_\rho)} \leq |  Q|^{1/p-1/2} \| T_j R_{Q,j}\|_{L^{2}(\Rn)} 
\\&\lesssim 
    2^{-k_{  Q} n(1/p-1/2)}2^{-j(\kappa+(n-\kappa)(1-\rho))(1/p-1/2)} \| R_{Q,j}\|_{L^{2}(\Rn)}
\lesssim 
    2^{-j(\kappa+(n-\kappa)(1-\rho))(1/p-1/2)}
}

Hence LHS of \eqref{eq:TLnewgoal1} is bounded by
\eq{
\sum_{j=\max\{1,k_Q\}}^\infty 2^{-j(\kappa+(n-\kappa)(1-\rho))(1/p-1/2)} \lesssim 1
}

\textbf{Step 3 -- Proof of \eqref{eq:TLnewgoal2} and \eqref{eq:TLnewgoal3}}\\
\eqref{eq:TLnewgoal2} follows directly from using Lemma \ref{lem:bstarcestimate} and Lemma \ref{lem:TLatomseq4}. \eqref{eq:TLnewgoal3} follows immediately from Lemma \ref{lem:TLatomseq5}. This concludes the proof.
\end{proof}

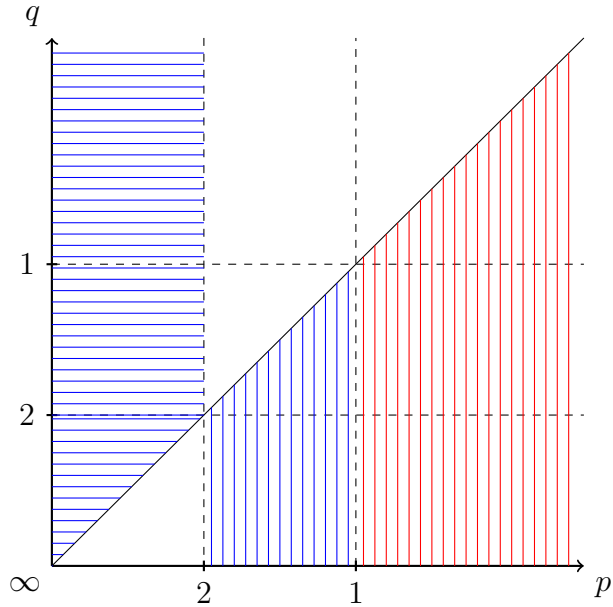
\begin{figure}\label{Firgure triebel lizorkin}
\begin{tikzpicture}

    \draw[thick,->] (0,0) -- (7,0) node[anchor=north west] {$p$};
    \draw[thick,->] (0,0) -- (0,7) node[anchor=south east] {$q$};

    \draw (0,0) -- (7,7);
    \draw[dashed] (2,7) -- (2,0);
    \draw[dashed] (0,2) -- (7,2);
    \draw[dashed] (4,0) -- (4,7);
    \draw[dashed] (0,4) -- (7,4);
    
    \draw[thick] (2 cm,2pt) -- (2 cm,-2pt) node[anchor=north] {$2$};
    \draw[thick] (4 cm,2pt) -- (4 cm,-2pt) node[anchor=north] {$1$};
    \draw[thick] (2pt,2 cm) -- (-2pt,2 cm) node[anchor=east] {$2$};
    \draw[thick] (2pt,4 cm) -- (-2pt,4 cm) node[anchor=east] {$1$};
    \draw (0,0) node[anchor=north east] {$\infty$};

    \foreach \x in {2.1,2.25,...,4}{
        \draw[blue] (\x cm,0) -- (\x cm,\x cm);
    }
    \foreach \x in {4.1,4.25,...,6.90}{
        \draw[red] (\x cm,0) -- (\x cm,\x cm);
    }
    \foreach \y in {0,0.15,...,2}{
        \draw[blue] (0,\y cm) -- (\y cm,\y cm);
    }
    \foreach \y in {2,2.15,...,6.8}{
        \draw[blue] (0,\y cm) -- (2 cm,\y cm);
    }

\end{tikzpicture}
\caption{Triebel-Lizorkin boundedness for FIO:s with their respective critical decay. The blue horizontal lines illustrate Theorem \ref{thm:PRS}. The blue vertical lines illustrate the boundedness results obtained by applying a duality argument for the case $p>2$. The red vertical lines illustrate Proposition \ref{prop:TLp<1} together with interpolation with the vertical blue area.}\label{pic:TLendpointresults}
\end{figure}

We now have all the tools we need to generalize \cite[Corollary 3.2]{PRS}.
\begin{Th}\label{thm:TLpgeq2FIO}
Let $1\leq \kappa \leq n-1,$ $0\leq\rho\leq1,$ $0\leq \delta<1$ and $a(x,\xi)\in S^{m}_{\rho,\delta}(\Rn)$ where $$ m=-\big(\kappa+(n-\kappa)(1-\rho)\big)\Big|\frac{1}{p}-\frac{1}{2}\Big| - n\max\set{0,\frac{\delta-\rho}2}.$$ 
Assume also that $\varphi(x,\xi) \in\Phi^2,$ is \emph{SND} on the support of $a(x,\xi)$ and is positively homogeneous of degree one in $\xi$ and has rank $\kappa$ on the support of $a(x, \xi)$. If $s\in \Rl$ and either one of the following cases hold 
\begin{enumerate}
\item[$(i)$] $2<p<\infty$ when $0<q\leq p$,
\item[$(ii)$] $\frac{n}{n+1}<p<2$ when $p\leq q$,
\item[$(iii)$] $p=q=2$,
\end{enumerate}
then it is true that $T_a^\varphi$ is bounded from $F_{p,q}^{s}(\Rn)$ to $F_{p,q}^s(\Rn)$.
\end{Th}

(These ranges of $p$ and $q$ values are illustrated in {\bf Figure \ref{pic:TLendpointresults}}.)

\begin{proof}
As before, we only treat the case when $0\leq \rho\leq \frac{2}{3}$, since the part when $\frac{2}{3}<\rho\leq 1$ can be done in a similar manner following the decomposition in \cite{SSS}.\\

We separate the operator into a low and a high-frequency part. The result for the low-frequency part follows from Theorem \ref{thm:low_freq_TL_BL_FIO}, so we only consider the high-frequency part from now on. \\

Observe also that the contents of $(iii)$ is contained in Theorem \ref{basicL2}. So from now on we only need to show $(i)$ and $(ii)$.\\

We split the proof into different ranges of $p$ and $q$, the two parts of the proof correspond to the blue and the red regions in {\bf Figure \ref{pic:TLendpointresults}}, respectively.\\

\textbf{Part 1 -- Proof when $\boldsymbol{p > 2}$ and $\boldsymbol{p\geq q>0}$}\\
We use Theorem \ref{thm:left composition with pseudo} to write 
\eq{
    \psi_j(D) T_a^\varphi f(x) = T_a^\varphi \psi_j(\nabla_x\varphi(x,D))f(x) + Rf(x),
}
The operator $R$ is bounded by Lemma \ref{lem:local and global nonendpoint TL}. So from now on we will only consider the first term. Denote $T_j:= T_a^\varphi \psi_j(\nabla_x\varphi(x,D)).$\\

Hence we can use Theorem \ref{thm:PRS} with $\mathcal E_{\tilde Q}:= \alpha_n Q^*_\rho$ (with an appropriately chosen $\alpha_n$ which yields $\Gamma\geq 1$ in that theorem), $b:=\kappa +(n-\kappa )(1-\rho)$ and $S_j:= T_j(1-\Delta)^{-b/2p}$ to prove the desired result. Observe that the $h^p\to L^p$ boundedness (Proposition \ref{prop:hptoLpboundednessFIO}) and the $L^2$ boundedness (Theorem \ref{basicL2}) of $T_a^\varphi$ yield \eqref{eq:PRS1} and \eqref{eq:PRS2} respectively.\\

We have that for $0<\rho\leq 1$
\eq{
    &\int_{\Rn \setminus Q^*_\rho} |K_j(x, y)| \dd y 
    \\
&\lesssim
    \sum_{\nu=1}^{\mathscr N_j}\norm{\frac1{g_j^\nu(x,\cdot)^{N}}}_{L^{2}(\Rn \setminus Q^*_\rho)}\| g_j^\nu(x,\cdot)^N \,K_j^\nu(x,\cdot)\|_{L^2(\Rn \setminus Q^*_\rho)} 
\\&\lesssim
    \sum_{\ell=0}^{6N} 2^{(j-k_Q)\ell/2}2^{-L(j-k_Q)/2}\,2^{-jn(1/p-1)}2^{k n(1/p-1)}\lesssim 2^{(j-k_Q)(3N-L/2)}\,2^{(k-j) n(1/p-1)}
}
where the first inequality comes from a second dyadic decomposition (as in section \ref{SSS decomposition}) and Hölder's inequality, while the second inequality follows by a similar argument as in \eqref{eq:equation6}.\\

Theorem \ref{thm:PRS} now yields that
\eq{
    &\norm{ \brkt{\sum_{j=0}^\infty 
    2^{jbq/p}|\Psi_j(D) T_jf_j|^q}^{1/q}}_{L^p(\Rn)}=\norm{ \brkt{\sum_{j=0}^\infty 
    2^{jbq/p}|\Psi_j(D) S_jf_j|^q}^{1/q}}_{L^p(\Rn)} 
    \\&
    \lesssim \brkt{ \sum_{j=0}^\infty \|f_j\|_{L^p(\Rn)}^p}^{1/p}.
}

Thus $S_j:B^0_{p,p}\to F^{b/p}_{p,q}$ which immediately implies that $T_j:B^{b/p}_{p,p}\to F^{b/p}_{p,q}$. Now the assertion follows from the facts that $F^0_{p,q}\xhookrightarrow{} F^0_{p,p}=B^0_{p,p}$ and the calculus using Bessel potentials and Theorem \ref{thm:left composition with pseudo}.\\

\textbf{Part 2 -- Proof when $\boldsymbol{1<p<2}$ and $\boldsymbol{p\leq q}$}\\
 Using \eqref{the_adjoint}, one can deduce an estimate analogue to \eqref{eq:kernelestimate1} for the kernel of the adjoint of $T_a^\varphi$, using which we can also obtain the estimate \eqref{eq:PRS} for the kernel of $(T_a^\varphi)^*$.\\
 
\textbf{Part 3 -- Proof when $\boldsymbol{0<p<1, p\leq q\leq\infty}$}\\
By Proposition \ref{prop:TLp<1} one obtains the result for $0<p<\frac{2n-2n\rho+2\kappa\rho}{2n-n\rho+\kappa\rho}$ and $q\geq p$.\\

\textbf{Part 3.1 -- Proof when $\boldsymbol{0<p\leq q<1}$}\\
In this case using Riesz-Thorin interpolation with $F^s_{p,p}=B^s_{p,p},$ yields the result. \\

\textbf{Part 3.2 -- Proof when $\boldsymbol{0<p<1\leq q\leq \infty}$}\\
Here Riesz-Thorin interpolation with {\bf{Part 2}} above, yields the result. \\

Now notice that $(1-\Delta)^{\frac s2} T^{\varphi}_a(1-\Delta)^{-\frac s2}$ is a similar operator associated to an amplitude in $S^{m_c(p)}(\Rn)$ and phase $\varphi$, and hence bounded from $F^0_{p,q}(\Rn)$ to itself. Therefore using the fact that the operator $(1-\Delta)^{\frac{s}{2}}$ is an isomorphism from $F^s_{p,q}(\Rn)$ to $F^0_{p,q}(\Rn)$ for $0<p\leq \infty$, we obtain the desired result. 
\end{proof}

\subsection{Triebel-Lizorkin estimates related to forbidden amplitudes}\label{subsec:FIO_Forbidden}

It is well known that the Fourier integral operators with amplitudes in $S^{m}_{1,1}(\Rl^n)$ fail to be $L^2$-bounded. However, one may show that these operators are Sobolev-bounded for $H^s(\Rn)$ for $s>0$ (see \cite{CIS}), and the pseudodifferential case goes back to E. Stein and independently by Y. Meyer. In this section, we state and prove two results about the boundedness of Fourier integral operators in Triebel-Lizorkin spaces with amplitudes in $S^{m}_{1,1}(\Rl^n)$, and also the much worse case of $S^{m}_{0,1}(\Rl^n)$, which will be the content of the following theorem.

\begin{Th}\label{thm:Sobolev_fio0}
Let $n\geq 1$ and
    assume that $r\in[1,2]$, $r<p< \infty$, $ r< q\leq \infty$, and $\varphi\in\Phi^2$. Then,  if  $a\in S^{-\frac{n}{r}}_{0,1}(\Rl^n)$, and $s>
    n(\frac{1}{\min\{1, p, q\}}-1)$ for $q<\infty$, and $s>\frac{n}{p}$ when $q=\infty$, and
the \emph{FIO} $T_a^\varphi$ is bounded from $F^s_{p,q}(\Rl^n)\to F^s_{p,q}(\Rl^n).$ \\
\end{Th}

\begin{proof}
Here we use the reduction of the phase to bring it to the form $x\cdot \xi+ \theta(x, \xi)$ where $\theta \in \Phi^1$. Then we note that for $a(x, \xi)\in S^{\frac{-n}{r}}_{0,1}(\Rn)$ one has that $ a(x, \xi) e^{i\theta(x,\xi)}\in S^{\frac{-n}{r}}_{0,1}(\Rn).$ Therefore we could reduce the study of the boundedness of $T^{\varphi}_a$ on $F^s_{p,q}(\Rl^n)$ to that of a pseudodifferential operator $a(x,D)$ with symbol $a(x,\xi)\in S^{\frac{-n}{r}}_{0,1}(\Rn).$

First, let us rewrite the operator as
\begin{equation*}
     T_a^\varphi f(x)= \int_{\Rn} a(x, \xi) \, e^{ix \cdot \xi}\, \widehat{f}(\xi)\ddd \xi.
\end{equation*}

Now, we decompose $a(x,\xi)$ as 
\begin{equation*}
    a(x,\xi)= \sum_{k, l} a_{k,l} (x,\xi)
\end{equation*}
where
\nm{eq1}{
    &a_{k,l} (x,\xi)= \int_{\Rn} e^{ix\cdot \eta}\,\psi_{k+l}(\eta)\,\psi_{k}(\xi)\,\int_{\Rn} e^{-iy\cdot\eta } \,a(y,\xi)\dd y\ddd \eta 
\\&
    = \psi_{k+l}^{\vee} \ast (a(\cdot, \xi) \psi_{k}(\xi)),\quad k\geq 0,l>0,\nonumber
}
and
\begin{equation*}
    a_{k,0} (x,\xi)= \int_{\Rn} e^{ix\cdot \eta}\big(\psi_{0}+\cdots+\psi_{k}\big)(\eta)\,\psi_{k}(\xi)\int_{\Rn} e^{-iy\cdot\eta }\,    a(y,\xi)\dd y\ddd \eta,\quad k\geq 0.
\end{equation*}
Notice that $a_{k,l} (x,\xi)$ has $\xi$-support in $|\xi|\sim 2^{l+k}$.\\

Now we claim that for all $\alpha$
\begin{equation}\label{hellish:claim1}
    |\d^\alpha_\xi    a_{k,l}(x,\xi)| \lesssim 2^{-lN} 2^{\frac{-kn}{r}},
\end{equation}
for $N>0$.\\

To see this we begin with the case when $l=0$. Using Young's convolution inequality one has that
\begin{align*}
    &\|\d^\alpha_\xi    a_{k,0}\|_{L^\infty(\Rn)}\\
    &=\|\d^\alpha_\xi (\FF\big(\psi_{0}+\cdots+\psi_{k})) * (\psi_{k}(\xi) \,a(\cdot,\xi))(x)\|_{L^\infty(\Rn)}\\
    &\lesssim \| \FF\big(\psi_{0}+\cdots+\psi_{k})\|_{L^1(\Rn)}\,\|\d^\alpha_\xi(\psi_{k}(\xi)\,a(x,\xi))\|_{L^\infty(\Rn)}\\
    &\lesssim \|\d^\alpha_\xi(\psi_{k}(\xi)\,a(x,\xi))\|_{L^\infty(\Rn)}
    \lesssim 2^{\frac{-kn}{r}}.
\end{align*}

Setting $\psi_k a:= a_k$, since $|\partial^{\alpha}_{\xi} \partial ^{\beta}_x a_k(x, \xi)|\lesssim 2^{k(\frac{-n}{r}+|\beta|)},$ then one has for any $M>0$ that
\begin{align*}
    &| \partial_{\xi}^{\alpha}   a_{k,l} (x, \xi)|
    =|(\FF\psi_{k+l}) *  \partial_{\xi}^{\alpha}    a_k(\cdot,\xi))(x)|\\&=
    \Big|\int_{\Rn} |\eta|^{-2M}\, \psi_{k+l}(\eta)\Big( \int_{\Rn} |\eta|^{2M} \partial_{\xi}^{\alpha}    a_k(y, \xi) \,e^{-iy\cdot \eta} \dd y\Big)\, e^{ix\cdot \eta}  \, \dd \eta\Big|
   \\&=  \Big|\int_{\Rn} |\eta|^{-2M} \psi_{k+l}(\eta)\Big( \int_{\Rn} \Delta_{y}^{M}\partial_{\xi}^{\alpha} a_k(y, \xi) \,e^{-iy\cdot \eta} \dd y\Big) \,e^{ix\cdot \eta}  \, \dd \eta \Big|\\&=
    2^{-2M(k+l)} 2^{(k+l)n} \Big| \int_{\Rn}  |\eta|^{-2M} \psi_{1}(\eta)\,\Big( \int_{\Rn} \Delta_{y}^{M}\partial_{\xi}^{\alpha}    a_k(y, \xi) \,e^{-i2^{k+l}y\cdot \eta} \dd y\Big) \,e^{i 2^{k+l}x\cdot \eta}  \, \dd \eta\Big|\\&=
     2^{-2M(k+l)} 2^{(k+l)n}\Big|\int_{\Rn} \theta^{\vee}(2^{k+l}y)\, \Delta_{y}^{M} \partial_{\xi}^{\alpha} a_k(x-y, \xi) \dd y\Big| \\&=
     2^{-2M(k+l)} 2^{(k+l)n}\, \sup_{y} |\Delta_{y}^{M}\partial_{\xi}^{\alpha}    a_k(y, \xi)| 
     \int_{\Rn} |\theta^{\vee}(2^{k+l}y)|  \dd y 
     \\&=
    2^{-2M(k+l)} 2^{(k+l)n} 2^{-(k+l)n} 2^{2Mk } 2^{\frac{-kn}{r}}\lesssim 2^{-2Ml} 2^{\frac{-kn}{r}},
\end{align*}

where $\theta(y):=  |y|^{-2M} \psi_{1}(y) \in \mathscr{S}(\Rn).$ Thus choosing $N:=2M$ proves \eqref{hellish:claim1}.\\

on the support of $a_{k,l}$.\\

Now let
\begin{equation*}
    K_{k,l}(x,y)=\frac{1}{(2\pi)^{n}}\int_{\Rn} a_{k,l}(x,\xi) \,e^{i\langle y,
    \xi\rangle}d\xi= \check{a}_{k,l}(x,y),
\end{equation*}
where $\check{a}_k$ here denotes the inverse Fourier transform
of $a_{k}(x,\xi)$ with respect to $\xi$. One observes that
\eq{
&|a_{k,l}(x,D)f_k(x)|^{r}=
\Big|\int_{\Rn} K_{k,l}(x,y)\,f_k(x-y)\dd y\Big|^{r}\\&=
\Big|\int{\Rn} K_{k,l}(x,y)\,\sigma(y)\frac{1}{\sigma(y)}\,f_k(x-y)\dd y\Big|^{r},
}
with weight functions $\sigma(y)$ which will be chosen
momentarily. Therefore, H\"older's inequality with $\frac{1}{r}+ \frac{1}{r'}=1$ yields
\begin{equation}\label{eq2.14}
    |a_{k,l}(x,D)f_k(x)|^{r}\leq
\Big( \int_{\Rn} |K_{k,l}(x,y)|^{r'} \,| \sigma(y)|^{r'} \dd y\Big)^{\frac{r}{r'}}
  \Big( \int_{\Rn} \frac{| f_k(x-y)|^{r}}{|\sigma(y)|^{r}}\dd y\Big).
\end{equation}
Now for an $\lambda>\frac{n}{2}$, we
define $\sigma$ by
\begin{equation*}
\sigma(y)=\begin{cases}
1, &| y| \leq 1; \\
| y|^{\lambda}, & | y| >1.
\end{cases}
\end{equation*}
By Hausdorff-Young's theorem and the estimate \eqref{hellish:claim1}, first
for $\alpha=0$ and then for $| \alpha|=\lambda$, we have
\begin{equation*}\label{eq2.17}
\int_{\Rn} |K_{k,l} (x,y)|^{r'} \dd y \leq
\int_{\Rn} |a_{k,l} (x,\xi)|^{r} d\xi
\lesssim \int_{| \xi| \sim 2^{k+l}} 2^{-rlN}2^{-nk}d\xi
\lesssim 2^{l(n-rN)},
\end{equation*}
and
\begin{equation*}\label{eq2.18}
\int_{\Rn} | K_{k,l} (x,y)|^{r'}\, |y|^{r'\lambda}\dd y
\lesssim \int_{\Rn} |\nabla_{\xi}^{\lambda} a_{k,l} (x,\xi)|^r \,d\xi
\lesssim  \int\limits_{|\xi| \sim 2^{k+l}}2^{-rNl}2^{-kn}\,d\xi
\lesssim  2^{l(\frac{n}{2}-rN)} .
\end{equation*}
Hence, splitting the integral into $| y| \leq 1$
and $,| y| > 1$ yields
\begin{equation*}
\int_{\Rn} | K_{k,l} (x,y)|^{r}\,| \sigma(y)|^{r}\dd y
\lesssim 2^{l(n-rN)}.
\end{equation*}

Furthermore, we also have
\begin{equation*}
\int_{\Rn} \frac{| f_k(x-y)|^{r}\dd y}{|\sigma(y)|^{r}}\lesssim \big(\mathcal{M}_{r} f_k(x)\big)^{r}
\end{equation*}
with a constant that only depends on the dimension $n$. Thus
\eqref{eq2.14} yields
\begin{equation}\label{pointwiseakestim}
 |a_{k,l}(x,D)f_k(x)|^{r} \lesssim 2^{l(n-rN)}\big(\mathcal{M}_{r} f_k(x)\big)^{r}
\end{equation}

Thus we obtain the pointwise estimate
\begin{equation*}
|a_{k,l}(x,D)f_k(x)|
\lesssim 2^{l(n-rN)/r}\mathcal{M}_{r} f_k(x)
\end{equation*}

Now recall that at the beginning of the proof we reduced the operator $T^{\varphi}_a$ to $a(\cdot,D)$, note that the same argument allows us to recover the original operator since upon examination of that argument one sees that 
\begin{equation*}
    T^{\varphi}_a f_k(x)=a(x,D)f_k(x)
\end{equation*}
and therefore 
\begin{equation}\label{interpolationstep1}
|T^{\varphi}_{a_{k,l}} f_k(x)|=|a_{k,l}(x,D)f_k(x)|
\lesssim 2^{l(n-rN)/r}\mathcal{M}_{r} f_k(x).
\end{equation}

Now taking $h_{k,l}(x)=T^{\varphi}_{a_{k,l}} f_k(x)$ and $g_l(x)=\sum_{k\geq 0}h_{k,l}(x)$ and $u_k=f_k$ one has by the definition of the Triebel-Lizorkin norm (Definition \ref{def:TLspace}) that
\begin{equation}\label{hypothessis for 4.8 1}
    \norm{ \{2^{ks}u_k\}_{k=0}^\infty}_{_{L^p(l^q)}} \lesssim \Vert f\Vert_{F^{s}_{p,q}(\Rl^n)}.
\end{equation}

Thus for $0<p,q<\infty$ and $s>n(\frac{1}{\min\{p, q,1\}}-1)$ one obtains using \eqref{hypothessis for 4.8 1} and \eqref{interpolationstep1} and Lemma \ref{Usingequivalenceofnormslemma} that
\begin{equation}
    \|T^{\varphi}_{a}f\|_{F^{s}_{p,q}(\Rn)}=\norm{\sum_{l\geq 0} g_l}_{_{L^p(l^q)}}\lesssim \Vert f\Vert_{F^{s}_{p,q}(\Rn)}.
\end{equation}
\end{proof}
The following result is the generalisation of the boundedness of FIOs with forbidden amplitudes on Sobolev spaces of positive order obtained in \cite{CIS}, which was in turn a generalisation on an unpublished result of Meyer and Stein on boundedness of pseudodifferential operators with forbidden symbols  on the aforementioned Sobolev spaces.
\begin{Th}\label{thm:Sobolev_fio1}
Let $n\geq 1$ and let $m=-\kappa\big|\frac{1}{p}-\frac{1}{2}\big|$ and $a\in S^m_{1,1}(\Rl^n)$. Assume that $0\leq \kappa\leq n-1$ and $\varphi\in\Phi^2$ is \emph{SND} and has rank $\kappa$ on the support of $a(x, \xi)$. If $s>n\big(\frac{1}{\min\{1, p, q\}}-1\big)$, $0<\kappa\leq n-1$
    and either one of the following cases hold 
    \begin{enumerate}
        \item[$(i)$] $2<p<\infty$ when $0<q\leq p,$\\
        \item[$(ii)$] $\frac{n}{n+1}<p<2$ when $p\leq q,$\\
        \item[$(iii)$] $p=q=2,$
    \end{enumerate}
    or if $s>\frac{n}{p}$ with $q=\infty$, then the \emph{FIO} $T_a^\varphi$ is bounded from $F^s_{p,q}(\Rl^n)\to F^s_{p,q}(\Rl^n).$ \\
    
    In the case of $\kappa=0$ $($to which pseudodifferential operators belong$)$, one can extend the boundedness results to the ranges $0<p<\infty$ and $0<q\leq \infty$ and the aforementioned ranges of $s$.
   
\end{Th}

\begin{proof}

Observe that $\mathbf{t}_x(\xi):=\nabla_x\varphi(x,\xi)$ is a global diffeomorphism.  This allows us to rewrite the operator as
\begin{equation*}
    T_a^\varphi f(x)= \int_{\Rl^n} e^{i\varphi(x,\xi)}\, \sigma (x,\mathbf{t}_x(\xi))\jap{\xi}^{m} \widehat{f} (\xi) \ddd\xi,
\end{equation*}
where $\sigma(x,\xi):=a(x,\mathbf{t}_{x}^{-1}(\xi))\jap{\mathbf{t}_{x}^{-1}(\xi)}^{-m} \in S^0_{1,1}(\Rl ^n)$ (which could be checked by using the facts that $a(x, \xi)\in S^m_{1,1}(\Rl ^n)$, $\varphi\in \Phi^2$ and $|\mathbf{t}_{x}(\xi)|\sim |\xi| $).\\

Now, observe that $S_{1, 1}^{0}(\Rl^n)\subset C_{*}^{r} S_{1, 1}^{0}(\Rl^n)$ for all $r>0$. Therefore, as was done in \cite{Bourdaud, Meyer} for the case of pseudodifferential operators with forbidden symbols, it is enough to show the result for elementary amplitudes in the class $C_{*}^{r} S_{1, 1}^{0}(\Rl^n)$ where $r>0$.\\

By definition, an elementary amplitude $\sigma(x,\xi)\in C_{*}^{r} S_{1, 1}^{0}(\Rl^n)$ is of the form
\eq{\sum_{k\geq 0} \alpha_k(x)\,\psi_k(\xi)
}
where $\psi_{k}$ was introduced in Definition \ref{def:LP} and $\alpha_{k}(x)$ satisfies
\begin{align}\label{estimates for Qk}
    &|\partial^\gamma\alpha_{k}(x)|\lesssim 2^{k|\gamma|} & \|\alpha_{k}(x)\|_{C_{*}^{r}(\Rn)}\lesssim 2^{kr}
\end{align}
for $|\gamma|\geq 0$ and the $C_{*}^{r}$-norm is given in Definition \ref{def:Zygmund}.\\

We define
\eq{
    \alpha_{k,l}(x) := \psi_{k+l}(D)\alpha_k(x)
}
(note that $\alpha_k$ is not a Littlewood-Paley piece). Then
\begin{equation*}
    \sigma(x,\mathbf{t}_{x}(\xi)) = \sum_{k\geq 0} \alpha_k(x)\,\psi_k(\nabla_x\varphi(x,\xi))= \sum_{k\geq 0} \sum_{l\geq 0}\alpha_{k,l}(x)\,\psi_k(\nabla_x\varphi(x,\xi)).
\end{equation*}

We claim that for $k\geq 1$ and $r>0$
\begin{equation}\label{third estimate for Qk}
\Vert \alpha_{k,l}\Vert_{L^\infty(\Rl^n)}\lesssim 2^{-rl}.
\end{equation}

To see \eqref{third estimate for Qk}, we observe that the second part of \eqref{estimates for Qk}, relations \eqref{Zygmundegenskap} and the fact that 
$$
    \int_{\Rl^n} 2^{(k+l)n}\, \psi^\vee\left(\frac{y}{2^{-(k+l)}}\right)  \dd y= \int_{\Rl^n} \psi^\vee\left(y\right) \dd y = 0,
$$ 
yield that \eqref{third estimate for Qk} follows from
\begin{align*}
    |\alpha_{k,l}(x)| 
    & \lesssim\Big|\int_{\Rl^n} 2^{(k+l)n}\, \psi^\vee\Big(\frac{y}{2^{-(k+l)}}\Big)\,
    \big(\alpha_{k}(x-y)-\alpha_{k} (x)\big) \dd y\Big| \\ 
    & \lesssim {2^{-l r}} \int_{\Rl^n}\Big|2^{(k+l)n} \psi^\vee \Big(\frac{y}{2^{-(k+l)}}\Big)\Big|\, \frac{|y|^r}{2^{-(k+l)r}} \dd y \lesssim 2^{-l r}, 
\end{align*}
for $k>0$ and $l\geq 0$. For $k=l=0$ this is a consequence of the $L^\infty$ boundedness of $\psi_0(D)$ and the first estimate in \eqref{estimates for Qk}.\\

If we now set
\begin{equation*}
    F_{k}(x):= \int_{\Rl^n} \psi_k(\nabla_x\varphi(x,\xi))\jap{\xi}^m\, e^{i\varphi(x,\xi)} \,\widehat{f}(\xi) \ddd \xi,
\end{equation*}
then
\begin{equation*}
    T^{\varphi}_{a} f(x) = \sum_{\substack{k\geq0\\l\geq 0}} \alpha_{k,l}(x)\,  F_{k} (x).
\end{equation*}

For the low frequency part when $k=0$ we just use Theorem \ref{thm:low_freq_TL_BL_FIO}, hence from now on we only consider the case $k>0$.\\

We claim that
\eq{
     |\alpha_{k,l}(x)\,  F_{k} (x)| \lesssim 2^{-rl}\mathcal{M}f_k
}

At this point, using Theorem \ref{thm:left composition with pseudo} we obtain the representation
\begin{align*}
F_{k}(x)
&\nonumber = \psi_k(D)T^\varphi_{\jap{\xi}^{m}} f(x) 
- 2^{-k\eps }\,T^\varphi_{r_{k}} f(x)
\end{align*}
with
\begin{equation*}
|\partial^{\beta}_{\xi} \partial^{\gamma}_{x}  r_{k} (x,\xi)| 
\lesssim  \,\bra{\xi}^{m-(1/2-\eps)-|\beta|},
 \end{equation*}
where the estimates above are uniform in $k$.  \\

Define also,
\begin{align}
\label{F1 och F2}
    &F^1_{k}(x)=\psi_k(D)F(x)\\
    &F^2_{k}(x)=
    -2^{-k\eps}\,T^\varphi_{r_{k}} f(x).
\end{align}

We start with $F^1_{k}$. \\

Now observe that if $\chi$ is a smooth cut-off function that is 1 on the support of $\psi$, and that $F_k^{1}(x)=\Psi_k(D)F_k^{1}(x)$ where $\Psi_k(D)$ is a fattened Littlewood-Paley operator, then we have that
\begin{align*}
    &|  F^{1}_{k} (x)| = |\chi_k(D)  F^1_{k} (x)|\\
    &\nonumber= \Big|\int_{\Rl^n}\int_{\Rl^n} e^{iy\cdot\xi} \chi_k(\xi) \ddd \xi \,F^{1}_{k}(x-y) \dd y \Big|\\
    &\nonumber= \Big|\int_{\Rl^n}\check{\chi}_{k}(y) \,F^{1}_{k}(x-y) \dd y \Big|\\
    &\nonumber= \Big|\int_{\Rl^n}\check{\chi}_{k}(y)\,(1+2^k|y|)^a \,\frac{\Psi_k(D)F^{1}_{k}(x-y)}{(1+2^k|y|)^a} \dd y \Big|\\
    &\nonumber\lesssim \int_{\Rl^n}|\check{\chi}_{k}(y)|\,(1+2^k|y|)^a \dd y \,\,  \mathfrak{M}_{a,2^k}(F^{1}_{k})(x)\\
    &\nonumber\lesssim \int_{\Rl^n}|\check{\chi}_{k}(y)\,|(1+2^k|y|)^a \dd y \,\,  \mathcal{M}_{p_0}(F^{1}_{k})(x),
\end{align*}
for $a\geq \frac{n}{p_0}$. At this point, we claim that for any $a>0$ we have
\begin{equation*}
    \int_{\Rl^n}|\check{\chi}_{k}(y)|\,(1+2^k|y|)^{a}\dd y \lesssim 1,
\end{equation*}
where the implicit constant is independent of $k$.\\

To see this, denote the integral part of a real number $a$ by $[a]$ and its fractional part by $\{a\}$. Then note that for $b>0$ such that $-b+ \{a+b\}<-n$ and $[a+b]$ is an even natural number, we have 
\begin{align*}
    &\int_{\Rl^n}|\check{\chi}_{k}(y)|\,(1+2^k|y|)^{a}\dd y\\
    &=\int_{\Rl^n} 2^{-kn}|\check{\chi}_{k}(2^{-k}z)|\,(1+|z|)^{a}\dd z\\
    &=\int_{\Rl^n} |\mathscr{F}_{\xi\to z}\chi_k(2^k\cdot)|\,(1+ |z|)^{a+b}\,(1+ |z|)^{-b}\dd z\\
    &\lesssim \int_{\Rl^n} |\mathscr{F}_{\xi\to z}\big\{(1-\Delta)^{[a+b]/2}\chi_k(2^k\cdot)\big\}|\,(1+ |z|)^{-b+\{a+b\}}\dd z\\
    &\leq\|(1-\Delta)^{[a+b]/2}\chi_k(2^k\cdot)\|_{L^1(\Rn)}\int_{\Rl^n}(1+ |z|)^{-b+\{a+b\}}\dd z\\
    &\lesssim  1
\end{align*}
where we have used that $ \chi_k(2^kz)$ is in the Schwartz class and has fixed compact support in $z$.\\

This shows that 
\begin{align}\label{interpolationstep2}
    &| \alpha_{k,l}(x) F^{1}_{k} (x)| \lesssim 2^{-rl} \mathcal{M}_{p_0}F^{1}_{k}(x),
\end{align}
for all $r,p_0,k>0$ and $l\geq 0$.\\

Now taking $h_{k,l}(x)=\alpha_{k,l}(x) F^{1}_{k} (x)$ and $g_l(x)=\sum_{k\geq 0}h_{k,l}(x)$ and $u_k=F^{1}_{k}(x)$ one has by the definition of the Triebel-Lizorkin norm (Definition \ref{def:TLspace}) and Theorem \ref{thm:TLpgeq2FIO} that

\begin{equation}\label{hypothessis for 4.8 2}
    \| \{2^{ks}u_k\}_{k=0}^\infty\|_{_{L^p(l^q)}}\lesssim \Vert T^\varphi_{\jap{\xi}^m} f\Vert_{F^{s}_{p,q}(\Rl^n)} \lesssim \Vert f\Vert_{F^{s}_{p,q}(\Rl^n)},
\end{equation}
for the ranges of $p$'s and $q$'s in that theorem.\\
However if $\kappa=0$ then using the boundedness of pseudodifferential operators with symbols in $S^m_{1,0}(\Rl^n)$ established in \cite{Paivarinta:pdo-TL} one can obtain the boundedness above for $0<p<\infty$ and $0<q\leq \infty.$\\
Thus one obtains using \eqref{hypothessis for 4.8 2} and Lemma \ref{Usingequivalenceofnormslemma} that
\begin{equation}
    \norm{\sum_{l\geq 0} g_l}_{_{L^p(l^q)}}\lesssim \Vert f\Vert_{F^{s}_{p,q}(\Rn)}
\end{equation}

Proceeding to  the analysis of the $F_k^{2}$-term, we observe that it is an operator with amplitude in $S^{m-(1/2-\eps)}_{1,1}(\Rn)$. Now iterating the above process for $F_k^{2}$, we would eventually obtain an operator $T_{b_k}^\varphi$ with $b_k\in S^{-N}_{1,1}(\Rn)$ for $N$ arbitrarily large, and hence conclude that 
$$T_{b_k}^\varphi: F^s_{p,q}\to F^s_{p,q}$$
by following the proof of \cite[Theorem 6.2 (ii)]{IRS}, which only uses finitely many derivatives of the amplitude. Note also that the negative power of $2^k$ in the representation \eqref{F1 och F2} is also helpful in producing a convergent sum.
\end{proof}

\subsection{Besov-Lipschitz estimates}
In this section we include the Besov-Lipschitz boundedness results of FIO's with both classical and exotic amplitudes. We also present two results about the Besov-Lipschitz boundedness of FIOs with amplitudes in $S^m_{0,1}(\Rn)$ and $S^m_{1,1}(\Rn)$.\\

We also include an abstract theorem that allows one to lift $h^p\to L^p$ boundedness to Besov-Lipschitz boundedness. This result can be applied to a wide class of oscillatory integral operators.

\begin{Th}\label{theorem:BL-booster theorem}
Let $0<p,q<\infty$, $m \leq 0$, $a\in S^m_{\rho,\delta}(\Rn)$ such that $a$ is supported in $\Rn\times\Rn\setminus B(0,1)$. Suppose $\varphi$ is a phase function that verifies the conditions \eqref{eq:composition_conditions} of \emph{Theorem \ref{thm:left composition with pseudo}}.
If $T^\varphi_a$ is $h^p\to L^p$-bounded, then $T^\varphi_a$ is bounded from $B^s_{p,q}$ to itself for all $s\in \Rl$. 
\end{Th}

\begin{proof}
    For a proof see \cite{IRS}.
\end{proof}

\begin{Th}\label{thm:BL_tdphase_FIO}
Let $\rho\in[0,1],\,\delta\in[0,1)$ and $0<q\leq\infty$.
Assume furthermore that $\varphi(x, \xi)\in\Phi^2$ is an \emph{SND} phase of rank $\kappa$ satisfying the conditions in \emph{Definition \ref{def:FIO}} and $a\in S^{m_c(p)}_{\rho, \delta}(\Rl^n)$, and let $T_a^\varphi$  be the associated \emph{FIO}. Then $T_a^\varphi$ is a bounded operator from $B^s_{p,q}(\Rl^n)$ to $B^s_{p,q}(\Rl^n)$ for all $s\in\Rl$ for $\frac{n}{n+1}<p<\infty$. 
\end{Th}

\begin{proof}
    Let $\chi\in C_c^\infty(\Rn)$ be supported in $\{|\xi|\lesssim 1\}$, and write
    \begin{align*}
        T_a^\varphi f(x) &: =
        \int_{\Rn} e^{i\varphi(x,\xi)}\, a(x,\xi) \widehat{f}(\xi)\, (1-\chi(\xi)) \dd\xi 
        +
        \int_{\Rn} e^{i\varphi(x,\xi)}\, a(x,\xi) \widehat{f}(\xi)\, \chi(\xi) \dd\xi\\
        &=T_{\text{high}} f(x)+T_{\text{low}} f(x).
    \end{align*}
    The boundedness of $T_{\text{low}}$ follows from the low frequency result Theorem \ref{thm:low_freq_TL_BL_FIO}. The boundedness of $T_{\text{high}}$ follows immediately by the Besov-Lipschitz lift Theorem \ref{theorem:BL-booster theorem} and Proposition \ref{prop:hptoLpboundednessFIO}.
\end{proof}

The following two theorems are concerned with forbidden amplitudes. However, the proofs of these results are very close to the ones given in section \ref{subsec:FIO_Forbidden}. Therefore we confine ourselves to the discussion of the main points of difference between the proofs.

\begin{Th}\label{thm:BS_fio0}
Let $n\geq 1$ and
    assume that $r\in[1,2]$, $r<p< \infty$, $ r< q\leq \infty$, and $\varphi\in\Phi^2$. Then,  if  $a\in S^{-\frac{n}{r}}_{0,1}(\Rl^n)$, and $s>
    n(\frac{1}{\min\{1, p, q\}}-1)$ for $q<\infty$, and $s>\frac{n}{p}$ when $q=\infty$, and
the \emph{FIO} $T_a^\varphi$ is bounded from $B^s_{p,q}(\Rl^n)\to B^s_{p,q}(\Rl^n).$ 
\end{Th}

\begin{proof}
    The proof of this result differs from that of Theorem \ref{thm:Sobolev_fio0} mainly in showing the Besov-Lipschitz analogue of \eqref{hypothessis for 4.8 1}. For that, one just uses the Besov-Lipschitz part of Lemma \ref{Usingequivalenceofnormslemma}.
\end{proof}

\begin{Th}\label{thm:BS_fio1}
Let $n\geq 1$ and
    assume that $0\leq \kappa\leq n-1$ and $\frac{n}{n+1}<p< \infty, \,  0< q\leq \infty$, $\varphi\in\Phi^2$ is \emph{SND} and $\mathrm{rank}\,\partial^{2}_{\xi\xi} \varphi(x, \xi) = \kappa$ on the support of $a(x, \xi)$. Then, $m=-\kappa\big|\frac{1}{p}-\frac{1}{2}\big|$, $s>n\big(\frac{1}{\min\{1, p, q\}}-1\big)$ for $q<\infty$, and $s>\frac{n}{p}$ when $q=\infty$, and $a\in S^m_{1,1}(\Rl^n)$,
the \emph{FIO} $T_a^\varphi$ is bounded from $B^s_{p,q}(\Rl^n)\to B^s_{p,q}(\Rl^n).$ 
\end{Th}

\begin{proof}
    This proof differs from the proof of Theorem \ref{thm:Sobolev_fio1} mainly in showing the Besov-Lipschitz analogue of \eqref{hypothessis for 4.8 2}. To that end one just uses the Besov-Lipschitz part of Lemma \ref{Usingequivalenceofnormslemma} and thereafter Theorem \ref{thm:BL_tdphase_FIO} instead of Theorem \ref{thm:TLpgeq2FIO}. The rest of the proof remains the same.
\end{proof}

\begin{Rem}
   
As a matter of fact, see e.g., \cite{Triebel1}, if $ -\infty < s_1 < s < s_2< \infty, $ $0< p <\infty$ and $0< q \leq \infty$, \(0<q_{1} \leq \infty \) and \(0<q_{2} \leq \infty\) and $s=\theta s_1 + (1-\theta) s_2$, then $(F^{s_1}_{p,q_1}(\Rn), F^{s_2}_{p,q_2}(\Rn))_{\theta, q} = B^{s}_{p, q}(\Rn),$ in the sense of real interpolation. Moreover, the Besov spaces are independent of \(q_{1}\) and \(q_{2}\) as well as of the choice of \(s_{1}\) and \(s_{2}\). This fact can be used to deduce Besov-Lipschitz boundedness results from Triebel-Lizorkin ones.
\end{Rem}

\section{Estimates for oscillatory integrals related to the Klein-Gordon equation}\label{Klein-Gordon}

In this section we shall consider boundedness of oscillatory integral operators with phase functions of the form $\varphi(x, \xi)= \varphi_0 (x, \xi) + \varphi_1(x, \xi)$ where $\varphi_0(x,\xi) \in \Phi^2$ and $\varphi_1 \in S^{-1}_{1,0}(\Rl^n),$ for large $|\xi|.$\\

A basic example of such phase functions is $x\cdot \xi+ \jap{\xi}$ which is the one appearing in the study of the constant coefficient Klein-Gordon equation. Indeed, it is easy to see that for $|\xi|\geq 1$, this phase function can be written as $\varphi_0(x,\xi)+ \varphi_1(x, \xi) $ with $\varphi_{0}(x,\xi)=x\cdot \xi+|\xi| \in \Phi^2$ and $\varphi_1(x,\xi)= \lambda(\xi)$ with $\lambda \in S^{-1}_{1,0}(\Rl^n).$\\

We begin with the Triebel-Lizorkin estimates for the Klein-Gordon equation with Forbidden amplitudes in either $S^{m}_{1,1}$ or $S^{m}_{0,1}$. The following results contain the Klein-Gordon equation's version of Theorem \ref{thm:Sobolev_fio0} and Theorem \ref{thm:Sobolev_fio1}.

\begin{Prop}
Let $\rho$ be either $1$ or $0$. The result of \emph{Theorem \ref{thm:Sobolev_fio0}} and \emph{Theorem \ref{thm:Sobolev_fio1}} are also valid (with their respective values of $p,q,s$, and $m$) for the inhomogeneous phase function $\varphi(x, \xi)=\varphi_0 (x,\xi)+\varphi_1(x,\xi).$
\end{Prop}

\begin{proof}
This is clear since $e^{i\varphi_1(x, \xi)} a(x, \xi)\in S^{m}_{\rho,1}(\Rn)$ if $a\in S^m_{\rho,1}(\mathbb{R}^n).$
\end{proof}

We now proceed to estimates related to the Klein-Gordon equation with classical amplitudes.

\begin{Th}\label{linearhpthm}
Assume that either one of the cases hold true
\begin{enumerate}
    \item[$(i)$] $2<p<\infty$ when $0<q\leq p,$\\
    \item[$(ii)$] $\frac{n}{n+1}<p<2$ when $p\leq q,$\\
    \item[$(iii)$] $p=q=2,$
\end{enumerate}

Assume further that $m=-(n-\rho)\abs{\frac{1}{p}-\frac{1}{2}}-n\max\{0,(\delta-\rho)/2\}$ and $0\leq\rho\leq 1$ and $0\leq\delta<1$. Then any oscillatory integral operators of the form 
$$T_a f(x)= \int_{\Rl^n} a(x,\xi)\, e^{ix\cdot \xi +i\langle \xi\rangle}\, \widehat{f}(\xi) \ddd\xi ,$$
with an amplitude $a(x,\xi)\in S^{m}_{\rho,\delta}$, satisfies the estimate
\[
	\Vert T_a f\Vert_{F^{s}_{p,q}(\Rl^n)}\lesssim \|f\|_{F^s_{p,q}(\Rl^n)}.
\]

\end{Th}
 \begin{proof}
First, we divide the operator into low and high-frequency portions, $T_{a_L}$ and $T_{a_H}$ with an amplitude that is compactly supported in $\xi$, and one whose frequency support is in the set $|\xi|\geq 1$, respectively. For $T_{a_L}$ the result follows from \ref{thm:low_freq_TL_BL_FIO}. For $T_{a_H}$ we note that the phase function $\langle \xi\rangle $ can be written (for $|\xi|\geq 1$) as $|\xi|+ \lambda(\xi)$ with $\lambda \in S_{1,0}^{-1}(\Rl^n).$ \\

Now since $b(x, \xi):=e^{i\lambda(\xi)}\, a(x, \xi) \in S_{\rho,\delta}^{m}(\Rl^n)$ and the rank of $x\cdot\xi+|\xi|$ is equal to $n-1$, Theorem \ref{thm:TLpgeq2FIO} yields that  the FIO
$$ T_{a_H}f(x)= \int_{\Rn} b(x, \xi)\, e^{ix\cdot \xi + |\xi|} \,\widehat{f}(\xi)  \ddd \xi, $$ is  bounded.
This ends the proof of the theorem.
\end{proof}

 \bibliographystyle{siam}
 \bibliography{references}

\end{document}